\documentclass[11pt,a4paper,reqno]{amsart} 

\usepackage[applemac]{inputenc}
\usepackage{amsmath}
\usepackage{amssymb}
\usepackage{euscript}
\usepackage{mathrsfs}  
\usepackage{multirow}
\usepackage{graphicx}
\usepackage{comment}
\usepackage{xcolor}
\usepackage[colorlinks=true,linktocpage=true,pagebackref=false, citecolor=black,linkcolor=black]{hyperref}
\usepackage{pifont}
\usepackage{tikz,pgfplots}
\usepackage{mathabx}
\usepackage[all]{xy}
\usepackage{enumitem}

\pgfplotsset{compat=1.10}
\usepgfplotslibrary{fillbetween}
\usetikzlibrary{cd,arrows,decorations.pathmorphing,backgrounds,automata,positioning,fit,matrix}

\pgfplotsset{soldot/.style={color=black,only marks,mark=*}} \pgfplotsset{holdot/.style={color=black,fill=white,only marks,mark=*}}

\newtheorem{thm}{Theorem}[section]

\newtheorem{lem}[thm]{Lemma}
\newtheorem{prop}[thm]{Proposition}

\newtheorem{cor}[thm]{Corollary}

\theoremstyle{definition}
\newtheorem{defn}[thm]{Definition}

\newtheorem{notation}[thm]{Notation}
\newtheorem*{ack}{Acknowledgements}

\theoremstyle{remark}
\newtheorem{remark}[thm]{Remark}
\newtheorem{remarks}[thm]{Remarks}

\newtheorem{examples}[thm]{Examples}

\numberwithin{equation}{section}
\numberwithin{figure}{section}

 \newcommand{\N}{{\mathbb N}}
\newcommand{\Z}{{\mathbb Z}} \newcommand{\R}{{\mathbb R}}
\newcommand{\Q}{{\mathbb Q}} \newcommand{\C}{{\mathbb C}}
 
\newcommand{\HH}{{\mathbb H}} 
\newcommand{\sph}{{\mathbb S}} 
\newcommand{\E}{{\mathbb E}} \newcommand{\PP}{{\mathbb P}}
\newcommand{\B}{{\mathbb B}}

\newcommand{\reg}{{\mathcal R}} 
 
 \newcommand{\I}{{\mathcal I}}
\newcommand{\G}{\mathbb{G}}

 \newcommand{\gtn}{{\mathfrak n}}

\newcommand{\Ll}{{\EuScript L}}

\newcommand{\Reg}{\operatorname{Reg}}
\newcommand{\Sing}{\operatorname{Sing}}

\newcommand{\cl}{\operatorname{Cl}}

\newcommand{\id}{\operatorname{id}}

\newcommand{\x}{{\tt x}}
\newcommand{\y}{{\tt y}}

\newcommand{\veps}{\varepsilon}

\newcommand{\ol}{\overline}

\newcommand{\qr}{{\ol{\Q}^r}}

\newcommand{\mr}{\mathrm}

\newcommand{\sss}{\scriptscriptstyle}
\newcommand{\II}{{\mc{I}}}
\newcommand{\ZZ}{{\mc{Z}}}
\newcommand{\mc}{\mathcal}

\newcommand{\mk}{\mathfrak}
\newcommand{\sfh}{\mathsf{h}}

\newcommand{\Nn}{\EuScript{N}}

\newcommand{\cc}{\mathtt{c}}
\newcommand{\cinfty}{{\EuScript{C}^\infty}}
\newcommand{\czero}{{\EuScript{C}^0}}
\newcommand{\cnu}{{\EuScript{C}^\nu}}

\makeatletter
\def\@tocline#1#2#3#4#5#6#7{\relax
  \ifnum #1>\c@tocdepth 
  \else
    \par \addpenalty\@secpenalty\addvspace{#2}%
    \begingroup \hyphenpenalty\@M
    \@ifempty{#4}{%
      \@tempdima\csname r@tocindent\number#1\endcsname\relax
    }{%
      \@tempdima#4\relax
    }%
    \parindent\z@ \leftskip#3\relax \advance\leftskip\@tempdima\relax
    \rightskip\@pnumwidth plus4em \parfillskip-\@pnumwidth
    #5\leavevmode\hskip-\@tempdima
      \ifcase #1
       \or\or \hskip 1em \or \hskip 2em \else \hskip 3em \fi%
      #6\nobreak\relax
    \dotfill\hbox to\@pnumwidth{\@tocpagenum{#7}}\par
    \nobreak
    \endgroup
  \fi}
\makeatother

\begin{document}

$\;$

\vspace{-2em}

\title[The Nash-Tognoli theorem over the rationals]{The Nash-Tognoli theorem over the rationals \\ and its version for isolated singularities}

\author{Riccardo Ghiloni}
\address{Dipartimento di Matematica, Universit\`a di Trento, Via Sommarive 14, 38123 Povo-Trento (ITALY)}
\email{riccardo.ghiloni@unitn.it}

\author{Enrico Savi}
\address{LAREMA of Universit\'e d'Angers, Angers (FRANCE)}
\email{enrico.savi@univ-angers.fr}


\begin{abstract}
Let $\Q$ be the field of rational numbers and let $X$ be a subset of $\R^n$. We say that $X$ is $\Q$-algebraic if it is the common zero set  in $\R^n$ of a family of polynomials in $\Q[\x_1,\ldots,\x_n]$. If $X$ is $\Q$-algebraic and of dimension $d$, then we say that $X$ is $\Q$-nonsingular if, for all $a\in X$, there exist a neighborhood $U$ of $a$ in $\R^n$ and $f_1,\ldots,f_{n-d}\in\Q[\x_1,\ldots,\x_n]$ such that $\nabla f_1(a),\ldots,\nabla f_{n-d}(a)$ are linearly independent and $X\cap U=\{x\in U:f_1(x)=0,\cdots,f_{n-d}(x)=0\}$. 

The celebrated Nash-Tognoli theorem asserts the following: if $M$ is a compact smooth mani\-fold of dimension $d$ and $\psi:M\to\R^{2d+1}$ is a smooth embedding, then $\psi$ can be approximated by an arbitrarily close smooth embedding $\phi:M\to\R^{2d+1}$ whose image $\phi(M)$ is a nonsingular algebraic subset of~$\R^{2d+1}$. In this article, we prove that $\phi$ can be chosen in such a way that $\phi(M)$ is a $\Q$-nonsingular $\Q$-algebraic subset of~$\R^{2d+1}$. This guarantees for the first time that, up to smooth diffeomorphisms, every compact smooth manifold $M$ can be described both globally and locally by means of finitely many exact data, such as a finite system of generators of the ideal of polynomials in $\Q[\x_1,\ldots,\x_{2d+1}]$ vanishing on $\phi(M)$.

We extend our result to the singular setting by proving that every real algebraic set with finitely many singularities is semialgebraically homeomorphic to a $\Q$-algebraic set with the same number of singularities. A first consequence is that every affine compact or noncompact Nash manifold of dimension $d$ is Nash diffeomorphic to a $\Q$-nonsingular $\Q$-algebraic subset of~$\R^{2d+1}$. Another consequence is that every germ of a real algebraic set with an isolated singularity is semialgebraically homeomorphic to the germ of a $\Q$-algebraic set with an isolated singularity.
\end{abstract}

$\;$

\vspace{-3em}

\maketitle 


$\;$

\vspace{-2em}

\renewcommand{\contentsname}{Table of Contents}

{\small
\setcounter{tocdepth}{2}
\tableofcontents
}


\section{Introduction and main results}

\subsection{Introduction}\label{subsec:intro}
The description of geometric objects in the simplest possible terms is one of the goals of geometry, in particular of real algebraic geometry. Indeed, the need to simplify the description of manifolds and even topological spaces with singularities, and to find increasingly rich structures on them, has contributed significantly to the birth and development of real algebraic geometry.

In what follows, $\R^n$ and each of its subsets are endowed with the Euclidean topology, `manifold' means `non-empty manifold without bounda\-ry', `smooth(ly)' means `of class $\cinfty$' and a `smooth submanifold' is understood in the usual sense (\cite[p.13]{hirsch:difftop}).

In 1936, Whitney \cite{whit:cobordism} proved that every smooth manifold $M$ of dimension $d$ can be smoothly embedded in $\R^{2d+1}$ and every smooth embedding $\psi:M\to\R^{2d+1}$ can be approximated by an arbitrarily close smooth embedding $\phi:M\to\R^{2d+1}$ whose image $M':=\phi(M)$ is a real analytic submanifold of $\R^{2d+1}$. It follows that $M$ can be described both globally and locally by means of real analytic equations in some Euclidean space. Indeed, one first identifies $M$ with $M'\subset\R^{2d+1}$ via $\phi$ and then observes that, by Cartan's Theorem B, $M'$ is the set of solutions of finitely many global real analytic equations $f_1=0,\ldots,f_s=0$, i.e., each $f_i$ is a (real-valued) real analytic function defined on the whole $\R^{2d+1}$. In addition, locally at each of its points~$a$, $M'$ is the set of solutions of good real analytic equations, good in the sense of the Implicit Function Theorem, i.e., there exist an open neighborhood $U$ of $a$ in $\R^n$ and real analytic functions $g_1,\ldots,g_{d+1}$ defined on $U$ such that the gradients $\nabla g_1(a),\ldots,\nabla g_{d+1}(a)$ are linearly independent and $M'\cap U=\{x\in U:g_1(x)=0,\ldots,g_{d+1}(x)=0\}$. In particular, the smooth manifold $M$ admits a real analytic structure. At this point it is natural to wonder if $M'$ also admits what we could call a real algebraic structure, obtained by requiring that the previous real analytic equations describing~$M'$, both globally and locally, are actually real polynomial equations, i.e., polynomial equations with coefficients in $\R$.

Recall that a subset $X$ of $\R^n$ is called algebraic if it is the set of solutions of real polynomial equations in $\R^n$. A point $a$ of the algebraic set $X\subset\R^n$ is said to be nonsingular of dimension $d$ if there exist an open neighborhood $U$ of $a$ in $\R^n$ and real polynomials $g_1,\ldots,g_{n-d}\in\R[\x_1,\ldots,\x_n]$ vanishing on the whole $X$ such that the gradients $\nabla g_1(a),\ldots,\nabla g_{n-d}(a)$ are linearly independent and $X\cap U=\{x\in U:g_1(x)=0,\ldots,g_{n-d}(x)=0\}$. If all points of $X$ are nonsingular of the same dimension, then $X$ is said to be nonsingular. A (nonsingular) real algebraic set is a (nonsingular) algebraic subset of some $\R^n$. These are basic concepts of real algebraic geometry. Our standard reference for real algebraic geometry is \cite{BCR}.

In his groundbreaking article \cite{nash} published in 1952, Nash proved that, if the given smooth manifold $M$ of dimension $d$ is compact, then we can assume that the real analytic submanifold $M'$ of $\R^{2d+1}$ approximating $\psi(M)$ is actually a union of some connected components of an algebraic subset $X$ of $\R^{2d+1}$ and all points of $X$ belonging to $M'$ are nonsingular of dimension~$d$. Nash~conjectured that $M'$ can be chosen to be a whole nonsingular algebraic subset of $\R^{2d+1}$, a so-called algebraic model of $M$. In 1973, Tognoli \cite{togn:algmodel} proved this conjecture to be true, i.e., if $M$ is a compact smooth manifold of dimension $d$ and $\psi:M\to\R^{2d+1}$ is a smooth embedding, then $\psi$ can be approximated by an arbitrarily close smooth embedding $\phi:M\to\R^{2d+1}$ whose image $M':=\phi(M)$ is an algebraic model of $M$, i.e., $M'$ is a nonsingular algebraic subset of~$\R^{2d+1}$. This is the celebrated Nash-Tognoli theorem. In 1976, King \cite{king76} proved the projective version of the latter result, in which $\R^{2d+1}$ is replaced with $\PP^{2d+1}(\R)$.

The Nash-Tognoli theorem marked the birth of modern real algebraic geometry. Quoting from \cite[p.4]{BCR}: ``A systematic study of real algebraic varieties started seriously only in 1973 after Tognoli's surprising discovery (based on earlier work of John Nash) that every compact smooth manifold is diffeomorphic to a nonsingular real algebraic set''. There is a wide literature devoted to this result and its extensions. We refer the reader to the books \cite[Sect.II.8]{akbking:tras}, \cite[Ch.14]{BCR}, \cite[Ch.5]{man}, the articles \cite{Be2022,GT2017,kucharz2011} and the references therein. See also the surveys \cite[Sect.1]{delellis} and \cite[Sects.1\&2]{kollar2017} on Nash's work for other enlightening presentations of this crucial result.

The problem of making topological spaces algebraic has also been studied in the singular case.

In 1981, Akbulut and King \cite{ak1981} obtained a complete description of the topology of real algebraic sets with isolated singularities. Their idea is to consider the family $\mc{T}$ of compact topological spaces $X$ that admits a topological {resolution of singularities} in the following sense: there exist a compact smooth mani\-fold $M$, a finite family $\mc{M}=\{M_1,\ldots,M_s\}$ of pairwise disjoint subsets of $M$ and a finite set $S=\{p_1,\ldots,p_t\}$ with $t\geq s$ such that each $M_i$ is a finite union of smooth hypersurfaces of $M$ in general position and the quotient topological space obtained from $M\sqcup S$ by blowing down each $M_i$ to $p_i$ is homeomorphic to $X$. Note that, if $\dim(M)>0$ and $t>s$, then $p_{s+1},\ldots,p_t$ correspond to the isolated points of $X$. The topological data $(M,\mc{M};S)$ represents a topological {resolution} of $X$. By Hironaka's {resolution of singularities} \cite{Hi64} (see also \cite{BM,Ko}), the family $\mc{T}$ includes all compact real algebraic sets with isolated singularities. Now, the strategy of Akbulut and King is first to make algebraic the topological data $(M,\mc{M};S)$ associated to $X$, obtaining real algebraic data $(M',\mc{M}';S)$ by mean of algebraic approximation techniques a l\`a Nash-Tognoli, and then to blow down these real algebraic data obtaining a real algebraic set homeomorphic to~$X$. It follows that, up to homeomorphisms, the family $\mc{T}$ coincides with the family of all compact real algebraic sets with isolated singularities.

As Alexandrov's compactification of a real algebraic set can be made algebraic, it follows that a (not necessarily compact) topological space $X$ is homeomorphic to a real algebraic set with isolated singularities if and only if it can be obtained from topological data $(M,\mc{M};S)$ such as the above by considering $M\sqcup S$, blowing down some of the $M_i$ to the corresponding point $p_i$ and removing the remaining $M_i$ and the corresponding point $p_i$ (see \cite[Sect.4]{ak1981}). As a consequence, a noncompact smooth manifold is smoothly diffeomorphic to a nonsingular real algebraic set if and only if it is smoothly diffeomorphic to the interior of a compact smooth manifold with non-empty boundary (see \cite[Cor.4.3]{ak1981}). Further developing the above topological resolution technique, Akbulut and King also proved that every compact real analytic set of dimension $\leq3$ is homeomorphic to a real algebraic~set (see \cite[Ch.VII]{akbking:tras}).

\vspace{1em}

Now it is natural to go a step further and ask whether the description of a geometric object admitting a real algebraic structure can be simplified by requiring that the coefficients of the describing real polynomial equations belong to a subfield $K$ of $\R$ as small as possible. Here the final goal is $K=\Q$, the field of rational numbers which is the smallest subfield of $\R$. We refer to the above open question for $K=\Q$ as the \emph{$\Q$-algebraicity problem}, see \cite[Sect.4.4.2]{P2021}.

The answer is affirmative if $K$ is the field $\qr$ of real algebraic numbers, the smallest real closed field. This is due to three of the most important results in semialgebraic and Nash geometry. Our standard reference for semialgebraic and Nash geometry is \cite{BCR} (see also \cite{benris:semi-alg,Sh}). Let $X$ be an algebraic subset of $\R^n$. Choose a description of $X$:
\[
X=\{x\in\R^n:f_1(a,x)=\ldots=f_s(a,x)=0\}
\]
for some polynomials $f_i\in\Z[\mathtt{a}_1,\ldots,\mathtt{a}_m,\x_1,\ldots,\x_n]$, where $a=(a_1,\ldots,a_m)\in\R^m$ is the vector of all coefficients (ordered in some way) that appear in a fixed real polynomial system $f_1=0,\ldots,f_s=0$ in $\R^n$ whose set of solutions is $X$. Define
\[
V:=\{(b,x)\in(\qr)^{m+n}:f_1(b,x)=\ldots=f_s(b,x)=0\}
\]
and denote by $\pi:V\to(\qr)^m$ the projection $(b,x)\mapsto b$. By Hardt's trivialization theorem \cite{hardt1980}, there exists a finite semialgebraic partition $\{M_i\}_{i=1}^t$ of $(\qr)^m$ and, for each $i\in\{1,\ldots,t\}$, an algebraic subset $F_i$ of $(\qr)^n$ and a semialgebraic homeomorphism $h_i:M_i\times F_i\to V\cap\pi^{-1}(M_i)$ compatible with $\pi$. By the Tarski-Seidenberg principle \cite[Ch.5]{BCR}, we can extend coefficients from $\qr$ to $\R$, obtaining a semialgebraic partition $\{(M_i)_\R\}_{i=1}^t$ of $\R^m$ and  semialgebraic homeomorphisms $(h_i)_\R:(M_i)_\R\times(F_i)_\R\to V_\R\cap\pi_\R^{-1}((M_i)_\R)$ compatible with $\pi_\R$, where $\pi_\R:\R^{m+n}\to\R^m$ is the projection $(b,x)\mapsto b$. It follows that $a$ belongs to $(M_j)_\R$ for a unique $j\in\{1,\ldots,t\}$, and hence $X$ is semialgebraically homeomorphic to $(F_j)_\R$. Note that the set $(F_j)_\R\subset\R^n$ is $\qr$-algebraic in the sense that there exist polynomials $f'_1,\ldots,f'_t\in\qr[\x_1,\ldots,\x_n]$ such that $(F_j)_\R=\{x\in\R^n:f'_1(x)=0,\ldots,f'_t(x)=0\}$. This proves that, up to semialgebraic homeomorphisms, every real algebraic set can be globally described as the set of solutions of finitely many polynomial equations with coefficients in $\qr$.

In \cite[Thm.A]{CS1992} Coste and Shiota proved a version of Hardt's trivialization theorem for Nash manifolds. As a consequence, if the algebraic set $X\subset\R^n$ is nonsingular, we can assume that the algebraic set $F_j\subset(\qr)^n$ is nonsingular and $X$ is Nash diffeomorphic to $(F_j)_\R$. Note that, in this situation, the set $(F_j)_\R\subset\R^n$ is not only $\qr$-algebraic but also $\qr$-nonsingular in the sense that, if $d$ is the dimension of $X$ and $a$ is any point of $X$, then there exist an open neighborhood $U$ of $a$ in $\R^n$ and polynomials $g_1,\ldots,g_{n-d}\in\qr[\x_1,\ldots,\x_n]$ vanishing on the whole $X$ such that the gradients $\nabla g_1(a),\ldots,\nabla g_{n-d}(a)$ are linearly independent and $X\cap U=\{x\in U:g_1(x)=0,\ldots,g_{n-d}(x)=0\}$. This gives us the following improved version of the Nash-Tognoli theorem, which we can call the Nash-Tognoli theorem over $\qr$: {\it Every compact smooth manifold $M$ of dimension $d$ has an algebraic model $M'\subset\R^{2d+1}$ that is $\qr$-algebraic and $\qr$-nonsingular}. 

As $\Q$ is not a real closed field, none of the above results by Hardt, Tarski-Seidenberg and Coste-Shiota are available in the case $K=\Q$. In \cite[Thm.11 \& Rmk.13]{PR2020}, using Zariski equisingular deformations of the coefficients $a=(a_1,\ldots,a_m)$, Parusi\'{n}ski and Rond proved that, if the field extension of $\Q$ obtained by adding $a_1,\ldots,a_m$ is purely transcendental, the algebraic set $X\subset\R^n$ is homeomorphic to an algebraic subset of $\R^n$ defined as the set of solutions of finitely many polynomial equations with coefficients in $\Q$. The reason why this approach does not provide a complete solution of the case $K=\Q$, even when $X$ is a compact nonsingular algebraic hypersurface of $\R^n$, is that the above Zariski equisingular deformations of the coefficients $a=(a_1,\ldots,a_m)$ preserve the polynomial relations over $\Q$ satisfied by $a_1,\ldots,a_m$. In general, the case $K=\Q$ remains open, even in the case of germs (see \cite[Open problems 1 \&~2, pp.199-200]{P2021}). For the links between the $\Q$-algebraicity problem and equisingular deformations, see also \cite{PP,Ro,Te90} and the references mentioned therein.

The goal of this article is to give a complete affirmative solution to the case $K=\Q$ for all nonsingular real algebraic sets and, more in general, for all real algebraic sets with isolated singularities, i.e., with a finite number of singularities. As a byproduct, we obtain an affirmative solution to the case $K=\Q$ for all real algebraic germs with an isolated singularity.

\subsection{Main results} To present our results, we need some preparation.

Let us recall some concepts introduced in \cite{FG}. Let $n\in\N^*:=\N\setminus\{0\}$, let $K:=\R$ or $\Q$, and let $K[\x]:=K[\x_1,\ldots,\x_n]$. Consider $\Q[\x]$ as a subset of $\R[\x]$. Given sets $F\subset\R[\x]$ and $X\subset\R^n$, define
\begin{align}
\ZZ_\R(F)&:=\{x\in\R^n: f(x)=0,\ \forall f\in F\},\label{1}\\
\II_K(X)&:=\{f\in K[\x]: f(x)=0,\ \forall x\in X\}.\label{2}
\end{align}

Observe that $\II_\R(X)$ is the usual vanishing ideal of $X$ in $\R[\x]$ and $\II_\Q(X)=\II_\R(X)\cap\Q[\x]$ is an ideal of $\Q[\x]$. For short, if $F=\{f_1,\ldots,f_s\}\subset\R[\x]$ is finite, we set $\ZZ_\R(f_1,\ldots,f_s):=\ZZ_\R(F)$. 

\begin{defn}[{\cite[Def.2.1.1]{FG}}]\label{def:K-alg}
Let $X$ be a subset of $\R^n$. We say that $X$ is \emph{$K$-algebraic}, or $X\subset\R^n$ is a \emph{$K$-algebraic set}, if $X=\ZZ_\R(F)$ for some $F\subset K[\x]$. $\sqbullet$
\end{defn}

An $\R$-algebraic subset of $\R^n$ is a usual algebraic subset of $\R^n$ and a $\Q$-algebraic subset of $\R^n$ is also an algebraic subset of $\R^n$. The singletons $\{\sqrt{2}\}$ and $\{q\}$, where $q$ is any transcendental real number, are examples of algebraic subsets of $\R$ that are not $\Q$-algebraic. The singleton $\{\sqrt[3]{2}\}=\ZZ_\R(\x_1^3-2)\subset\R$ is $\Q$-algebraic.

\begin{notation}
In what follows, for every set $X\subset\R^n$, we denote $\dim(X)$ the dimension of the Zariski closure of $X$ in $\R^n$, i.e., $\dim(X)$ is the Krull dimension of the ring $\R[\x]/\II_\R(X)$. We say that $X$ has dimension $\dim(X)$. $\sqbullet$
\end{notation}

\begin{defn}\label{13}
Let $X\subset\R^n$ be a $K$-algebraic set of dimension~$d$ and let $a\in X$. We say that $a$ is a \emph{$K$-nonsingular point of $X$} if either $d=n$ (or, equivalently, $X=\R^n$) or $d<n$ and there exist polynomials $f_1,\ldots,f_{n-d}\in\II_K(X)$ and a neighborhood $U$ of $a$ in $\R^n$ such that the gradients $\nabla f_1(a),\ldots,\nabla f_{n-d}(a)$ are linearly independent as vectors in $\R^n$ and $X\cap U=\ZZ_\R(f_1,\ldots,f_{n-d})\cap U$. We say that $a$ is a \emph{$K$-singular point of $X$} if it is not $K$-nonsingular. We denote $\Reg^K(X)$ the set of $K$-nonsingular points of~$X$ and $\Sing^K(X):=X\setminus\Reg^K(X)$ the set of $K$-singular points of $X$. If all the points of $X$ are $K$-nonsingular, i.e., $X=\Reg^K(X)$, then we say that $X$ is \emph{$K$-nonsingular}. $\sqbullet$
\end{defn}

Let $X\subset\R^n$ be an algebraic set of dimension $d$. Observe that an $\R$-nonsingular/$\R$-singular point of $X$ is a usual nonsingular/singular point of $X$, $\Reg^\R(X)$ is the usual set $\Reg(X)$ of nonsingular points of $X$, $\Sing^\R(X)$ is the usual set $\Sing(X)$ of singular points of $X$ and {an $\R$-nonsingular} $\R$-algebraic subset of $\R^n$ is a usual nonsingular algebraic subset of~$\R^n$, see \cite[Sect.3.3]{BCR}. Recall that $\Sing(X)\subset\R^n$ is an algebraic set of dimension $<d$ (see \cite[Prop.3.3.14]{BCR}) and $\Reg(X)$ is a Nash submanifold of~$\R^n$ of dimension $d$ (see \cite[Def.2.9.9 \& Prop.3.3.11]{BCR}). We say that $X$ has \emph{isolated singularities} if it has finitely many singular points, i.e., $\Sing(X)$ is discrete in $\R^n$ and therefore finite.

\begin{remark}\label{rem1}
Let $X\subset\R^n$ be a $\Q$-algebraic set. By the very definitions, every $\Q$-nonsingular point of $X$ is also nonsingular, i.e., $\Reg^\Q(X)\subset\Reg(X)$ or, equivalently, $\Sing^\Q(X)\supset\Sing(X)$. These inclusions are strict in general. Let us give two simple examples.

Consider the polynomials $f_1:=\x_1^2+\x_2^2-\sqrt[3]{2}\x_1\in\R [\x]:=\R[\x_1,\x_2]$ and $g_1:=(\x_1^2+\x_2^2)^3-2\x_1^3\in\Q[\x]$, and the circumference $C_1:=\ZZ_\R(f_1)$ of $\R^2$, which is $\Q$-algebraic as $C_1$ coincides with $\ZZ_\R(g_1)$. Since $f_1$ is irreducible in $\R[\x]$, $g_1$ is irreducible in $\Q[\x]$ and both $f_1$ and $g_1$ change sign in $\R^2$, we have that $\II_\R(C_1)=(f_1)\R[\x]$ and $\II_\Q(C_1)=(g_1)\Q[\x]$, see \cite[Thm.4.5.1]{BCR} and \cite[Prop.3.2.4]{FG}. It follows that each point of $C_1$ is nonsingular as $\nabla f_1$ never vanishes on $C_1$. By contrast, the origin $O$ of $\R^2$ is a $\Q$-singular point of $C_1$ since the restriction of $\nabla g_1$ to $C_1$ vanishes (exactly) at $O$. Therefore, $\Reg^\Q(C_1)=C_1\,\setminus\{O\}\subsetneqq C_1=\Reg(C_1)$ so $C_1$ is nonsingular but not $\Q$-nonsingular.

Here is a variant of the previous example. Let $f_2:=\x_2^2+\sqrt[3]{2}\x_1(\x_1-1)^3\in\R[\x]$, let $g_2:=\x_2^6+2\x_1^3(\x_1-1)^9\in\Q[\x]$ and let $C_2$ be the algebraic curve $\ZZ_\R(f_2)$ of $\R^2$, which is $\Q$-algebraic as $C_2=\ZZ_\R(g_2)$. The curve $C_2$ is homeomorphic to $C_1$ and has a cusp at the point $p:=(1,0)$. Actually, the zero sets of the restrictions of $\nabla f_2$ and $\nabla g_2$ to $C_2$ are $\{p\}$ and $\{O,p\}$, respectively. Thus, $\Reg^\Q(C_2)=C_2\setminus\{O,p\}\subsetneqq C_2\setminus\{p\}=\Reg(C_2)$.

See \cite[Ex.5.4.2]{FG} for additional examples. $\sqbullet$ 
\end{remark}

\begin{defn}
Let $M$ be a smooth manifold and let $M'\subset\R^m$ be a $\Q$-nonsingular $\Q$-algebraic set. If $M$ is smoothly diffeomorphic to $M'$, then we say that $M'$ is a \emph{$\Q$-algebraic model of~$M$}.~$\sqbullet$
\end{defn}

Let $\PP^n(\R)$ be the real projective $n$-space, let $[x_0,x]=[x_0,x_1,\ldots,x_n]$ be its homogeneous coordinates, let $\Q[\x_0,\x]:=\Q[\x_0,\x_1,\ldots,\x_n]$ and let $S$ be a subset of $\PP^n(\R)$. According to \cite[Def.2.6.1]{FG}, we say that $S$ is a \emph{$\Q$-algebraic subset of $\PP^n(\R)$} if it is the common zero set of a family of homogeneous polynomials in $\Q[\x_0,\x]$. Denote $\xi:\R^n\to\PP^n(\R)$ the standard parametrization $\xi(x):=[1,x]$.

\begin{defn}\label{def:proj-Q-closed}
Given a set $X\subset\R^n$, we say that $X$ is \emph{projectively $\Q$-closed} if $\xi(X)$ is a $\Q$-algebraic subset of $\PP^n(\R)$. $\sqbullet$ 
\end{defn}

Given a smooth manifold $M$, we denote $\cinfty(M,\R^m)$ the set of smooth maps from $M$ to $\R^m$ endowed with the usual weak $\cinfty$ topology, see the definition of `$C^\infty_W(M,N)$' in \cite[p.36]{hirsch:difftop}.  

For every $n,m\in\N$ with $m>n$, we identify $\R^n$ with the subset $\R^n\times\{0\}$ of $\R^n\times\R^{m-n}=\R^m$, so we can write $\R^n\subset\R^m$ and every subset of $\R^n$ is also a subset of $\R^m$.

Our first main result is the following improved version of the Nash-Tognoli theorem over $\qr$, which we refer to as the Nash-Tognoli theorem over the rationals.

\begin{thm}\label{thm:NTQ}
Every compact smooth manifold has a $\Q$-algebraic model.

More precisely, if $M$ is a compact smooth manifold of dimension $d$, $\psi:M\to\R^{2d+1}$ is a smooth embedding and $\,\mc{V}$ is a neighborhood of $\psi$ in $\cinfty(M,\R^{2d+1})$, then there exists a smooth embedding $\phi:M\to\R^{2d+1}$ belonging to $\mc{V}$ such that $M':=\phi(M)\subset\R^{2d+1}$ is a $\Q$-nonsingular $\Q$-algebraic set. In addition, we can assume that $M'\subset\R^{2d+1}$ is projectively $\Q$-closed.
\end{thm}

It is worth mentioning a few remarks.

\begin{remarks}\label{rem17}
$(\mr{i})$ If $M$ is a compact Nash submanifold of some $\R^n$ (for example, a compact nonsingular algebraic subset of $\R^n$) of dimension $d$, then we can assume that $\phi$ is a Nash embedding, i.e., a smooth embedding that is also a Nash map. To do that, it suffices to combine Theorem~\ref{thm:NTQ} with the standard Nash approximation result \cite[Cor.8.9.7]{BCR}. In particular, $M$ is Nash diffeomorphic to a $\Q$-nonsingular $\Q$-algebraic subset of $\R^{2d+1}$. The same holds in the noncompact case, see Corollary \ref{cor:main-3}$(\mr{i})$.

$(\mr{ii})$ In the statement of Theorem \ref{thm:NTQ}, we may replace $2d+1$ with any natural number $\geq2d+1$. See Remark \ref{rem112}.

$(\mr{iii})$ If $M$ is a compact smooth hypersurface of $\R^n$, then we can assume that the $\Q$-algebraic model $M'$ of $M$ is a subset of $\R^n$. To prove this, consider a compact neighborhood $K$ of $M$ in $\R^n$ and a smooth function $f\in\cinfty(\R^n)$ such that $\ZZ_\R(f)=M$, $0$ is a regular value of $f$ and $f(x)=2$ for all $x\in\R^n\setminus K$. Then there exist polynomials $p,q\in\Q[\x]$ such that $\ZZ_\R(q)=\varnothing$, the regular function $g:=\frac{p}{q}$ on $\R^n$ is arbitrarily close to $f-2$ in $\cinfty(\R^n)$ and $\sup_{x\in\R^n\setminus K}|g(x)|<1$ (see Lemma \ref{lem:L} below with $L:=\varnothing$). It follows that the regular function $g':=g+2=\frac{p+2q}{q}$ on $\R^n$ is arbitrarily close to $f$, $\inf_{x\in\R^n\setminus K}g'(x)>1$ and $M':=\ZZ_\R(g')=\ZZ_\R(p+2q)$ is the $\Q$-algebraic model of $M$ in $\R^n$ we were looking for.

$(\mr{iv})$ Theorem \ref{thm:NTQ} guarantees for the first time that it is theoretically possible to store any compact smooth manifold $M$, up to smooth diffeomorphisms, on a computer system by means of a finite number of exact data, such as a finite system of generators of $\II_\Q(M')$ in $\Q[\x_1,\ldots,\x_{2d+1}]$.~$\sqbullet$   
\end{remarks}

Let us introduce the concept of $\Q$-determined $\Q$-algebraic set.

\begin{defn}\label{Q-determined}
Let $X\subset\R^n$ be a $\Q$-algebraic set. We say that $X$ is \emph{$\Q$-determined} if every nonsingular point of $X$ is also $\Q$-nonsingular, i.e., $\Reg^\Q(X)=\Reg(X)$. $\sqbullet$ 
\end{defn}

If $X\subset\R^n$ is a $\Q$-determined $\Q$-algebraic set, then the set $\Sing(X)=X\setminus\Reg(X)$ of singular points of $X$ coincides with the set $\Sing^\Q(X)=X\setminus\Reg^\Q(X)$ of $\Q$-singular points of $X$.

The curves $C_1$ and $C_2$ of $\R^2$ defined in Remark \ref{rem1} are examples of $\Q$-algebraic sets that are not $\Q$-determined, the first nonsingular and the second with only one singular point.

Given a topological space $S$, we denote $\czero(S,\R^m)$ the set of continuous maps from $S$ to $\R^m$ endowed with the usual compact-open topology, see \cite[Sect.7-5]{Munkres}. A map between topological spaces is a continuous embedding if it is a homeomorphism onto its image.

Our second main result is as follows.

\begin{thm}\label{thm:main}
Every real algebraic set with isolated singularities is semialgebraically homeo\-morphic to a $\Q$-determined $\Q$-algebraic set with isolated singularities. 

More precisely, the following holds. Let $X$ be an algebraic subset of $\R^n$ of dimension $d$ with isolated singularities, let $m:=n+2d+4$, and let $i:X\hookrightarrow\R^m$ and $j:\Reg(X)\hookrightarrow\R^m$ be the inclusions of $X$ and $\Reg(X)$ in $\R^m$. Then, for every neighborhood $\mc{U}$ of $i$ in $\czero(X,\R^m)$ and for every neighborhood $\mc{V}$ of $j$ in $\cinfty(\Reg(X),\R^m)$, there exists a semialgebraic continuous embedding $\phi:X\to\R^m$ with the following properties:
\begin{itemize}
 \item[$(\mr{i})$] $X':=\phi(X)$ is a $\Q$-determined $\Q$-algebraic subset of $\R^m$ such that $\Sing(X')\subset\phi(\Sing(X))$ or, equivalently, $\phi(\Reg(X))\subset\Reg(X')$.
 \item[$(\mr{ii})$] The restriction $\phi|_{\Reg(X)}:\Reg(X)\to\R^m$ of $\phi$ to $\Reg(X)$ is a Nash embedding. In particular, the Nash manifolds $\Reg(X)$ and $\phi(\Reg(X))$ are Nash diffeomorphic.
 \item[$(\mr{iii})$] $\phi\in\mc{U}$ and $\phi|_{\Reg(X)}\in\mc{V}$.
\end{itemize}
\end{thm}

\begin{remarks}
$(\mr{i})$ If the nature of a singular point of $X$ is not topological, i.e., the link of $X$ at that point is homeomorphic to a $(d-1)$-sphere, then there could exist a $\Q$-approximation $X'$ of $X$ in which this singularity disappears, i.e., the inclusion $\Sing(X')\subset\phi(\Sing(X))$ can be strict.

Here is an example. Consider the algebraic cusp $X:=\ZZ_\R(\x_1^2-\sqrt{2}\x_2^3)\subset\R^2$, which has a unique singularity at the origin $O$ of $\R^2$ and is homeomorphic to $\R$. Let $\alpha>0$ be a small rational number, let $\beta$ be a real number close to $1$ such that $\gamma:=\beta\sqrt{2}\in\Q$, let $\psi:\R^2\to\R^2$ be the homeomorphism $\psi(x_1,x_2):=(x_1,\sqrt[3]{\beta(\alpha x_2+x_2^3)})$, let $X':=\psi^{-1}(X)$ and let $\phi:X\to\R^2$ be the semialgebraic continuous embedding $\phi(x):=\psi^{-1}(x)$. Since $X'=\ZZ_\R(\x_1^2-\gamma(\alpha\x_2+\x_2^3))$, it follows that $X'=\phi(X)\subset\R^2$ is a $\Q$-nonsingular $\Q$-algebraic set. In addition, if $\alpha$ is sufficiently small and $\beta$ is sufficiently close to $1$, then $\phi$ is arbitrarily $\czero$ close to the inclusion $X\hookrightarrow\R^2$ and $\phi|_{X\setminus\{O\}}$ is a Nash embedding arbitrarily $\cinfty$ close to the inclusion $X\setminus\{O\}\hookrightarrow\R^2$.

$(\mr{ii})$ It is always possible to require that $X'$ has the same number of singularities of $X$, i.e., $\Sing(X')=\phi(\Sing(X))$ and thus $\Reg(X')=\phi(\Reg(X))$, see Remarks \ref{rem:degen} and \ref{rem:degen2}. $\sqbullet$
\end{remarks}

If we are willing to lose approximation properties $(\mr{iii})$ in the statement of Theorem \ref{thm:main}, then we can drop the embedding dimension of $X'$.

\begin{thm}\label{thm:main-2}
Every real algebraic set $X$ of dimension $d$ with isolated singularities is semi\-algebraically homeomorphic to a $\Q$-determined $\Q$-algebraic subset $X'$ of $\R^{2d+5}$.

More precisely, there exists a semialgebraic homeomorphism $\eta:X\to X'$ such that $\eta(\Reg(X))=\Reg(X')$ and the restriction of $\eta$ from $\Reg(X)$ to $\Reg(X')$ is a Nash diffeomorphism.
\end{thm}

In \cite[Rem.VI.2.11, p.208]{Sh}, Shiota proved that a noncompact Nash submanifold of some $\R^n$ is Nash diffeomorphic to a nonsingular real algebraic set. This result, the noncompact nonsingular case of the last two theorems and the $\R|\Q$-generic projection theorem \cite[Cor.6.3.5]{FG} have the following immediate consequence:

\begin{cor}\label{cor:main-3}
Let $M$ be a noncompact Nash submanifold of $\R^n$ of dimension $d$. We have:
\begin{itemize}
 \item[$(\mr{i})$] $M$ is Nash diffeomorphic to a $\Q$-nonsingular $\Q$-algebraic subset of $\R^{2d+1}$.
 \item[$(\mr{ii})$] If $M$ is a noncompact nonsingular algebraic subset of $\R^n$ of dimension $d$, $m:=n+2d+4$ and $\mc{V}$ is a neighborhood of the inclusion map $M\hookrightarrow\R^m$ in $\cinfty(M,\R^m)$, then there exists a Nash embedding $\phi:M\to\R^m$ belonging to $\mc{V}$ such that $\phi(M)\subset\R^m$ is a $\Q$-nonsingular $\Q$-algebraic set.
\end{itemize}
\end{cor}

Let $X\subset\R^n$ be a semialgebraic set and let $a$ be a point of $X$. Recall that the local dimension $\dim(X_a)$ of $X$ at $a$ is the dimension of a sufficiently small semialgebraic neighborhood of $a$ in~$X$ (see \cite[Def.2.8.11]{BCR}). Another consequence of Theorem \ref{thm:main-2} (or better to say of its proof) is the following result concerning real algebraic germs with an isolated singularity.

\begin{thm}\label{thm:main-germs}
Every real algebraic set germ with an isolated singularity is semialgebraically homeomorphic to a real $\Q$-determined $\Q$-algebraic set germ with an isolated singularity. 

More precisely, the following holds. Let $X\subset\R^n$ be an algebraic set containing the origin~$O$ of $\R^n$, and let $d:=\dim(X_O)$. Suppose that $O$ is an isolated point of $\Sing(X)$. Then there exist a compact $\Q$-determined $\Q$-algebraic set $X'\subset\R^{2d+4}$ that contains the origin $O'$ of $\R^{2d+4}$, a semialgebraic open neighborhood $U$ of $O$ in $X$, a semialgebraic open neighborhood $U'$ of $O'$ in $X'$ and a semialgebraic homeomorphism $\phi:U\to U'$ such that $\Sing(X)\cap U=\{O\}$, $\Sing(X')\cap U'=\{O'\}$, $\phi(O)=O'$ and the restriction of $\phi$ from $U\setminus\{O\}$ to $U'\setminus\{O'\}$ is a Nash diffeomorphism.
\end{thm}

\subsection{Strategy of the proofs and structure of the article} The proof of the Nash-Tognoli theorem can be sketched as follows.

\begin{proof}[Idea of the proof of the Nash-Tognoli theorem]
Let $M$ be a compact smooth submanifold of $\R^n$ of dimension $d$. Since the unoriented cobordism group $\mk{N}_*$ is generated by nonsingular real algebraic sets, increasing the ambient dimension $n$ if necessary and doubling the cobordism, we can find a compact nonsingular algebraic set $P\subset\R^n$ and a compact smooth submanifold $S$ of $\R^{n+1}$ such that $S$ is transverse to $\R^n\times\{0\}$ in $\R^{n+1}=\R^n\times\R$ and $S\cap(\R^n\times\{0\})$ is the disjoint union $(M\times\{0\})\sqcup(P\times\{0\})$ of $M\times\{0\}$ and $P\times\{0\}$. Let $\rho:V\to S$ be a smooth open tubular neighborhood of $S$ in $\R^{n+1}$, let $\E\subset\R^N$ be the total space of the universal vector bundle over the Grassmannian $\G$ of $n-d$-dimensional vector subspaces of $\R^{n+1}$, and let $B:S\to\G$ be the smooth map sending $a\in S$ to the vector subspace of $\R^{n+1}$ perpendicular to the tangent space of $S$ at $a$. Define the smooth map $\theta:V\to\E$ by $\theta(p):=((B\circ\rho)(p),p-\rho(p))$. We have that $\theta$ is transverse in $\E$ to the zero section $\G\times\{0\}$, $\theta^{-1}(\G\times\{0\})=S$ and the restriction of $\theta$ to $P\times\{0\}$ is a regular map. Applying a relative version of the Weierstrass approximation theorem to $\theta$, we obtain an algebraic set $Z\subset\R^{n+1+N}=\R^{n+1}\times\R^N$ of dimension $n+1$, a smooth map $v:V\to\R^N$ and a regular map $\eta:Z\to\E$ with the following properties: the graph $\widehat{V}$ of $v$ is an open subset of the nonsingular locus $\Reg(Z)$ of $Z$, $v$ vanishes on $P\times\{0\}$ so $P':=P\times\{0\}\times\{0\}\subset\widehat{V}$, $v$ is arbitrarily $\cinfty$ close to the zero map from $V$ to $\R^N$ and, if $(x,x_{n+1},y)$ are the coordinates of $\R^{n+1+N}=\R^n\times\R\times\R^N$ and $\pi:\widehat{V}\to V$ is the smooth diffemorphism $(x,x_{n+1},y)\mapsto(x,x_{n+1})$, then $\eta|_{\widehat{V}}$ is arbitrarily $\cinfty$ close to $\theta\circ\pi$, and $\eta(x,0,0)=\theta(x,0)$ for all $x\in P$.

Denote $Z_1$ the union of all the $d+1$-dimensional irreducible components of the algebraic set $\eta^{-1}(\G\times\{0\})\subset\R^{n+1+N}$, and define $\widehat{S}_1:=\widehat{V}\cap Z_1$. Recall that transversality is an open condition that induces small diffeotopies between the inverse images of sufficiently $\cinfty$ close maps \cite[Thm.14.1.1]{BCR}. It follows that $\eta|_{\widehat{V}}$ is transverse to $\G\times\{0\}$ in $\E$, so $\widehat{S}_1\subset\Reg(Z_1)$. In addition, since $\eta(x,0,0)=\theta(x,0)$ for all $x\in P$, there exists a smooth diffeomorphism $\nu:\widehat{V}\to\widehat{V}$ arbitrarily $\cinfty$ close to the identity $\mr{id}_{\widehat{V}}$ such that $\nu(\pi^{-1}(S))=\widehat{S}_1$ and $\nu=\mr{id}_{\widehat{V}}$ on $P'$. The map $\widehat{\nu}:S\to\R^{n+1+N}$, $(x,x_{n+1})\mapsto\nu(\pi^{-1}(x,x_{n+1}))$ is a smooth embedding arbitrarily $\cinfty$ close to the inclusion map $S\hookrightarrow\R^{n+1+N}$, $(x,x_{n+1})\mapsto(x,x_{n+1},0)$ such that $\widehat{\nu}(S)=\widehat{S}_1$ and $\widehat{\nu}(x,0)=(x,0,0)$ for all $x\in P$. As a consequence, if $\pi_{n+1}:\R^{n+1+N}\to\R$ is the projection $\pi_{n+1}(x,x_{n+1},y):=x_{n+1}$, then $0$ is a regular value of $\pi_{n+1}|_{\widehat{S}_1}$, $\widehat{S}_1\cap\pi_{n+1}^{-1}(0)=M_1\sqcup P'$ for some compact smooth submanifold $M_1$ of $\R^{n+1+N}$ and there exists a smooth embedding $\phi_1:M\to\R^{n+1+N}$ arbitrarily $\cinfty$ close to the inclusion map $M\hookrightarrow\R^{n+1+N}$, $x\mapsto(x,0,0)$ such that $\phi_1(M)=M_1$.

Since $\widehat{S}_1$ is compact in $\R^{n+1+N}$, $Z_1\setminus\widehat{S}_1$ is closed in $\R^{n+1+N}$ and $P'\subset\R^{n+1+N}$ is a compact nonsingular algebraic set contained in $\widehat{S}_1$, we can apply the aforementioned relative Weierstrass approximation theorem to $\pi_{n+1}$, obtaining a regular function $u:Z_1\to\R$ such that $u|_{\widehat{S}_1}$ is arbitrarily $\cinfty$ close to $\pi_{n+1}|_{\widehat{S}_1}$, $u$ vanishes on $P'$ and $u>1$ on $Z_1\setminus\widehat{S}_1$. It follows that $0$ is a regular value of $u$ (thus, $u^{-1}(0)\subset\R^{n+1+N}$ is a nonsingular algebraic set of dimension~$d$), $u^{-1}(0)=M'\sqcup P'$ for some compact smooth submanifold $M'$ of $\R^{n+1+N}$ and there exists a smooth embedding $\widehat{\phi}_1:M_1\to\R^{n+1+N}$ arbitrarily $\cinfty$ close to the inclusion map $M_1\hookrightarrow\R^{n+1+N}$ such that $\widehat{\phi}_1(M_1)=M'$. The map $\phi:M\to\R^{n+1+N}$, $x\mapsto\widehat{\phi}_1(\phi_1(x))$ is a smooth embedding arbitrarily $\cinfty$ close to the inclusion map $M\hookrightarrow\R^{n+1+N}$, $x\mapsto(x,0,0)$ such that $\phi(M)=M'$. Since $u^{-1}(0)$ and $P'$ are nonsingular algebraic subsets of $\R^{n+1+N}$ of the same dimension $d$ such that $P'\subsetneqq u^{-1}(0)$, their difference $M'=u^{-1}(0)\setminus P'$ is also a nonsingular algebraic subset of $\R^{n+1+N}$ of dimension $d$. By the generic projection theorem, we can assume that $M'$ belongs to $\R^{2d+1}$. In fact, if $n+1+N>2d+1$, then the generic projection theorem assures that a generic linear projection $\Pi:\R^{n+1+N}\to\R^{2d+1}$ has the following properties: $\Pi(M')\subset\R^{2d+1}$ is an algebraic set and the restriction of $\Pi$ from $M'$ to $\Pi(M')$ is a biregular isomorphism, so the algebraic set $\Pi(M')\subset\R^{2d+1}$ is also nonsingular.
\end{proof}

Some of the key tools used in the proof we have just outlined are as follows:
\begin{itemize}
 \item[$(\mr{I})$] The algebraicity of the unoriented cobordism group $\mk{N}_*$.
 \item[$(\mr{II})$] The stability of the family of nonsingular real algebraic sets under inverse images of tran\-sverse regular maps: If $X$ and $Y$ are nonsingular real algebraic sets, $W$ is a nonsingular algebraic subset of $Y$ of codimension $c$ and $f:X\to Y$ is a regular map transverse to~$W$, then $f^{-1}(W)$ is a nonsingular algebraic subset of $X$ of codimension $c$.
 \item[$(\mr{III})$] The relative Weierstrass approximation theorem:   If $L\subset\R^n$ is a compact nonsingular algebraic set and $f:\R^n\to\R$ is a smooth function such that $f|_L$ is regular, then there exists a regular function $g:\R^n\to\R$ arbitrarily $\cinfty$ close to $f$ such that $g|_L=f|_L$.
 \item[$(\mr{IV})$] The difference $X\setminus Y$ between nonsingular real algebraic sets $X$ and $Y$ of the same dimension $d$ with $Y\subsetneqq X$ is still a nonsingular real algebraic set of dimension $d$.
 \item[$(\mr{V})$] The generic projection theorem.
\end{itemize}

Our strategy to prove Theorem \ref{thm:NTQ} is to adapt the previous proof of the Nash-Tognoli theorem to the $\Q$-algebraic setting. This is not an easy task. Some time ago, when we started thinking about how to perform this adaptation, we quickly realized that there was no adequate theory of $\Q$-algebraic sets that contained at least the counterparts of the five previous tools $(\mr{I})$-$(\mr{V})$.

An initial issue is to understand what is the `right' counterpart over~$\Q$ of the concept of nonsingular real algebraic set. This issue is rather subtle and lends itself to different solutions. In Definition \ref{13}, we gave a possible solution. Another is the one of `nonsingular algebraic set $X\subset\R^n$ defined over $\Q$'. According to \cite[Def.3, p.30]{togn:instmat}, an algebraic set $X\subset\R^n$ is said to be \emph{defined over $\Q$} if $\II_\R(X)=\II_\Q(X)\R[\x]$, i.e., $\II_\R(X)$ is generated by $\II_\Q(X)$ in $\R[\x]$. The algebraic set $X\subset\R^n$ is said to be a \textit{nonsingular algebraic set defined over $\Q$} if it is both nonsingular and defined over $\Q$. In Proposition \ref{prop:implications}, we will prove that, if $X\subset\R^n$ is a nonsingular algebraic set defined over $\Q$, then $X\subset\R^n$ is also a $\Q$-nonsingular $\Q$-algebraic set. The converse is not true in general. A simple example is the singleton $X:=\{\sqrt[3]{2}\}=\ZZ_\R(\x_1^3-2)\subset\R$, which is $\Q$-nonsingular $\Q$-algebraic but is not defined over $\Q$ as $\II_\R(X)=(\x_1-\sqrt[3]{2})\R[\x_1]\supsetneqq(\x_1^3-2)\R[\x_1]=\II_\Q(X)\R[\x_1]$. In \cite[Thm.2, p.56]{togn:instmat} and \cite[Thm.0.1]{BaTo92}, the authors extended the proof of the Nash-Tognoli theorem to the $\Q$-algebraic setting by using the concept of `nonsingular algebraic set $X\subset\R^n$ defined over $\Q$'. They obtained the following result: \emph{In the statement of Theorem \ref{thm:NTQ}, we can assume that $\phi(M)\subset\R^{2d+1}$ is a nonsingular algebraic set defined over~$\Q$}. This stronger version of Theorem \ref{thm:NTQ} cannot be considered valid. The main reason is that the proofs presented in \cite{togn:instmat,BaTo92} use implicitly the following $\Q$-version of $(\mr{II})$: `If $X\subset\R^n$, $Y\subset\R^m$ and $W\subset\R^m$ are nonsingular algebraic sets defined over~$\Q$ with $W\subset Y$, and $f:X\to Y$ is a $\Q$-regular map transverse to~$W$ (see Definition \ref{def:Q-reg-function} for the concept of $\Q$-regular map), then $f^{-1}(W)\subset\R^n$ is also a nonsingular algebraic set defined over $\Q$'. Unfortunately, the latter assertion is not correct. A counterexample is as follows: the singleton $\{0\}\subset\R$ is a nonsingular algebraic set defined over $\Q$ and the $\Q$-regular function $f:\R\to\R$, $x\mapsto x^3-2$ has $0$ as a regular value, however the algebraic set $f^{-1}(0)=\{\sqrt[3]{2}\}\subset\R$ is not defined over~$\Q$. Another reason is that the proofs presented in \cite{togn:instmat,BaTo92} also explicitly use the following $\Q$-version of $(\mr{III})$, which is a particular case of Lemma \ref{lem:basic-Q-pair}$(\mr{ii})(\mr{iii})$ below: `If $L\subset\R^n$ is a compact nonsingular algebraic set defined over~$\Q$ and $f:\R^n\to\R$ is a smooth function such that $f|_L$ is $\Q$-regular, then there exists a $\Q$-regular function $g:\R^n\to\R$ arbitrarily $\cinfty$ close to $f$ such that $g|_L=f|_L$'. This result is used when $L$ is just a nonsingular $\Q$-algebraic set, without verifying that $L$ is defined over~$\Q$. Unfortunately, it does not work when $L$ is just a nonsingular $\Q$-algebraic set. For example, if $f_1:=\x_1^2+\x_2^2-\sqrt[3]{2}\x_1\in\R [\x]:=\R[\x_1,\x_2]$, $g_1:=(\x_1^2+\x_2^2)^3-2\x_1^3\in\Q[\x]$ and $L$ is the circumference $C_1:=\ZZ_\R(f_1)=\ZZ_\R(g_1)$ of $\R^2$ defined as in Remark \ref{rem1}, then the polynomial function $f_1:\R^2\to\R$ cannot be $\cinfty$ approximated by $\Q$-regular functions $g:\R^2\to\R$ vanishing on $C_1$, as $\nabla f_1(0,0)\neq0$ and $g=hg_1$ for some $\Q$-regular function $h$ so $\nabla g(0,0)=0$.

The phenomenon described above {concerning the existence of real $\Q$-algebraic sets that are not defined over $\Q$} does not occur in the complex case. In fact, Hilbert's Nullstellensatz implies that every complex $\Q$-algebraic set is defined over~$\Q$, see Remarks \ref{rem210}$(\mr{ii})$.

The concepts of real $\Q$-algebraic set, real algebraic set defined over $\Q$, and real $\Q$-determined $\Q$-algebraic set given in Definition \ref{Q-determined} are three possible formal meanings that can be attributed to the statement: the algebraic set $X\subset\R^n$ is `defined over $\Q$'. In Appendix~\ref{appendix-C}, we will introduce two other concepts of this type and study the relationships between all these concepts, proving in particular that they are not equivalent to each other. This study could be used in the future to try to improve the results obtained in this article, primarily Theorem \ref{thm:NTQ}.

The foundational elements and results of a theory of real $\Q$-algebraic sets have recently been developed in \cite{FG}. In the present article, we further develop this theory by identifying counterparts over~$\Q$ of each tool $(\mr{I})$-$(\mr{V})$ that are `consistent with each other' and thus allow us to extend the proof of the Nash-Tognoli theorem to prove Theorem \ref{thm:NTQ}.

Our strategy to prove Theorems \ref{thm:main}, \ref{thm:main-2} and \ref{thm:main-germs} is to adapt to the $\Q$-algebraic setting the topological resolution technique of Akbulut and King \cite{ak1981} we sketched in Subsection~\ref{subsec:intro}. This is not an easy task either. Given a compact real algebraic set $X$ with isolated singularities, our idea is to consider real algebraic data $(M,\mc{M};S)$ given by Hironaka's resolution of singularities applied to $X$. Then we approximate $(M,\mc{M};S)$ by $\Q$-algebraic data $(M',\mc{M}';S')$ and blow down $(M',\mc{M}';S')$ obtaining a $\Q$-algebraic set $X'$, which is $\Q$-determined and semialgebraically homeomorphic to $X$. To implement this idea, we will introduce and study appropriate counterparts over $\Q$ of key tools used in \cite{ak1981}, such as the following:

\begin{itemize}
 \item[$(\mr{VI})$] The algebraic homology of a nonsingular real algebraic set $W$ and its link with the unoriented bordism group $\mk{N}_*(W)$.
 \item[$(\mr{VII})$] The real algebraic blowing down operation.
\end{itemize}

In Sections \ref{sec:basic-Q-theory} and \ref{sec:Q-alg-approx}, we present the mentioned counterparts over $\Q$ of the tools $(\mr{I})$-$(\mr{VII})$, together with our $\Q$-algebraic approximation results and the proof of Theorem \ref{thm:NTQ}. Section \ref{sec:proofs} is devoted to the proofs of Theorems \ref{thm:main}, \ref{thm:main-2} and \ref{thm:main-germs}: in the first two subsections, we develop additional tools that we will then use in the proofs, presented in the third and last subsection.

The proofs of some preliminary results are known to experts or can be obtained by carefully modifying standard arguments. For completeness, we include these proofs in Appendix~\ref{appendix-A}. Appendix \ref{appendix-B} contains the proof of Theorem \ref{thm:1.7}, which is the main result of Subsection \ref{BFR}. 


\section{Some $\Q$-algebraic geometry}\label{sec:basic-Q-theory} 
\emph{Throughout this section, unless otherwise indicated, the sets $\R^n$ and $\C^n=\R^{2n}$ and each of their subsets are endowed with the Euclidean topology.}

This section deals with real $\Q$-algebraic sets and is organized into four subsections. In the first, we recall some useful concepts and results obtained in \cite{FG}. In the second, we obtain new results concerning $\Q$-nonsingular/$\Q$-singular points of $\Q$-algebraic sets, $\Q$-regular maps, projectively $\Q$-closed $\Q$-algebraic sets, and we prove a $\R|\Q$-version of the generic projection theorem. The third subsection is devoted to the study of certain $\Q$-algebraic models of smooth manifolds such as Milnor's generators of the unoriented cobordism group. We use this study in the fourth subsection to investigate the new concepts of projectively $\Q$-algebraic unoriented bordism and homology.


\subsection{Affine and projective preliminaries}\label{subsec:review} \emph{Throughout this subsection, we assume that $L|K$ is one of the following four extensions of fields: $L|K=\R|\R$, $\R|\Q$, $\C|\C$ and $\C|\Q$.}

\subsubsection{Affine preliminaries} Let $n\in\N^*$, let $L[\x]$ be the ring of polynomials $L[\x_1,\dots,\x_n]$ and let $K[\x]$ be its subring $K[\x_1,\dots,\x_n]$. Given $F\subset L[\x]$ and $X\subset L^n$, we define
\begin{align*}
\ZZ_L(F)&:=\{x\in L^n:f(x)=0,\,\forall f\in F\},\\
\I_K(X)&:=\{f\in K[\x]:f(x)=0,\,\forall x\in X\}.
\end{align*}

Observe that $\I_K(X)$ is an ideal of $K[\x]$. If the subset $F$ of $L[\x]$ consists of a finite number of polynomials $f_1,\dots,f_s$, we set $\ZZ_L(f_1,\dots,f_s):=\ZZ_L(F)=\{x\in L^n:f_1(x)=0,\cdots, f_s(x)=0\}$.

In the next definition, we collect some $L|K$-algebraic/topological notions.

\begin{defn}[{\cite[Def.2.1.1\;\&\;2.1.3\;\&\;2.1.4\;\&\;2.1.5\;\&\;2.1.9]{FG}}]\label{K-algebraic-set}
Let $X$ be a subset of $L^n$. We say that $X$ is a \emph{$K$-algebraic subset of $L^n$}, or $X\subset L^n$ is a \emph{$K$-algebraic set}, if $X=\ZZ_L(F)$ for some $F\subset K[\x]$. We say that $X$ is \emph{$K$-constructible} if it is a Boolean combination of $K$-algebraic subsets of $L^n$. The $K$-algebraic subsets of $L^n$ are the closed sets of a topology called \emph{$K$-Zariski topology of $L^n$}. We say that $X$ is \emph{$K$-Zariski closed}, \emph{$K$-Zariski open} or \emph{$K$-irreducible} if $X$ is closed, open or irreducible with respect to the $K$-Zariski topology of $L^n$, respectively. We call \emph{$K$-Zariski closure of $X$ (in $L^n$)} the closure of $X$ in $L^n$ with respect to the $K$-Zariski topology.

We define the \emph{$K$-Zariski topology} of the set $X\subset L^n$ as the relative topology of $X$ induced by the $K$-Zariski topology of $L^n$. A subset of $X$ is \emph{$K$-Zariski closed} or \emph{$K$-Zariski open} if it is closed or open with respect to the $K$-Zariski topology of $X$, respectively.~$\sqbullet$
\end{defn}

\begin{remarks}
$(\mr{i})$ The concepts of $L$-algebraic/$L$-Zariski closed/$L$-Zariski open/$L$-irreducible subsets of $L^n$ and $L$-Zariski closure of a subset of $L^n$ coincide with the usual ones without the prefix~`$L$-'. If $L=\C$, a $L$-constructible subset of $L^n$ is a usual constructible subset of $\C^n$.

$(\mr{ii})$ A $K$-algebraic subset of $L^n$ is also an algebraic subset of $L^n$.

$(\mr{iii})$ If $X\subset L^n$ is $K$-algebraic, then $X$ is $K$-irreducible if it is non-empty and there do not exist $K$-algebraic subsets $X_1$ and $X_2$ of $L^n$ such that $X_1\subsetneqq X$, $X_2\subsetneqq X$ and $X=X_1\cup X_2$. Given any subset $X$ of $L^n$, the set $\ZZ_L(\II_K(X))$ is the $K$-Zariski closure of $X$ in $L^n$. Thus, $X\subset L^n$ is $K$-algebraic if and only if $X=\ZZ_L(\II_K(X))$. For example, if $L|K=\R|\Q$, $X:=\{\sqrt{2}\}$ and $\overline{X}$ is the $\Q$-Zariski closure of $X$ in $\R$, then $\overline{X}=\{-\sqrt{2},\sqrt{2}\}$. Observe that $\overline{X}$ is a $\Q$-irreducible $\Q$-algebraic subset of $\R$, but it has two irreducible components as an algebraic subset of $\R$. $\sqbullet$
\end{remarks}

The $K$-Zariski topology of $L^n$ is Noetherian {since} it is coarser than the usual Zariski topology of $L^n$, which is Noetherian. Thus, by \cite[Prop.1.5, p.5]{ha}, we have:

\begin{lem}
For every $K$-algebraic subset $X$ of $L^n$, there exists a unique finite family $\{X_1,\ldots,X_r\}$ of $K$-irreducible $K$-algebraic subsets of $L^n$ such that $X=\bigcup_{i=1}^rX_i$ and $X_i\not\subset\bigcup_{j\in\{1,\ldots,r\}\setminus\{i\}}X_j$ for all $i\in\{1,\ldots,r\}$. The sets $X_i$ are said to be the \emph{$K$-irreducible components} of~$X$. $\sqbullet$ 
\end{lem}

Let us recall the notion of $K$-dimension.

\begin{defn}[{\cite[Def.2.1.11]{FG}}]\label{def:K-dim}
Let $X$ be a subset of $L^n$. We define the \emph{$K$-dimension $\dim_K(X)$ of $X$ (in $L^n$)} as the Krull dimension of the ring $K[\x]/\II_K(X)$, i.e., $\dim_K(X):=\dim(K[\x]/\II_K(X))$. If $X=\varnothing$, then $\dim_K(X):=-1$.~$\sqbullet$
\end{defn}

If $X$ is an algebraic subset of $L^n$, then $\dim_L(X)$ is the usual (algebraic) dimension of $X$. 

An elementary but useful property of the $K$-dimension is as follows.

\begin{lem}[{\cite[Lem.2.1.13]{FG}}]\label{dimirred}
If $X,Y\subset L^n$ are $K$-algebraic sets such that $Y\subsetneqq X$ and $X$ is $K$-irreducible, then $\dim_K(Y)<\dim_K(X)$.
\end{lem}

The dimension does not depend on the field used to compute it:

\begin{thm}[{\cite[Thm.2.4.2]{FG}}]\label{dimension}
If $X\subset L^n$ is a $K$-algebraic set, then $\dim_L(X)=\dim_K(X)$.
\end{thm}

Thanks to the previous theorem, we can introduce the following notation.

\begin{notation}\label{dime}
For every algebraic subset $X$ of $L^n$, we denote $\dim(X)$ the usual dimension $\dim_L(X):=\dim(L[\x]/\II_L(X))$ of $X$ and we say that $\dim(X)$ is the dimension of $X$. If $X$ is a $K$-algebraic subset of $L^n$, then $X$ is also an algebraic set of $L^n$ and therefore its dimension is $\dim(X):=\dim_L(X)=\dim_K(X)$. $\sqbullet$
\end{notation}

The next definition deals with $K$-nonsingular/$K$-singular points and $K$-Zariski tangent spaces.

\begin{defn}[{\cite[Def.4.1.1\;\&\;4.2.1\;\&\;5.1.1]{FG}}]\label{28}
Let $X\subset L^n$ be a $K$-algebraic set of dimension $d$, let $a=(a_1,\ldots,a_n)\in X$, let $\gtn_a$ be the maximal ideal $(\x_1-a_1,\ldots,\x_n-a_n)L[\x]$ of $L[\x]$ and let $e\in\{0,\ldots,d\}$. We define the \emph{$K$-local ring $\reg^K_{X,a}$ of $X$ at $a$} by $
\reg^K_{X,a}:=L[\x]_{\gtn_a}/(\II_K(X)L[\x]_{\gtn_a})$. We say that $a$ is a \emph{$K$-nonsingular point of $X$ of dimension $e$} if the local ring $\reg^K_{X,a}$ is regular and of dimension $e$. We denote {by} $\Reg^K(X,e)$ the set of all $K$-nonsingular points of $X$ of dimension $e$. A \emph{$K$-nonsingular point of $X$} is a $K$-nonsingular point of dimension $d$. We denote $\Reg^K(X)$ the set of all $K$-nonsingular points of $X$, i.e., $\Reg^K(X):=\Reg^K(X,d)$. If $a$ is not a $K$-nonsingular point of $X$, then we say that $a$ is a \emph{$K$-singular point of $X$}. We denote $\Sing^K(X)$ the set of all $K$-singular points of $X$, i.e., $\Sing^K(X):=X\setminus\Reg^K(X)$. We say that $X$ is \emph{$K$-nonsingular} if $X=\Reg^K(X)$ and that $X$ is \emph{$K$-determined} if $\Reg(X)=\Reg^K(X)$.

The \emph{$K$-Zariski tangent space $T^K_a(X)$ of $X$ at $a$} is the $L$-vector subspace of $L^n$ defined by
\begin{equation*}\label{eq:E|K-Zariskitangspace}
\textstyle
T^K_a(X):=\big\{v\in L^n:\langle\nabla g(a),v\rangle=0\,\text{ for all $g\in\II_K(X)$}\big\},
\end{equation*}
where $\langle\nabla g(a),v\rangle:=\sum_{j=1}^n\frac{\partial g}{\partial \x_j}(a)v_j$ if $v=(v_1,\ldots,v_n)\in L^n$.

If $L=K$, we simplify the notations by setting $\reg_{X,a}:=\reg^L_{X,a}$, $\Reg(X,e):=\Reg^L(X,e)$, $\Reg(X):=\Reg^L(X)$, $\Sing(X):=\Sing^L(X)$ and $T_a(X):=T^L_a(X)$. $\sqbullet$
\end{defn}

\begin{remark}
If $X$ is a $K$-algebraic subset of $L^n$ (and therefore also an algebraic subset of $L^n$), then $T_a(X)=T^L_a(X)\subset T^K_a(X)$ {as} $\II_K(X)\subset\II_L(X)$. $\sqbullet$
\end{remark}

The following $L|K$-Jacobian criterion holds true.

\begin{prop}[{\cite[Thm.5.1.9\;\&\;Prop.5.2.6]{FG}}]\label{jacobian-criterion}
Let $X\subset L^n$ be a $K$-algebraic set and let $a\in X$. Then $a$ is a $K$-nonsingular point of $X$ of dimension $e$ in the sense of Definition \ref{28} if and only if either $e=n$ (so $X=L^n$) or $e<n$ and there exist polynomials $f_1,\ldots,f_{n-e}\in\II_K(X)$ and a Euclidean neighborhood $U$ of $a$ in $L^n$ such that the gradients $\nabla f_1(a),\ldots,\nabla f_{n-e}(a)$ are linearly independent vectors of $L^n$ and $X\cap U=\ZZ_L(f_1,\ldots,f_{n-e})\cap U$.
\end{prop}

\begin{remarks}\label{rem210}
$(\mr{i})$ If $L|K=\R|\R$ or $\C|\C$, the concepts introduced in Definition \ref{28} coincides with the usual ones: if $X$ is an algebraic subset of $L^n$, then the $L$-local ring $\reg_{X,a}=\reg^L_{X,a}$ of~$X$ at $a$ is the usual local ring of $X$ at $a$, a $L$-nonsingular/$L$-singular point of $X$ is a usual nonsingular/singular point of $X$, $\Reg(X)=\Reg^L(X)$ is the usual set of nonsingular points of~$X$, $\Sing(X)=\Sing^L(X)$ is the usual set of singular points of $X$ and the $L$-Zariski tangent space $T_a(X)=T^L_a(X)$ of $X$ at $a$ is the usual Zariski tangent space of $X$ at $a$. Furthermore, the concept of $L$-nonsingular point of $X\subset L^n$ of dimension $e$ coincides with the usual one of nonsingular point of $X\subset L^n$ of dimension $e$, see \cite[p.23]{akbking:tras} and \cite[Def.3.3.9]{BCR}.

$(\mr{ii})$ Also in the case $L|K=\C|\Q$, the concepts introduced in Definition \ref{28} coincides with the usual ones. Let $X$ be a $\Q$-algebraic subset of $\C^n$. Since the ideal $\II_\Q(X)$ of $\Q[\x]$ is radical, the ideal $\II_\Q(X)\otimes_\Q\C=\II_\Q(X)\C[\x]$ of $\C[\x]$ is also radical by \cite[Ch.V,\! \S15,\! Prop.5]{b}. Thus, Hilbert's Nullstellensatz implies that $\II_\C(X)=\II_\Q(X)\C[\x]$ {since} $\II_\C(X)=\II_\C(\ZZ_\C(\II_\Q(X)))=\II_\C(\ZZ_\C(\II_\Q(X)\C[\x]))=\sqrt{\II_\Q(X)\C[\x]}=\II_\Q(X)\C[\x]$. It follows that $\reg_{X,a}=\reg^\Q_{X,a}$, $T_a(X)=T^\Q_a(X)$, $\Reg(X,e)=\Reg^\Q(X,e)$, $\Reg(X)=\Reg^\Q(X)$ and $\Sing(X)=\Sing^\Q(X)$.

$(\mr{iii})$ If $L|K=\R|\Q$, the concepts introduced in Definition \ref{28} do not coincide with the usual ones, they are new: if $X\subset\R^n$ is a $\Q$-algebraic set, then it {may} happen that $\II_\R(X)\supsetneqq\II_\Q(X)\R[\x]$, $\reg_{X,a}$ is not isomorphic to $\reg^\Q_{X,a}$, $\Reg^\Q(X)\subsetneqq\Reg(X)$ and $T_a(X)\subsetneqq T^\Q_a(X)$. For~example, if $X$ is one of the $\Q$-algebraic subsets $C_i=\ZZ_\R(f_i)=\ZZ_\R(g_i)$ of $\R^2$ (for~$i=1,2$) defined in Remark \ref{rem1} and $O\in C_i$ is the origin of $\R^2$, then $\II_\R(C_i)=(f_i)\R[\x]\supsetneqq(g_i)\R[\x]=\II_\Q(C_i)\R[\x]$, $T_O(C_i)=\nabla f_i(O)^\perp=\{0\}\times\R\subsetneqq\R^2=\nabla g_i(O)^\perp=T^\Q_O(C_i)$ and $O\in\Reg(C_i)\setminus\Reg^\Q(C_i)$. The latter condition also implies that the local ring $\reg_{C_i,O}$ is not isomorphic to $\reg^\Q_{C_i,O}$ {since} the former is regular of dimension $1$, the latter is not. See \cite[Ex.5.4.2]{FG} for other examples.

$(\mr{iv})$ By Proposition \ref{jacobian-criterion}, if $L|K=\R|\R$ or $\R|\Q$, the notions of $K$-nonsingular/$K$-singular points of a $K$-algebraic subset of $L^n$ introduced in Definitions \ref{13} and \ref{28} coincide. In addition, if $X\subset\R^n$ is a $\Q$-algebraic set (and therefore also an algebraic subset of $\R^n$), then $\Reg^\Q(X,e)\subset\Reg(X,e)$ for all $e\in\{0,\ldots,\dim(X)\}$. $\sqbullet$
\end{remarks}

\subsubsection{Projective preliminaries}\label{subsub:proj} Let us extend some of the affine notions presented above to the projective case, following \cite[Subsec.2.6\;\&\;Rmk.C.1.2]{FG}.

Let $\PP^n(L)$ be the usual projective $n$-space over $L$, let $[x_0,x]=[x_0,x_1,\ldots,x_n]$ be its homogeneous coordinates, let $K[\x_0,\x]:=K[\x_0,\x_1,\ldots,\x_n]$ and let $K[\x_0,\x]_\sfh$ be the subset of $K[\x_0,\x]$ of homogeneous polynomials. Given any $f\in K[\x_0,\x]_\sfh$, $\PP\ZZ_L(f)$ denotes the zero set of $f$, i.e., $\PP\ZZ_L(f):=\{[x_0,x]\in\PP^n(L): f(x_0,x)=0\}$.

Let $X$ be a subset of $\PP^n(L)$. We say that $X$ is \emph{$K$-algebraic} if $X=\bigcap_{f\in F}\PP\ZZ_L(f)$ for some $F\subset K[\x_0,\x]_\sfh$. We say that $X$ is \emph{$K$-constructible} if it is a Boolean combination of $K$-algebraic subsets of $\PP^n(L)$. The $K$-algebraic subsets of $\PP^n(L)$ are the closed sets of a topology called \emph{$K$-Zariski topology of $\PP^n(L)$}. Observe that an $L$-algebraic subset of $\PP^n(L)$ is a usual algebraic subset of $\PP^n(L)$ and the $L$-Zariski topology of $\PP^n(L)$ is the usual Zariski topology of $\PP^n(L)$. The $K$-Zariski topology of $\PP^n(L)$ is coarser than the usual Zariski topology, so it is Noetherian. As in the affine case, using the $K$-Zariski topology of $\PP^n(L)$, we can introduce the concepts of \emph{$K$-Zariski closed}/\emph{$K$-Zariski open}/\emph{$K$-irreducible} subsets of $\PP^n(L)$ and \emph{$K$-irreducible components} of a $K$-algebraic subset of $\PP^n(L)$, that coincides with the usual ones when $L=K$. The \emph{$K$-Zariski closure of $X$ in $\PP^n(L)$} is the closure of $X$ in $\PP^n(L)$ with respect to the $K$-Zariski topology.

For each $i\in\{0,\ldots,n\}$, let $U^n_i\subset\PP^n(L)$ be the $K$-Zariski open set $\{[x_0,x]\in\PP^n(L):x_i\neq0\}$ and let $\theta^n_i:U^n_i\to L^n$ be the affine chart $\theta^n_i([x_0,x]):=\big(\frac{x_0}{x_i},\ldots,\frac{x_{i-1}}{x_i},\frac{x_{i+1}}{x_i},\ldots,\frac{x_n}{x_i}\big)$. Endow $L^n$ and $\PP^n(L)$ with their $K$-Zariski topologies, and each $U^n_i\subset\PP^n(L)$ with the relative topology. It can be easily verified that each $\theta^n_i$ is a homeomorphism. Since $\{U^n_i\}_{i=0}^n$ is an open cover of $\PP^n(L)$, it follows that the subset $X$ of $\PP^n(L)$ is $K$-Zariski closed/$K$-constructible if and only if $\theta^n_i(U^n_i\cap X)$ is a $K$-Zariski closed/$K$-constructible subset of $L^n$ for all $i\in\{0,\ldots,n\}$.

The \emph{$K$-dimension of $X$} (in $\PP^n(L)$) is defined by $\dim_K(X):=\max_{i\in\{0,\ldots,n\}}\{\dim_K(\theta^n_i(U^n_i\cap X))\}$. If $X\subset\PP^n(L)$ is $K$-algebraic, then $\dim_L(X)=\dim_K(X)$ by Theorem \ref{dimension}. Evidently, the same equality continues to hold for all $K$-constructible sets $X\subset\PP^n(L)$. If $X$ is any $L$-constructible subset of $\PP^n(L)$, we set $\dim(X):=\dim_L(X)$ and we say that \emph{$\dim(X)$ is the dimension of $X$}, so $\dim(X)=\dim_L(X)=\dim_K(X)$ when $X\subset\PP^n(L)$ is $K$-constructible.

Consider now the case of a finite product of projective spaces over $L$. Let $\ell\geq1$, let $n_1,\ldots,n_\ell\in\N^*$ and let $[x^i]=[x^i_0,\ldots,x^i_{n_i}]$ be the homogeneous coordinates of $\PP^{n_i}(L)$. Let $S$ be a subset of $\PP^{n_1}(L)\times\ldots\times\PP^{n_\ell}(L)$. We say that $S$ is \emph{$K$-algebraic} if it is the zero set of a family of polynomials in $K[\x^1,\ldots,\x^\ell]=K[\x^1_0,\ldots,\x^\ell_{n_\ell}]$ that are homogeneous separately in the variables $\x^i=(\x^i_0,\ldots,\x^i_{n_i})$ for all $i\in\{1,\ldots,\ell\}$. We say that $S$ is \emph{$K$-constructible} if it is a Boolean combination of $K$-algebraic subsets of $\PP^{n_1}(L)\times\ldots\times\PP^{n_\ell}(L)$. The \emph{$K$-Zariski topo\-logy of $\PP^{n_1}(L)\times\ldots\times\PP^{n_\ell}(L)$} is the topology whose closed sets are the $K$-algebraic subsets of $\PP^{n_1}(L)\times\ldots\times\PP^{n_\ell}(L)$. As in the case $\ell=1$, using this topology, we can define the concepts of \emph{$K$-Zariski closed}/\emph{$K$-Zariski open}/\emph{$K$-irreducible} subsets of $\PP^{n_1}(L)\times\ldots\times\PP^{n_\ell}(L)$ and \emph{$K$-irreducible components} of a $K$-algebraic subset of $\PP^{n_1}(L)\times\ldots\times\PP^{n_\ell}(L)$.

Endow $L^{n_1+\ldots+n_\ell}$ and $\PP^{n_1}(L)\times\ldots\times\PP^{n_\ell}(L)$ with their $K$-Zariski topologies and,  for each {multi-index} $(j_1,\ldots,j_\ell)$ with $j_i\in\{0,\ldots,n_i\}$, {endow} the set $U^{n_1}_{j_1}\times\ldots\times U^{n_\ell}_{j_\ell}\subset\PP^{n_1}(L)\times\ldots\times\PP^{n_\ell}(L)$ with the relative topology. It is easy to verify that each map $\theta^{n_1}_{j_1}\times\ldots\times\theta^{n_\ell}_{j_\ell}:U^{n_1}_{j_1}\times\ldots\times U^{n_\ell}_{j_\ell}\to L^{n_1+\ldots+n_\ell}$ is a homeomorphism. Therefore, a subset $X$ of $\PP^{n_1}(L)\times\ldots\times\PP^{n_\ell}(L)$ is $K$-Zariski closed if and only if each set $(\theta^{n_1}_{j_1}\times\ldots\times\theta^{n_\ell}_{j_\ell})(X\cap(U^{n_1}_{j_1}\times\ldots\times U^{n_\ell}_{j_\ell}))$ is a $K$-algebraic subset of $L^{n_1+\ldots+n_\ell}$.

The \emph{$K$-dimension $\dim_K(S)$} of the set $S\subset\PP^{n_1}(L)\times\ldots\times\PP^{n_\ell}(L)$ is defined as the maximum over all the multi-indices $(j_1,\ldots,j_\ell)$ with $j_i\in\{0,\ldots,n_i\}$ of the $K$-dimension of the set $(\theta^{n_1}_{j_1}\times\ldots\times\theta^{n_\ell}_{j_\ell})(S\cap(U^{n_1}_{j_1}\times\ldots\times U^{n_\ell}_{j_\ell}))\subset L^{n_1+\ldots+n_\ell}$. As in the case $\ell=1$, if $S\subset\PP^{n_1}(L)\times\ldots\times\PP^{n_\ell}(L)$ is $K$-constructible, then $\dim_L(S)=\dim_K(S)$.

In what follows, we will {freely use} the previous projective notions.


\subsection{Real $\Q$-algebraic sets}

\subsubsection{$\Q$-nonsingular points} Here we present some results related to the concept of $\Q$-nonsingular point introduced in Definition \ref{13}, which is equivalent to that introduced in Definition \ref{28} with $L|K=\R|\Q$. The first is a useful lemma, whose proof is postponed to Appendix \ref{appendix-A}.

\begin{lem}\label{X}
Let $X\subset\R^n$ be an algebraic set of dimension $d<n$, let $a\in\Reg(X,e)$ for some $e\in\{0,\ldots,d\}$ and let $f_1,\ldots,f_{n-e}\in\II_\R(X)$ be such that $\nabla f_1(a),\ldots,\nabla f_{n-e}(a)$ are li\-nearly independent. Then there exists an open neighborhood $U$ of $a$ in $\R^n$ such that $X\cap U=\ZZ_\R(f_1,\ldots,f_{n-e})\cap U$.
\end{lem}

A criterion for detecting $\Q$-nonsingular points among nonsingular points reads as follows.

\begin{prop}\label{Y}
Let $X\subset\R^n$ be a $\Q$-algebraic set of dimension $d$ and let $a\in\Reg(X,e)$ for some $e\in\{0,\ldots,d\}$. Then $a\in\Reg^\Q(X,e)$ if and only if $T_a(X)=T^\Q_a(X)$.
\end{prop}
\begin{proof}
If $e=n$, the statement is evident {as} $X=\Reg(X,e)=\Reg^\Q(X,e)=\R^n$ and $T_a(X)=T^\Q_a(X)=\R^n$. Assume $e<n$. If $a\in\Reg^\Q(X,e)$, then Proposition \ref{jacobian-criterion} assures the existence of polynomials $f_1,\ldots,f_{n-e}\in\II_\Q(X)$ such that $\nabla f_1(a),\ldots,\nabla f_{n-e}(a)$ are linearly independent, so $\dim(T^\Q_a(X))\leq n-(n-e)=e$. Since $a\in\Reg(X,e)$ and $T_a(X)\subset T^\Q_a(X)$, we also have $e=\dim(T_a(X))\leq\dim(T^\Q_a(X))\leq e$, so $\dim(T_a(X))=\dim(T^\Q_a(X))=e$ and $T_a(X)=T^\Q_a(X)$. Conversely, if $T_a(X)=T^\Q_a(X)$, again $\dim(T^\Q_a(X))=e$ so there exist $f_1,\ldots,f_{n-e}\in\II_\Q(X)$ such that $\nabla f_1(a),\ldots,\nabla f_{n-e}(a)$ are linearly independent. Proposition \ref{jacobian-criterion} and Lemma \ref{X} imply that $a\in\Reg^\Q(X,e)$, as required.
\end{proof}

The next result will play an important role in the proof of Theorem \ref{thm:NTQ} presented in Section~\ref{sec:Q-alg-approx}.

\begin{prop}\label{prop:Q_setminus}
If $X,Y\subset\R^n$ are two $\Q$-nonsingular $\Q$-algebraic sets of dimension~$d<n$, such that $Y\subsetneqq X$, then $X\setminus Y\subset\R^n$ is also a $\Q$-nonsingular $\Q$-algebraic set of dimension~$d$.
\end{prop}

To prove this result, we need a preliminary lemma.

\begin{lem}\label{K-difference-preparation}
Let $X\subset\R^n$ be a $\Q$-algebraic set, let $X_1,\ldots,X_r$ be the $\Q$-irreducible components of $X$ and let $a\in\Reg^\Q(X)$. Then there exists an index $i\in\{1,\ldots,r\}$ such that $a\in X_i\setminus\bigcup_{j\in\{1,\ldots,r\}\setminus\{i\}}X_j$.
\end{lem}
\begin{proof}
If $r=1$, the result is evident. Suppose that $r\geq2$ and the statement is false. Rearranging the indices if necessary, we can assume that $a\in X_1\cap X_2$. Since $a\in\Reg^\Q(X)$, by Proposition \ref{jacobian-criterion}, the local ring $\reg^\Q_{X,a}$ is regular so it is an integral domain by \cite[Cor.1, p.302]{zs2}.

Let us show that $\reg^\Q_{X,a}$ is not an integral domain, completing the proof. Let $X':=\bigcup_{\ell=2}^rX_\ell$, let $b_2\in X_2\setminus X_1$, let $b_1\in X_1\setminus X'$ and let $f_1,f_2\in\Q[\x]$ be such that $f_1(b_2)\neq0$, $X_1\subset\ZZ_\R(f_1)$, $f_2(b_1)\neq0$ and $X'\subset\ZZ_\R(f_2)$. Define $\alpha_\ell\in\reg^\Q_{X,a}$ by $\alpha_\ell:=f_\ell+\II_\Q(X)\R[\x]_{\gtn_a}$ for $\ell\in\{1,2\}$. Observe that $\alpha_1\alpha_2=0$ in $\reg^\Q_{X,a}$, so it suffices to prove that $\alpha_1\neq0$ and $\alpha_2\neq0$, i.e., $f_1,f_2\not\in\II_\Q(X)\R[\x]_{\gtn_a}$.

Suppose that $f_1\in\II_\Q(X)\R[\x]_{\gtn_a}$. Let $g_1,\ldots,g_s\in\Q[\x]$ be generators of $\II_\Q(X)$ in $\Q[\x]$. Since $f_1\in\II_\Q(X)\R[\x]_{\gtn_a}$, there exist $h_0,\ldots,h_s\in\R[\x]$ such that $h_0(a)\neq0$ and $h_0f_1=\sum_{k=1}^sh_kg_k$ in $\R[\x]$. Write each $h_k$ explicitly: $h_k=\sum_{\beta\in\N^n}c_{k,\beta}\x^\beta$, where only finitely many of the real coefficients $c_{k,\beta}$ are non-zero. Let $\{u_p\}_{p\in P}$ be a basis of $\R$ as a $\Q$-vector space. Write each real coefficient $c_{k,\beta}$ as follows: $c_{k,\beta}=\sum_{p\in P}c_{k,\beta,p}u_p$, where $c_{k,\beta,p}\in\Q$ and only finitely many of the $c_{k,\beta,p}$ are non-zero. For every $(k,p)\in\{0,\ldots,s\}\times P$, define the polynomial $h_{k,p}\in\Q[\x]$ by $h_{k,p}:=\sum_{\beta\in\N^n}c_{k,\beta,p}\x^\beta$. Observe that only finitely many of the $h_{k,p}$ are non-zero. It follows that $h_k=\sum_{p\in P}h_{k,p}u_p$ for every $k\in\{0,\ldots,s\}$, since $h_k=\sum_{\beta\in\N^n}c_{k,\beta}\x^\beta=\sum_{\beta\in\N^n}\sum_{p\in P}c_{k,\beta,p}u_p\x^\beta
=\sum_{p\in P}h_{k,p}u_p$. We deduce:
$$\textstyle
\sum_{p\in P}(h_{0,p}f_1)u_p=h_0f_1=\sum_{k=1}^sh_kg_k=\sum_{k=1}^s(\sum_{p\in P}h_{k,p}u_p)g_k=\sum_{p\in P}(\sum_{k=1}^sh_{k,p}g_k)u_p.
$$
Observe that the polynomials $\{h_{0,p}f_1\}_{p\in P}$ and $\{\sum_{k=1}^sh_{k,p}g_k\}_{p\in P}$ belongs to $\Q[\x]$, and only finitely many of these polynomials are non-zero. Since $\{u_p\}_{p\in P}$ is a $\Q$-vector basis of $\R$, it follows that $h_{0,p}f_1=\sum_{k=1}^sh_{k,p}g_k$ for all $p\in P$. Since $\sum_{p\in P}h_{0,p}(a)u_p=h_0(a)\neq0$, there exists $q\in P$ such that $h_{0,q}(a)\neq0$. Set $h:=h_{0,q}$. Observe that $h\in\Q[\x]$, $h(a)\neq0$ and $hf_1=\sum_{k=1}^sh_{k,q}g_k\in\II_\Q(X)$.

Endow $X_2$ with the relative topology induced by the $\Q$-Zariski topology of $\R^n$, making $X_2$ an irreducible topological space. Define the open subset $U$ of $X_2$ by $U:=X_2\setminus\ZZ_\R(f_1)$. Observe that $U$ is non-empty because $b_2\in U$ by construction, so $U$ is dense in $X_2$. Since $h\in\Q[\x]$ and $hf_1\in\II_\Q(X)\subset\II_\Q(X_2)$, the intersection $\ZZ_\R(h)\cap X_2$ is a closed subset of $X_2$ containing $U$. By the density of $U$ in $X_2$, we deduce that $X_2\subset\ZZ_\R(h)$, which is impossible {since} $a\in X_2$ and $h(a)\neq0$. This proves that $f_1\not\in\II_\Q(X)\R[\x]_{\gtn_a}$.

A similar argument shows that $f_2\not\in\II_\Q(X)\R[\x]_{\gtn_a}$. This completes the proof.
\end{proof}

Lemma \ref{K-difference-preparation} has been improved in \cite[Lem.5.2.1]{FG}.

\begin{proof}[Proof of Proposition \ref{prop:Q_setminus}]
We can assume $Y\neq\varnothing$, otherwise the result is evident. Denote $X_1,\ldots,X_r$ the $\Q$-irreducible components of $X\subset\R^n$. By Lemma \ref{K-difference-preparation}, the sets $X_i$ are pairwise disjoint {as} $X=\Reg^\Q(X)$ by hypothesis. Since $X\subset\R^n$ is a nonsingular algebraic set of dimension $d$ and the sets $X_i$ are closed in $\R^n$, by \cite[Props.3.3.10\;\&\;3.3.11]{BCR}, we deduce that each $X_i$ is a Nash submanifold of $\R^n$ of dimension $d$, so $\dim_\R(X_i)=d$. Let $Y_i:=Y\cap X_i$ for every $i\in\{1,\ldots,r\}$ and let $I$ be the (non-empty) set of all indices $i\in\{1,\ldots,r\}$ such that $Y_i\neq\varnothing$.

Fix $i\in I$ and set $X'_i:=\bigcup_{j\in\{1,\ldots,r\}\setminus\{i\}}X_j$. Since $X_i\cap X'_i=\varnothing$ and $Y_i\subset X_i$, we deduce $Y_i\cap X'_i=\varnothing$ so $Y_i=Y\cap(\R^n\setminus X'_i)$. Using again \cite[Props.3.3.10\;\&\;3.3.11]{BCR}, we have that $Y$ is a Nash submanifold of $\R^n$ of dimension $d$. Since $\R^n\setminus X'_i$ is a semialgebraic open subset of $\R^n$ and $Y_i=Y\cap(\R^n\setminus X'_i)\neq\varnothing$, it follows that $Y_i$ is also a Nash submanifold of $\R^n$ of dimension~$d$, so $\dim_\R(Y_i)=d$. Since $X_i$ and $Y_i$ are $\Q$-algebraic subsets of $\R^n$, by Theorem \ref{dimension}, we deduce $\dim_\Q(X_i)=\dim_\Q(Y_i)$ {since} $\dim_\Q(X_i)=\dim_\R(X_i)=d=\dim_\R(Y_i)=\dim_\Q(Y_i)$. Now Lemma \ref{dimirred} assures that $X_i=Y_i$. This proves that $Y=\bigsqcup_{i\in I}X_i$. Since $Y\subsetneqq X$, we deduce that $I\subsetneqq\{1,\ldots,r\}$, so $X\setminus Y=\bigsqcup_{j\in \{1,\ldots,r\}\setminus I}X_j$ is a $\Q$-nonsingular $\Q$-algebraic subset of $\R^n$ of dimension $d$, as required.
\end{proof}

Proposition \ref{prop:Q_setminus} has been improved in \cite[Cor.5.2.9]{FG}.

Using, among other results, the above-mentioned improved version of Lemma \ref{K-difference-preparation} obtained in \cite[Lem.5.2.1]{FG}, one can prove the following structure theorem for the $\Q$-nonsingular and $\Q$-singular loci of a $\Q$-algebraic subset of~$\R^n$.

\begin{prop}[{\cite[Prop.5.2.6\;\&\;Cor.5.2.8]{FG}}]\label{prop211}
Let $X\subset\R^n$ be a $\Q$-algebraic set. We have:
\begin{itemize}
 \item[$(\mr{i})$] A point $a$ of $X$ is $\Q$-nonsingular of dimension $e$ if and only if $a$ belongs to only one $\Q$-irreducible component $X'$ of $X$ such that $\dim(X')=e$ and $a\in\Reg^\Q(X')$.  
 \item[$(\mr{ii})$] $\Sing^\Q(X)\subset\R^n$ is a $\Q$-algebraic set of dimension $<\dim(X)$. In particular, $\Reg^\Q(X)$ is a non-empty $\Q$-Zariski open subset of $X$.
\end{itemize}
\end{prop}


\subsubsection{$\Q$-regular maps} \label{subsec:Q-regular-maps}
Let $K:=\R$ or $\Q$. Recall that $n$ denotes a positive natural number and $K[\x]$ the ring of polynomials $K[\x_1,\ldots,\x_n]$. 

The next definition is a reformulation of \cite[Def.4, p.30]{togn:instmat} (see also \cite[Def.4.3.1]{FG}).

\begin{defn}\label{def:Q-reg-function}
Let $X$ be a subset of $\R^n$, let $Y$ be a subset of $\R^m$, let $f:X\to Y$ be a map and let $a$ be a point of $X$. We say that $f$ is \emph{$K$-regular at $a$} if there exist a $K$-Zariski open neighborhood $U$ of $a$ in $\R^n$ and polynomials $p_1,\ldots,p_m,q\in K[\x]$ such that $q(x)\neq0$ and $f(x)=\big(\frac{p_1(x)}{q(x)},\ldots,\frac{p_m(x)}{q(x)}\big)$ for each $x\in X\cap U$. We say that $f$ is \emph{$K$-regular} if $f$ is $K$-regular at each point of $X$. We denote $\reg^K(X,Y)$ the set of $K$-regular maps from $X$ to $Y$, and $\reg^K(X)$ the set $\reg^K(X,\R)$ of $K$-regular functions on $X$. We say that $f:X\to Y$ is a \emph{$K$-biregular isomorphism} if $f$ is bijective and both $f$ and $f^{-1}$ are $K$-regular.

The same definition can be rephrased in the complex case by simply replacing $\R$ {by} $\C$. $\sqbullet$
\end{defn}

An $\R$-regular function on $X$ is a usual regular function on $X$, $\reg^\R(X)$ is the usual set $\reg(X)$ of regular functions on $X$ and $\reg^\R(X,Y)$ is the usual set $\reg(X,Y)$ of regular maps from $X$ to $Y$, see \cite[Sect.3.2]{BCR} and \cite[p.19]{akbking:tras}.

In the next lemma we collect some elementary but very useful properties of $\Q$-regular maps.

\begin{lem}\label{lem:Q-basic}
Let $X\subset\R^n$  and $Y\subset\R^m$ be sets and let $f:X\to Y$ be a map. We have:
\begin{itemize}
 \item[$(\mr{i})$] The map $f$ is $\Q$-regular if and only if there exist polynomials $p_1,\ldots, p_m,q\in\Q[\x]$ such that $q(x)\neq0$ and $f(x)=\big(\frac{p_1(x)}{q(x)},\ldots,\frac{p_m(x)}{q(x)}\big)$ for all $x\in X$.
 \item[$(\mr{ii})$] Suppose that $f$ is $\Q$-regular. If $Z\subset\R^\ell$ is a set and $g:Y\to Z$ is a $\Q$-regular map, then the composition $g\circ f:X\to Z$ is also a $\Q$-regular map.
\end{itemize}

In the remaining part of the statement, we assume that $X\subset\R^n$ is $\Q$-algebraic and $f\in\reg^\Q(X,Y)$. We have:
\begin{itemize}
 \item[$(\mr{iii})$] There exists $F\in\reg^\Q(\R^n,\R^m)$ such that $F(x)=f(x)$ for all $x\in X$.
 \item[$(\mr{iv})$] Let $X'$ be another $\Q$-algebraic subset of $\R^n$ with $X\cap X'=\varnothing$, let $f'\in\reg^\Q(X',Y)$ and let $f\sqcup f':X\sqcup X'\to Y$ be the disjoint union map of $f$ and $f'$ defined by $(f\sqcup f')(x):=f(x)$ for all $x\in X$ and $(f\sqcup f')(x):=f'(x)$ for all $x\in X'$. Then $f\sqcup f'\in\reg^\Q(X\sqcup X',Y)$.
 \item[$(\mr{v})$] If $W$ is a $\Q$-Zariski closed subset of $Y$, then the set $f^{-1}(W)\subset\R^n$ is $\Q$-algebraic. 
\end{itemize}
\end{lem}

An elementary proof of this result can be obtained by combining standard arguments with the Noetherianity of the $\Q$-Zariski topology and the fact that $\Q$ is a subfield of $\R$: some care is needed. For completeness, we include this proof in Appendix~\ref{appendix-A}.

In what follows, the symbol $\mr{rk}(A)$ denotes the rank of a square matrix $A$. If $V$ is a vector subspace of $\R^n$ and $v_1,\ldots,v_\ell$ are vectors of $\R^n$, then we denote $V^\perp$ the vector subspace of $\R^n$ that is orthogonal to $V$ with respect to the usual Euclidean scalar product, and $\mr{Span}(v_1,\ldots,v_\ell)$ the vector subspace of $\R^n$ generated by $v_1,\ldots,v_\ell$.

The next result deals with transversality of $\Q$-regular maps. 

\begin{prop}\label{prop:Q-transverse}
Let $X\subset\R^n$, $Y\subset\R^m$ and $W\subset\R^m$ be $\Q$-algebraic sets, let $f:X\to Y$ be a $\Q$-regular map and let $a\in X$ be such that $W\subset Y$, $a$ is a $\Q$-nonsingular point of $X$ of dimension~$d$, $f(a)$ is both a $\Q$-nonsingular point of $Y$ of dimension $e$ and a $\Q$-nonsingular point of $W$ of dimension $e'<e$, and $f$ is transverse to $W$ in~$Y$ at $a$, i.e., $d_af(T_a(X))+T_{f(a)}(W)=T_{f(a)}(Y)$. Then $a$ is a $\Q$-nonsingular point of the $\Q$-algebraic set $f^{-1}(W)\subset\R^n$ of dimension $d-(e-e')$.
\end{prop}
\begin{proof}
Since $X\subset\R^n$, $Y\subset\R^m$ and $W\subset\R^m$ are algebraic sets, $f$ is regular, $a$ is a nonsingular point of $X$ of dimension~$d$, $f(a)$ is both a nonsingular point of $Y$ of dimension $e$ and a nonsingular point of $W$ of dimension $e'<e$ in the usual sense, a standard transversality result (see \cite[Lem.2.2.13]{akbking:tras}) assures that $a$ is a nonsingular point of the algebraic set $f^{-1}(W)\subset\R^n$ of dimension $d-(e-e')$. In addition, Lemma \ref{lem:Q-basic}$(\mr{v})$ implies that $f^{-1}(W)\subset\R^n$ is $\Q$-algebraic.

Let $b:=f(a)$. Suppose for a moment that $d<n$ and $e<m$. Since $a$ is a $\Q$-nonsingular point of $X$ of dimension~$d$ and $b$ is a $\Q$-nonsingular point of $Y$ of dimension $e$, there exist $f_1,\ldots,f_{n-d}\in\II_\Q(X)$, an open neighborhood $U$ of $a$ in $\R^n$, $g_1,\ldots,g_{m-e}\in\II_\Q(Y)$ and an open neighborhood $V$ of $b$ in $\R^m$ such that $\nabla f_1(a),\ldots,\nabla f_{n-d}(a)$ are linearly independent in~$\R^n$, $\nabla g_1(b),\ldots,\nabla g_{m-e}(b)$ are linearly independent in $\R^m$, $X\cap U=\ZZ_\R(f_1,\ldots,f_{n-d})\cap U$ and $Y\cap V=\ZZ_\R(g_1,\ldots,g_{m-e})\cap V$. By Proposition \ref{Y}, we know that $T^\Q_b(Y)=T_b(Y)$ and $T^\Q_b(W)=T_b(W)$ so $\dim(T^\Q_b(Y))=e$ and $\dim(T^\Q_b(W))=e'$. Thus, there exist $g_{m-e+1},\ldots,g_{m-e'}\in\II_\Q(W)$ such that 
$\nabla g_1(b),\ldots,\nabla g_{m-e'}(b)$ are linearly independent in $\R^m$. By Lemma \ref{X}, shrinking $V$ around $b$ if necessary, we can assume that $W\cap V=\ZZ_\R(g_1,\ldots,g_{m-e'})\cap V$. By Lemma \ref{lem:Q-basic}$(\mr{i})(\mr{iii})$, there exist $p_1,\ldots,p_m,q\in\Q[\x]$ such that $\ZZ_\R(q)=\varnothing$ and $f(x)=\big(\frac{p_1(x)}{q(x)},\ldots,\frac{p_m(x)}{q(x)}\big)$ for all $x\in X$.

Define $F:\R^n\to\R^m$ by $F(x):=\big(\frac{p_1(x)}{q(x)},\ldots,\frac{p_m(x)}{q(x)}\big)$ for all $x\in\R^n$. Shrinking $U$ around~$a$ if necessary, we can assume that $F(U)\subset V$. Let $i\in\{1,\ldots,e-e'\}$, let $c_i:=\deg(g_{m-e+i})$ and let $h_i\in\Q[\x]$ be the polynomial such that $h_i(x)=q(x)^{c_i}g_{m-e+i}(F(x))$ for all $x\in\R^n$. Each $h_i$ belongs to $\II_\Q(f^{-1}(W))$ since $h_i(x)=q(x)^{c_i}g_{m-e+i}(f(x))=0$ for all $x\in f^{-1}(W)$. Define the $\Q$-regular map $G:\R^n\to\R^{n-d+e-e'}$ by $G:=(f_1,\ldots,f_{n-d},h_1,\ldots,h_{e-e'})$. We have:
\begin{align*}
f^{-1}(W)\cap U&=\{x\in X\cap U:f(x)\in W\}=\{x\in X\cap U:f(x)\in W\cap V\}\nonumber\\
&=\{x\in X\cap U:f(x)\in\ZZ_\R(g_1,\ldots,g_{m-e},g_{m-e+1},\ldots,g_{m-e'})\}\nonumber\\
&=(X\cap U)\cap\{x\in X:g_{m-e+1}(f(x))=0,\ldots,g_{m-e'}(f(x))=0\}\nonumber\\
&=\ZZ_\R(f_1,\ldots,f_{n-d})\cap\{x\in U:g_{m-e+1}(F(x))=0,\ldots,g_{m-e'}(F(x))=0\}\nonumber\\
&=G^{-1}(0)\cap U.\label{1}
\end{align*}
Since $f$ is transverse to $W$ in $Y$ at $a$, we have $d\geq e-e'$ and $\dim((d_af)^{-1}(T_b(W)))=d-(e-e')$. Let $J_G(a)$ be the Jacobian matrix of $G$ at $a$ and let $G_i$ be the $i^{\mr{th}}$-row of $J_G(a)$. Observe that $G_i=\nabla f_i(a)$ if $i\in\{1,\ldots,n-d\}$ and $G_{n-d+i}=q(a)^{c_i}\nabla g_{m-e+i}(b)J_G(a)$ if $i\in\{1,\ldots,e-e'\}$. Thus, a vector $v\in\R^n$ belongs to $\ker(J_G(a))$ if and only if $v\in\mr{Span}(\nabla f_1(a),\ldots,\nabla f_{n-d}(a))^\perp=T_a(X)$ and $d_af(v)=J_G(a)v\in\mr{Span}(\nabla g_{m-e+1}(b),\ldots,\nabla g_{m-e'}(b))^\perp\cap T_b(Y)=T_b(W)$. This proves that $\ker(J_G(a))=(d_af)^{-1}(T_b(W))$ so $\dim(\ker(J_G(a)))=d-(e-e')$ or, equivalently, $\mr{rk}(J_G(a))=n-d+e-e'$. The latter equality and the fact that $a$ is a nonsingular point of $f^{-1}(W)$ of dimension $d-(e-e')$ imply that $\nabla f_1(a),\ldots,\nabla f_{n-d}(a),\nabla h_1(a),\ldots,\nabla h_{e-e'}(a)$ are linearly independent and $f^{-1}(W)\cap U=\ZZ_\R(f_1,\ldots,f_{n-d},h_1,\ldots,h_{e-e'})\cap U$. Since the polynomials $f_1,\ldots,f_{n-d},h_1,\ldots,h_{e-e'}$ belong to $\II_\Q(f^{-1}(W))$, Proposition \ref{jacobian-criterion} assures that $a$ is a $\Q$-nonsingular point of $f^{-1}(W)$ of dimension $d-(e-e')$, as required. Finally, if $d=n$ and/or $e=m$, then $X=\R^n$ and/or $Y=\R^m$ and the previous argument continues to work, provided we omit the polynomials $f_1,\ldots,f_{n-d}$ and/or $g_1,\ldots,g_{m-e}$.
\end{proof}

The following lemma assures that $\Q$-biregular isomorphisms preserve nonsingular points.

\begin{lem}\label{861}
Let $X\subset\R^n$ and $Y\subset\R^m$ be $\Q$-algebraic sets of dimension $d$, let $e\in\{0,\ldots,d\}$, let $X'$ be a $\Q$-Zariski open subset of $X$, let $Y'$ be a $\Q$-Zariski open subset of $Y$, let $f:X'\to Y'$ be a $\Q$-biregular isomorphism and let $a\in X'$. Then $a\in\Reg^\Q(X,e)$ if and only if $f(a)\in\Reg^\Q(Y,e)$.  
\end{lem}
\begin{proof}
Let $b:=f(a)$ and let $f^*:\reg^\Q_{Y,b}\to\reg^\Q_{X,a}$ be the pullback homomorphism induced by $f$ on $\Q$-local rings in the natural way, see definition (4.7) in \cite[Subsec.4.3]{FG}. It is immediate to verify that $f^*$ is a ring isomorphism whose inverse is $(f^{-1})^*$. It follows that $\reg^\Q_{X,a}$ is a regular local ring of dimension $e$ if and only if $\reg^\Q_{Y,b}$ is {so}. 
\end{proof}

\subsubsection{Projectively $\Q$-closed $\Q$-algebraic sets}
Let $H:=\{[x_0,x]\in\PP^n(\R):x_0=0\}$. Recall that $\theta_0^n:U^n_0\to\R^n$ denotes the affine chart $\theta^n_0([x_0,x])=\big(\frac{x_1}{x_0},\ldots,\frac{x_n}{x_0}\big)$, where $U^n_0=\PP^n(\R)\setminus H$. Set $\theta:=\theta^n_0$ for short and observe that $\theta^{-1}(x):=[1,x]$. According to \cite[p.427]{ak1981}, a set $X\subset\R^n$ is projectively closed if $\theta^{-1}(X)$ is Zariski closed in $\PP^n(\R)$. By Definition \ref{def:proj-Q-closed}, $X$ is projectively $\Q$-closed if $\theta^{-1}(X)$ is $\Q$-Zariski closed in $\PP^n(\R)$. Clearly, if $X$ is projectively $\Q$-closed, then it is also projectively closed and compact in $\R^n$. Let us recall the concept of overt polynomial.

\begin{defn}[{\cite[p.427]{ak1981}}]
A polynomial $f\in\R[\x]$ is said to be \emph{overt} if either $f$ is a non-zero constant or $f$ is {non-constant} and, if we write $f=\sum_{\ell=0}^ef_\ell$, where $e=\deg(f)$ and each $f_\ell$ is a homogeneous polynomial of degree $\ell$, then the leading homogeneous term $f_e$ of $f$ {only vanishes} at $0$, i.e., $\ZZ_\R(f_e)=\{0\}$. $\sqbullet$
\end{defn}

If $f\in\R[\x]$, we denote $f^\sfh(\x_0,\x)\in\R[\x_0,\x]_\sfh$ the homogenization of $f$. Recall that $f^\sfh:=f$ if $f$ is constant. If $f$ is {non-constant} and we write $f=\sum_{\ell=0}^ef_\ell$ as above, then $f^\sfh(\x_0,\x):=\sum_{\ell=0}^e\x_0^{e-\ell}f_\ell(\x)$. Observe that $\theta^{-1}(\ZZ_\R(f))=\PP\ZZ_\R(f^\sfh)\setminus H$.

\begin{lem}\label{lem:overt}
A $\Q$-algebraic set $X\subset\R^n$ is projectively $\Q$-closed if and only if there exists an overt polynomial $f\in\Q[\x]$ such that $X=\ZZ_\R(f)$. 
\end{lem}

Before proving this lemma we present an elementary but useful result.

\begin{lem}\label{lem:Z}
If $X\subset\PP^n(\R)$ is a $\Q$-algebraic set, then there exists $f\in\Q[\x_0,\x]_\sfh$ such that $\PP\ZZ_\R(f)=X$.
\end{lem}
\begin{proof}
If $X=\varnothing$, it suffices to set $f:=1$. Assume $X\neq\varnothing$. Since the $\Q$-Zariski topology of $\PP^n(\R)$ is Noetherian, there exist finitely many non-constant homogeneous polynomials $f_1,\ldots,f_s\in\Q[\x_0,\x]_\sfh$ such that $X=\bigcap_{i=1}^s\PP\ZZ_\R(f_i)$. Let $e_i:=\deg(f_i)\in\N^*$, let $e:=\prod_{i=1}^se_i\in\N^*$, let $e'_i:=\frac{e}{e_i}\in\N^*$ and let $f:=\sum_{i=1}^sf_i^{2e'_i}\in\Q[\x_0,\x]_\sfh$. Evidently, it holds $\PP\ZZ_\R(f)=X$.
\end{proof}

\begin{proof}[Proof of Lemma \ref{lem:overt}]
If $X=\ZZ_\R(f)$ for some overt polynomial $f\in\Q[\x]$, then $X$ is projectively $\Q$-closed {as} $\theta^{-1}(X)=\PP\ZZ_\R(f^\sfh)$.

Suppose that $X$ is projectively $\Q$-closed. If $X=\varnothing$, it suffices to set $f:=1$. Assume $X\neq\varnothing$. By Lemma \ref{lem:Z}, there exists a non-constant polynomial $g\in\Q[\x_0,\x]_\sfh$ such that $\theta^{-1}(X)=\PP\ZZ_\R(g)$. Write $g(\x_0,\x)=\sum_{\ell=0}^e\x_0^{e-\ell}g_\ell(\x)$, where $e$ is the positive degree of $g$ and each $g_\ell$ is {a} homogeneous polynomial in $\Q[\x]_\sfh$ of degree $\ell$. Define $f(\x):=g(1,\x)=\sum_{\ell=0}^eg_\ell(\x)\in\Q[\x]$. Observe that $X=\ZZ_\R(f)$, $g_e$ is the leading term of $f$ and $\PP\ZZ_\R(g_e)=H\cap\theta^{-1}(X)=\varnothing$, i.e., $\ZZ_\R(g_e)=\{0\}$. This proves that $f$ has the required properties.
\end{proof}

The next lemma collects some elementary properties of projectively $\Q$-closed $\Q$-algebraic sets, which we will use later.

\begin{lem}\label{lem:projective}
Let $X\subset\R^n$ be a projectively $\Q$-closed $\Q$-algebraic set. We have:
\begin{itemize}
 \item[$(\mr{i})$] If $X'\subset\R^n$ is another projectively $\Q$-closed $\Q$-algebraic set, then $X\cup X'\subset\R^n$ is also projectively $\Q$-closed.
 \item[$(\mr{ii})$] If $Z\subset\R^n$ is a $\Q$-algebraic set contained in $X$, then $Z\subset\R^n$ is projectively $\Q$-closed.
 \item[$(\mr{iii})$] For every $v\in\Q^n$, the translated set $X+v:=\{x+v\in\R^n:x\in X\}\subset\R^n$ is a projectively $\Q$-closed $\Q$-algebraic set. In addition, $X+v\subset\R^n$ is $\Q$-nonsingular if $X$ is {so}.
 \item[$(\mr{iv})$] If $Y\subset\R^m$ is a projectively $\Q$-closed $\Q$-algebraic set, then $X\times Y\subset\R^{n+m}$ is also a projectively $\Q$-closed $\Q$-algebraic set. In particular, $X=X\times\{0\}$ is also a projectively $\Q$-closed $\Q$-algebraic subset of $\R^{n+m}=\R^n\times\R^m$. If in addition both $X\subset\R^n$ and $Y\subset\R^m$ are $\Q$-nonsingular, then $X\times Y\subset\R^{n+m}$ is also $\Q$-nonsingular.
\end{itemize}
\end{lem}
\begin{proof}
$(\mr{i})$ If $\theta^{-1}(X)$ and $\theta^{-1}(X')$ are $\Q$-algebraic subsets of $\PP^n(\R)$, their union $\theta^{-1}(X\cup X')$ {is so}.

$(\mr{ii})$ Choose $g\in\Q[\x]$ such that $\ZZ_\R(g)=Z$ and define $g^*\in\Q[\x_0,\x]_\sfh$ by $g^*:=\x_0g^\sfh(\x_0,\x)$. Observe that $\PP\ZZ_\R(g^*)=H\cup\theta^{-1}(Z)$ and thus $\theta^{-1}(Z)=\PP\ZZ_\R(g^*)\cap\theta^{-1}(X)$. Since $\theta^{-1}(X)\subset\PP^n(\R)$ is $\Q$-algebraic by hypothesis, we deduce that $\theta^{-1}(Z)\subset\PP^n(\R)$ is also $\Q$-algebraic.

$(\mr{iii})$ Choose $v\in\Q^n$. If $X=\varnothing$, then $X+v=\varnothing$ and the statement is evident. Suppose that  $X\neq\varnothing$. By Lemma \ref{lem:overt}, there exists an overt polynomial $f\in\Q[\x]$ of positive degree such that $\ZZ_\R(f)=X$. 
Let $F$ be the unique polynomial in $\Q[\x]$ such that $F(x)=f(x-v)$ for all $x\in\R^n$. Observe that $f$ and $F$ have the same degree and the same leading homogeneous term. In particular, $F$ is overt. Since $\ZZ_\R(F)=X+v$, the set $X+v\subset\R^n$ is $\Q$-algebraic and Lemma \ref{lem:overt} assures that the set $X+v\subset\R^n$ is also projectively $\Q$-closed. The last assertion follows immediately from Definition \ref{13} of $\Q$-nonsingular point.

$(\mr{iv})$ If either $X=\varnothing$ or $Y=\varnothing$, then $X\times Y=\varnothing$ and the statement is evident. Suppose $X$ and $Y$ are both non-empty. By Lemma \ref{lem:overt}, there exist non-constant overt polynomials $f\in\Q[\x]$ and $g\in\Q[\y]:=\Q[\y_1,\ldots,\y_m]$ such that $\ZZ_\R(f)=X$ and $\ZZ_\R(g)=Y$. Define $e:=\deg(f)>0$, $d:=\deg(g)>0$ and $u\in\Q[\x,\y]$ by $u(\x,\y):=f(\x)^{2d}+g(\y)^{2e}$. Observe that $u$ is overt and $\ZZ_\R(u)=X\times Y$. Using Lemma \ref{lem:overt} again, we deduce that $X\times Y\subset\R^{n+m}$ is a projectively $\Q$-closed $\Q$-algebraic set. If we set $g:=\sum_{j=1}^m\y_j^2\in\Q[\y]$, then $Y=\ZZ_\R(g)=\{0\}$ and thus $X=X\times\{0\}\subset\R^{n+m}$ is projectively $\Q$-closed. Finally, if in addition both $X\subset\R^n$ and $Y\subset\R^m$ are $\Q$-nonsingular, then the fact that $X\times Y\subset\R^{n+m}$ is also $\Q$-nonsingular follows easily from Definition \ref{13}.
\end{proof}

\subsubsection{An $\R|\Q$-generic projection theorem}
Let $n,r\in\N^*$ with $r<n$. Write $\R^n=\R^r\times\R^{n-r}$ and denote $x'$ the coordinates $(x_1,\ldots,x_r)$ of $\R^r$ and $x''$ the coordinates $(x_{r+1},\ldots,x_n)$ of $\R^{n-r}$, so $x=(x',x'')=(x_1,\ldots,x_n)$ are the coordinates of $\R^n$. For convenience, we consider $x'$ and $x''$ as column vectors. Denote $\mc{M}_{r,n-r}(\Q)$ the $\Q$-vector space of all $r\times(n-r)$-matrices with coefficients in~$\Q$. Identify $\mc{M}_{r,n-r}(\Q)$ with $\Q^{r(n-r)}$ and endow $\mc{M}_{r,n-r}(\Q)=\Q^{r(n-r)}$ with the usual Zariski topo\-logy. We also interpret each element $A$ of $\mc{M}_{r,n-r}(\Q)$ as a $r\times(n-r)$-matrix with real coef\-ficients, so the matrix product $Ax''$ is a well-defined vector of $\R^r$ for every $x''\in\R^{n-r}$. As in \cite[Def.6.3.1]{FG}, for each $A\in\mc{M}_{r,n-r}(\Q)$, we define the $\R$-vector subspace $V_A$ of $\R^n$ by $V_A:=\{(x',x'')\in\R^n:x'=Ax''\}$ and the corresponding $\R$-linear projection $\pi_A:\R^n\to\R^r$ of $\R^n$ onto $\R^r$ in the direction of $V_A$ by $\pi_A(x',x''):=x'-Ax''$.

Our next result reads as follows: {it} is a $\R|\Q$-version of the {classical} generic projection theorem.

\begin{thm}\label{thm:generic-projection}
Let $X\subset\R^n$ be a projectively $\Q$-closed $\Q$-nonsingular $\Q$-algebraic set of dimension~$d$ and let $r:=2d+1$. If $r<n$, then there exists a non-empty Zariski open subset $\Omega$ of $\mc{M}_{r,n-r}(\Q)$ such that, for every $A\in\Omega$, we have:
\begin{itemize}
 \item[$(\mr{i})$] $\pi_A(X)\subset \R^r$ is a projectively $\Q$-closed $\Q$-nonsingular $\Q$-algebraic set.
 \item[$(\mr{ii})$] The restriction of $\pi_A$ from $X$ to $\pi_A(X)$ is a $\Q$-biregular isomorphism.
\end{itemize}
\end{thm}

This statement without the prefix `$\Q$-' is a version of the classical generic projection theorem. The same statement, stripped of the words `projectively $\Q$-closed', coincides with that of \cite[Cor.6.3.5]{FG} when $L|K=\R|\Q$. Therefore, to prove Theorem \ref{thm:generic-projection}, it suffices to show that the next result holds.

\begin{prop}\label{213}
Let $X\subset\R^n$ be a projectively $\Q$-closed $\Q$-algebraic set of dimension~$d$ and let $r:=2d+1$. If $r<n$, then there exists a non-empty Zariski open subset $\Omega$ of $\mc{M}_{r,n-r}(\Q)$ such that $\pi_A(X)\subset \R^r$ is a projectively $\Q$-closed $\Q$-algebraic set for every $A\in\Omega$.
\end{prop}

Suppose for a moment that the previous result is true.

\begin{proof}[Proof of Theorem \ref{thm:generic-projection}]
By \cite[Cor.6.3.5]{FG}, there exists a non-empty Zariski open subset $\Omega_1$ of $\mc{M}_{r,n-r}(\Q)$ such that, for every $A\in\Omega_1$,  $\pi_A(X)\subset \R^r$ is a $\Q$-nonsingular $\Q$-algebraic set and the restriction of $\pi_A$ from $X$ to $\pi_A(X)$ is a $\Q$-biregular isomorphism. By Proposition \ref{213}, there exists a non-empty Zariski open subset $\Omega_2$ of $\mc{M}_{r,n-r}(\Q)$ such that $\pi_A(X)\subset \R^r$ is a projectively $\Q$-closed
$\Q$-algebraic set for every $A\in\Omega_2$. Since $\mc{M}_{r,n-r}(\Q)=\Q^{r(n-r)}$ (equipped with the Zariski topology) is an irreducible topological space, both $\Omega_1$ and $\Omega_2$ are dense in $\mc{M}_{r,n-r}(\Q)$. Thus, their intersection $\Omega:=\Omega_1\cap\Omega_2$ is a non-empty Zariski open subset of $\mc{M}_{r,n-r}(\Q)$ with the required properties.
\end{proof}

It remains to prove Proposition \ref{213}. We need some {preparation}. First, we recall a $\C|\Q$-elimination theorem proved in \cite{FG}.

\begin{thm}[{\cite[Thm.2.6.10 \& Rmk.C.1.2]{FG}}]\label{thm:K-elimination}
Let $n_1,\ldots,n_\ell$ be positive natural numbers and let 
$\rho:\PP^{n_1}(\C)\times\ldots\times\PP^{n_\ell}(\C)\to\PP^{n_\ell}(\C)$ be the canonical projection onto the last factor. We~have:
\begin{itemize}
 \item[$(\mr{i})$] $\rho$ is a $\Q$-Zariski closed map, i.e., it maps $\Q$-algebraic sets in $\Q$-algebraic sets.
 \item[$(\mr{ii})$] $\rho$ maps $\Q$-constructible sets in $\Q$-constructible sets.
\end{itemize}
\end{thm}

Secondly, we recall the notion of $\Q$-regular map in the complex projective setting.  

\begin{defn}[{\cite[Def.6.3.11]{FG}}] \label{def:331-projection}
Let $X$ be a subset of $\PP^n(\C)$, let $Y$ be a subset of $\PP^m(\C)$ and let $f:X\to Y$ be a map. Endow $X$ and $Y$ with the relative topology induced by the $\Q$-Zariski topology of $\PP^n(\C)$ and $\PP^m(\C)$, respectively. Assume that $f$ is continuous. Let $a\in X$, let $i\in\{0,\dots,n\}$ and let $j\in\{0,\ldots,m\}$ be such that $a\in U^n_i$ and $f(a)\in U^m_j$. Define the map $f_{ij}:\theta^n_i(X\cap U^n_i\cap f^{-1}(Y\cap U^m_j))\to\theta^m_j(Y\cap U^m_j)$ by $f_{ij}(x)=\theta^m_j(f((\theta^n_i)^{-1}(x)))$. We say that $f$ is \emph{$\Q$-regular at $a$} if $f_{ij}$ is $\Q$-regular at $a$ in the sense of Definition \ref{def:Q-reg-function}. We say that $f$ is \emph{$\Q$-regular} if $f$ is $\Q$-regular at each point of $X$. $\sqbullet$
\end{defn}

As above, we denote $\theta:\PP^n(\R)\setminus H\to\R^n$ the affine chart $\theta([x_0,x])=\big(\frac{x_1}{x_0},\ldots,\frac{x_n}{x_0}\big)$, where $H:=\{[x_0,x]\in\PP^n(\R):x_0=0\}$. Similarly, for short, we denote $\theta_\C:\PP^n(\C)\setminus H_\C\to\C^n$ the affine chart $\theta_\C([x_0,x])=\big(\frac{x_1}{x_0},\ldots,\frac{x_n}{x_0}\big)$, where $H_\C:=\{[x_0,x]\in\PP^n(\C):x_0=0\}$. Identify $\R^n$ with $\PP^n(\R)\setminus H$ via $\theta^{-1}$, and $\C^n$ with $\PP^n(\C)\setminus H_\C$ via $\theta_\C^{-1}$. Therefore, we can write $\R^n\subset\PP^n(\R)$ and $\C^n\subset\PP^n(\C)$. Let $j_n:\PP^n(\C)\to\PP^n(\C)$ be the conjugation map $j_n([x_0,\ldots,x_n]):=[\overline{x_0},\ldots,\overline{x_n}]$. Identify $\PP^n(\R)$ with the fixed point set of $j_n$, so we can also write $\PP^n(\R)\subset\PP^n(\C)$ and $\R^n\subset\C^n$.

Let $x=(x',x'')$ be the coordinates of $\C^r\times\C^{n-r}=\C^n$. Consider $x'=(x_1,\ldots,x_r)$, $x''=(x_{r+1},\ldots,x_n)$ and $x$ as column vectors. For every $A\in\mc{M}_{r,n-r}(\Q)$, denote $\widehat{V}_A$ the $(n-r-1)$-dimensional projective subspace of $\PP^n(\C)$ obtained intersecting $H_\C$ with the $\Q$-Zariski closure of $V_A=\{(x',x'')\in\R^n:x'=Ax''\}$ in $\PP^n(\C)$, i.e., $\widehat{V}_A:=\{[0,x',x'']\in\PP^n(\C):x'=Ax''\}$..

Next result collects some basic properties of real and complex projective $\Q$-algebraic sets. 

\begin{prop}\label{Z}
We have:
\begin{itemize}
 \item[$(\mr{i})$] Let $X$ be a $\Q$-algebraic subset of $\PP^n(\R)$ and let $Y$ be the $\Q$-Zariski closure of $X$ in $\PP^n(\C)$. Then $X=Y\cap\PP^n(\R)$.
 \item[$(\mr{ii})$] Let $Z\subset\C^n$ be a $\Q$-algebraic set with $\dim_\C(Z)=d>0$ and let $\overline{Z}$ be the $\Q$-Zariski closure of $Z$ in $\PP^n(\C)$. Then $\overline{Z}\cap H_\C\subset\PP^n(\C)$ is a $\Q$-algebraic set such that $\dim_\C(\overline{Z}\cap H_\C)=d-1$.
 \item[$(\mr{iii})$] Let $A\in\mc{M}_{r,n-r}(\Q)$ and let $\Pi_A:\PP^n(\C)\setminus\widehat{V}_A\to\PP^r(\C)$ be the projection of $\PP^n(\C)$ onto $\PP^r(\C)$ with center $\widehat{V}_A$, i.e., $\Pi_A([x_0,x]):=[x_0,x'-Ax'']$. Then $\Pi_A$ has the following property: if $S$ is a $\Q$-algebraic subset of $\PP^n(\C)$ with $S\cap\widehat{V}_A=\varnothing$, then $\Pi_A(S)$ is a $\Q$-algebraic subset of $\PP^r(\C)$.
  \item[$(\mr{iv})$] Let $W\subset\PP^n(\C)$ be a $\Q$-algebraic set such that $W\subset H_\C$ and $\dim_\C(W)\leq r-1$, and let $\Omega:=\{A\in\mc{M}_{r,n-r}(\Q):W\cap\widehat{V}_A=\varnothing\}$. Then $\Omega$ is a non-empty Zariski open subset of $\mc{M}_{r,n-r}(\Q)$.
  \item[$(\mr{v})$] Let $S$ be a $\Q$-constructible subset of $\PP^n(\C)$ and let $F:S\to\PP^r(\C)$ be a $\Q$-regular map. Then $F(S)$ is a $\Q$-constructible subset of $\PP^r(\C)$.
\end{itemize}
\end{prop}
\begin{proof}
$(\mr{i})$ The inclusion $X\subset Y\cap\PP^n(\R)$ is evident since $X\subset Y$. Let us prove the converse inclusion. By Lemma \ref{lem:Z}, there exists $f\in\Q[\x_0,\x]_\sfh$ such that $\PP\ZZ_\R(f)=X$. Since $X\subset\PP\ZZ_\C(f)$, it follows that $Y\subset\PP\ZZ_\C(f)$ {and therefore} $Y\cap\PP^n(\R)\subset\PP\ZZ_\C(f)\cap\PP^n(\R)=\PP\ZZ_\R(f)=X$. This proves that $X=Y\cap\PP^n(\R)$.

$(\mr{ii})$ By decomposing $Z$ into its $\Q$-irreducible components if necessary, we can assume that $Z\subset\C^n$ is $\Q$-irreducible and therefore so is $\overline{Z}\subset\PP^n(\C)$. Since $\overline{Z}\not\subset H_\C$, we have that $\overline{Z}\cap H_\C$ is a proper $\Q$-Zariski closed subset of $\overline{Z}\subset\PP^n(\C)$. Therefore, using Lemma \ref{dimirred} and Theorem \ref{dimension} via each affine chart $\theta^n_i$ of $\PP^n(\C)$, we obtain $\dim_\C(\overline{Z}\cap H_\C)=\dim_\Q(\overline{Z}\cap H_\C)<\dim_\Q(\overline{Z})=\dim_\C(\overline{Z})=\dim_\C(Z)=d$. Since $\overline{Z}$ is also a usual algebraic subset of $\PP^n(\C)$, we know that $\dim_\C(\overline{Z}\cap H_\C)\geq d+(n-1)-n=d-1$ (see \cite[Prop.3.28\,\&\,Cor.3.30]{mu}), so $\dim_\C(\overline{Z}\cap H_\C)=d-1$.

$(\mr{iii})$ Let $M(x_0,x,y_0,y)$ be the $(r+1)\times2$-matrix whose columns are $(x_0,x'-Ax'')$ and $(y_0,y)=(y_0,y_1,\ldots,y_r)$, let $D_1,\ldots,D_\ell\in\Q[x_0,x,y_0,y]$ be the determinants of all $2\times 2$-submatrices of $M(x_0,x,y_0,y)$, let $E$ be the $\Q$-algebraic subset of $\PP^n(\C)\times\PP^r(\C)$ defined by the polynomial equations $D_1=0,\ldots,D_\ell=0$ and let $\rho:\PP^n(\C)\times\PP^r(\C)\to\PP^r(\C)$ be the projection $(\alpha,\beta)\mapsto\beta$. Pick any $\Q$-algebraic subset $S$ of $\PP^n(\C)$ with $S\cap\widehat{V}_A=\varnothing$. Since $(S\times\PP^r(\C))\cap E\subset\PP^n(\C)\times\PP^r(\C)$ is $\Q$-algebraic and $\rho((S\times\PP^r(\C))\cap E)=\Pi_A(S)$, Theorem \ref{thm:K-elimination}$(\mr{i})$ implies that $\Pi_A(S)\subset\PP^r(\C)$ is also $\Q$-algebraic.

$(\mr{iv})$ This statement coincides with \cite[Lem.6.3.10]{FG} in the case $L|K=\C|\Q$.

$(\mr{v})$ This item is a particular case of \cite[Prop.6.3.12]{FG}. For the sake of completeness, we provide the proof below. Consider again the projection $\rho:\PP^n(\C)\times\PP^r(\C)\to\PP^r(\C)$ defined by $(\alpha,\beta)\mapsto\beta$. Since $S$ is a $\Q$-constructible subset of $\PP^n(\C)$ and $F$ is a $\Q$-regular map, the graph $\Gamma$ of $F$ is a $\Q$-constructible subset of $\PP^n(\C)\times\PP^r(\C)$. By Theorem \ref{thm:K-elimination}$(\mr{ii})$, we have that $F(S)=\rho(\Gamma)\subset\PP^r(\C)$ is also $\Q$-constructible, as required.
\end{proof}

\begin{proof}[Proof of Proposition \ref{213}] 
Let $Z\subset\C^n$ be the $\Q$-Zariski closure of $X$ in $\C^n$ and let $\overline{Z}$ be the $\Q$-Zariski closure of $X$ in $\PP^n(\C)$. Since $\II_\Q(Z)=\II_\Q(X)$, by Definition \ref{def:K-dim} and Theorem \ref{dimension}, we have $\dim_\C(Z)=\dim_\Q(Z)=\dim_\Q(X)=\dim_\R(X)=d$. Observe that $\overline{Z}$ is also the $\Q$-Zariski closure of $Z$ in $\PP^n(\C)$, so $\dim_\C(\overline{Z})=\dim_\C(Z)=d$. Since $X$ is $\Q$-Zariski closed in $\PP^n(\R)$ by hypothesis, Proposition \ref{Z}$(\mr{i})$ implies that $\overline{Z}\cap\PP^n(\R)=X$.

Denote $Z^*$ the $\Q$-algebraic subset $\overline{Z}\cap H_\C$ of $\PP^n(\C)$. By Proposition \ref{Z}$(\mr{ii})$, we know that $\dim_\C(Z^*)=d-1<r-1$. Identify $H_\C$ with $\PP^{n-1}(\C)$ via the map $[0,x]\mapsto[x]$. If $v$ and $w$ are two distinct points of $\PP^{n-1}(\C)$, we denote $\overline{vw}$ the projective line of $\PP^{n-1}(\C)$ connecting them. Define:
\begin{align*}
&S_1:=(Z\times Z)\setminus\{(x,y)\in Z\times Z:x=y\}\subset\C^{2n},\\
&F:S_1\to\PP^{n-1}(\C) \,\text{ by }\, F(x,y):=[x-y],\\
&S^*_1:=(Z^*\times Z^*)\setminus\{(\alpha,\beta)\in Z^*\times Z^*:\alpha=\beta\}\subset(\PP^{n-1}(\C))^2,\\
&S_2\textstyle:=\bigcup_{(\alpha,\beta)\in S^*_1}\overline{\alpha\beta}\subset\PP^{n-1}(\C),\\
&S_3:=F(S_1)\cup S_2,\\
&\text{$\overline{S_3}$ is the $\Q$-Zariski closure of $S_3$ in $\PP^{n-1}(\C)$,}\\
&W:=Z^*\cup\overline{S_3}.
\end{align*}

Since $F$ is a $\Q$-regular map, $S_1\subset\C^{2n}$ is $\Q$-constructible and $\dim_\C(S_1)=2d$, Proposition \ref{Z}$(\mr{v})$ ensures that $F(S_1)\subset\PP^{n-1}(\C)$ is $\Q$-constructible and $\dim_\C(F(S_1))\leq2d$.

We claim that $S_2\subset\PP^{n-1}(\C)$ is $\Q$-constructible and $\dim_\C(S_2)\leq 2d-1$. Let $T:=\{(\alpha,\beta,\gamma)\in(\PP^{n-1}(\C))^3:(\alpha,\beta)\in S^*_1,\gamma\in\overline{\alpha\beta}\}$ and let $\rho:(\PP^{n-1}(\C))^3\to\PP^{n-1}(\C)$ be the projection $(\alpha,\beta,\gamma)\mapsto\gamma$. Write: $\alpha=[a]$, $\beta=[b]$ and $\gamma=[c]$ with $a=(a_1,\ldots,a_n)$, $b=(b_1,\ldots,b_n)$ and $c=(c_1,\ldots,c_n)$. Let $M(a,b,c)$ be the $n\times3$-matrix whose columns are $a$, $b$ and $c$, let $D_1,\ldots,D_\ell\in\Q[a,b,c]$ be the determinants of all the $3\times 3$-submatrices of $M(a,b,c)$ and let $E$ be the $\Q$-algebraic subset of $(\PP^{n-1}(\C))^3$ defined by the equations $D_1=0,\ldots,D_\ell=0$. Since $T=(S^*_1\times\PP^{n-1}(\C))\cap E$, it follows that $T\subset(\PP^{n-1}(\C))^3$ is $\Q$-constructible. Moreover, we have $\dim_\C(T)\leq2\dim_\C(Z^*)+1=2(d-1)+1=2d-1$. Since $S_2=\rho(T)$, Theorem \ref{thm:K-elimination}$(\mr{ii})$ also ensures that $S_2\subset\PP^{n-1}(\C)$ is $\Q$-constructible. Moreover, $\dim_\C(S_2)\leq 2d-1$ as claimed.

We have just proved that $F(S_1)$ and $S_2$ are $\Q$-constructible subsets of $\PP^{n-1}(\C)$, $\dim_\C(F(S_1))\leq2d$ and $\dim_\C(S_2)\leq2d-1$. Thus, $S_3=F(S_1)\cup S_2$ is a $\Q$-constructible subset of~$\PP^{n-1}(\C)$ and $\dim_\C(S_3)\leq2d$. It follows that $\overline{S_3}\subset\PP^{n-1}(\C)$ is a $\Q$-algebraic set with $\dim_\C(\overline{S_3})\leq2d\leq r-1$. Since $\dim_\C(Z^*)<r-1$, we deduce that $W=Z^*\cup\overline{S_3}\subset\PP^{n-1}(\C)$ is also a $\Q$-algebraic set with $\dim_\C(W)\leq r-1$.

By Proposition \ref{Z}$(\mr{iv})$, the set $\Omega:=\{A\in\mc{M}_{r,n-r}(\Q):W\cap\widehat{V}_A=\varnothing\}$ is a non-empty Zariski open subset of $\mc{M}_{r,n-r}(\Q)$. Pick any $A\in\Omega$. Let $\Pi_A:\PP^n(\C)\setminus\widehat{V}_A\to\PP^r(\C)$ be the projection of $\PP^n(\C)$ onto $\PP^r(\C)$ with center $\widehat{V}_A$, i.e., $\Pi_A([x_0,x]):=[x_0,x'-Ax'']$, let $Y:=\Pi_A(\overline{Z})$ and let $\Pi'_A:\overline{Z}\to Y$ be the restriction of $\Pi_A$ from $\overline{Z}$ to $Y$. Observe that the restriction of $\Pi_A$ from $\R^n$ onto $\R^r$ coincides with $\pi_A$. Since $Z^*\cup S_3\subset W$ and $W\cap\widehat{V}_A=\varnothing$, we have both $Z^*\cap\widehat{V}_A=\varnothing$ and $S_3\cap\widehat{V}_A=\varnothing$. The equality $Z^*\cap\widehat{V}_A=\varnothing$ and Proposition \ref{Z}$(\mr{iii})$ ensure that $Y\subset\PP^r(\C)$ is $\Q$-algebraic. The equality $S_3\cap\widehat{V}_A=\varnothing$ implies that $\Pi'_A$ is injective. Since $\overline{Z}\cap\PP^n(\R)=X\subset\R^n$, we deduce:
$$
\pi_A(X)=\Pi_A(\overline{Z}\cap\PP^n(\R))\subset\Pi_A(\overline{Z})\cap\Pi_A(\PP^n(\R)\setminus\widehat{V}_A)=Y\cap\PP^r(\R).
$$ 

We claim that $\pi_A(X)=Y\cap\PP^r(\R)$. Let $q\in Y\cap\PP^r(\R)$. Since $Y=\Pi_A(\overline{Z})$, there exists $p\in\overline{Z}$ such that $\Pi_A(p)=q$. We have to prove that $p\in\PP^n(\R)$. Suppose this is not {the case}, i.e., $j_n(p)\neq p$. Since $\overline{Z}\subset\PP^n(\C)$ is $\Q$-algebraic, $\overline{Z}$ is the set of solutions in $\PP^n(\C)$ of some homogeneous polynomial equations with coefficients in $\Q$ and thus in $\R$. It follows that $j_n(p)$ is a point of $\overline{Z}$ {as $p$ is so}. Observe that $j_r(q)=q$ and $\Pi_A(j_n(p))=j_r(\Pi_A(p))$ so
$$
\Pi'_A(j_n(p))=\Pi_A(j_n(p))=j_r(\Pi_A(p))=j_r(q)=q=\Pi_A(p)=\Pi'_A(p),
$$
which is a contradiction since $\Pi'_A$ is injective. This proves that $\pi_A(X)=Y\cap\PP^r(\R)$ and {therefore} $\pi_A(X)\subset\PP^r(\R)$ is $\Q$-algebraic as $Y\subset\PP^r(\C)$ is so.

Since $\pi_A(X)\subset\PP^r(\R)$ is $\Q$-algebraic and $\pi_A(X)$ is contained in $\R^r$, we have that $\pi_A(X)$ is a projectively $\Q$-closed $\Q$-algebraic subset of $\R^r$, as required.
\end{proof}

\begin{remark}\label{rem111}
If we replace $r=2d+1$ with a larger natural number $r$, the statements of \cite[Cor.6.3.5]{FG}, Proposition \ref{213} and therefore Theorem \ref{thm:generic-projection} remain valid since their proofs clearly are not effected by this change. $\sqbullet$ 
\end{remark}


\subsection{$\Q$-algebraic models of some special manifolds}\label{subsec:Q-embeddings} Let $n,k\in\N$ be such that $n>0$ and $k\leq n$. Identify the set of real $n\times n$-matrices with~$\R^{n^2}$ and the Grassmannian $\G_{n,k}$ of $k$-dimensional real vector subspaces of $\R^n$ with the set
\begin{equation}\label{grassmannian}
\G_{n,k}:=\big\{X\in\R^{n^2}:X^T=X,X^2=X,\mr{tr}(X)=k\big\},
\end{equation}
where $X^T$ is the transpose of the matrix $X$ and $\mr{tr}(X)$ is the trace of $X$. It is well-known that $\G_{n,k}\subset\R^{n^2}$ is a nonsingular algebraic set of dimension $k(n-k)$, see \cite[Prop.3.4.3\;\&\;Thm.3.4.4]{BCR}. Evidently, $\G_{n,k}\subset\R^{n^2}$ is also $\Q$-algebraic {since} the polynomial equations $X^T-X=0$, $X^2-X=0$ and $\mr{tr}(X)-k=0$ in $\R^{n^2}$ have rational coefficients.

Consider each element $y\in\R^n$ as a column vector. Denote $\E_{n,k}$ the (total space of the) universal vector bundle over $\G_{n,k}$, i.e., $\E_{n,k}:=\{(X,y)\in\G_{n,k}\times\R^n:Xy=y\}$. Observe that $\E_{n,k}$ is a nonsingular $\Q$-algebraic subset of $\R^{n^2+n}$ of dimension $k(n-k+1)$.

\begin{lem} \label{lem:Q-grassmannians}
We have:
\begin{itemize}
 \item[$(\mr{i})$] $\G_{n,k}\subset\R^{n^2}$ is a projectively $\Q$-closed $\Q$-nonsingular $\Q$-algebraic set.
 \item[$(\mr{ii})$] $\E_{n,k}\subset\R^{n^2+n}$ is a $\Q$-nonsingular $\Q$-algebraic set.
 \end{itemize}
\end{lem}
\begin{proof}
$(\mr{i})$ If $X\in\G_{n,k}$, then $\mr{tr}(XX^T)=\mr{tr}(X^2)=\mr{tr}(X)=k$ so $\G_{n,k}$ is contained in the Euclidean sphere $\sph:=\{X\in\R^{n^2}:\mr{tr}(XX^T)-k=0\}$ of $\R^{n^2}$ centered at the origin of radius $\sqrt{k}$. Since $\sph$ is projectively $\Q$-closed by Lemma \ref{lem:overt}, $\G_{n,k}$ is projectively $\Q$-closed by Lemma~\ref{lem:projective}$(\mr{ii})$. 

It remains to show that $\G_{n,k}$ is $\Q$-nonsingular. Since $\G_{n,k}$ is nonsingular and of dimension $k(n-k)$, by Proposition \ref{Y}, it suffices to prove that $\dim(T^\Q_A(\G_{n,k}))=k(n-k)$ for all $A\in\G_{n,k}$. Define the polynomial map $\phi:\R^{n^2}\to\R^{n^2}\times\R^{n^2}=\R^{2n^2}$ by $\phi(X):=(X^T-X,X^2-X)$ and denote $J_\phi(A)$ the Jacobian matrix of $\phi$ at $A\in\G_{n,k}$. To complete the proof, it is enough to see that $\mr{rk}(J_\phi(A))\geq n^2-k(n-k)$. Indeed, if that were the case, then $\dim(T^\Q_A(\G_{n,k}))\leq k(n-k)$ and {therefore} $\dim(T^\Q_A(\G_{n,k}))=k(n-k)$ {since} $\dim(T^\Q_A(\G_{n,k}))\geq\dim(T_A(\G_{n,k}))\geq k(n-k)$.

First, we prove that $\mr{rk}(J_\phi(D))\geq n^2-k(n-k)$, where $D\in\G_{n,k}$ is the diagonal matrix in $\R^{n^2}$ having $1$ in the first $k$ diagonal positions and $0$ otherwise. Then we will show that the same inequality holds for every $A\in\G_{n,k}$. Denote $x=(x_{ij})_{i,j\in\{1,\ldots,n\}}$ the coordinates of $\R^{n^2}$ and $\x=(\x_{ij})_{i,j\in\{1,\ldots,n\}}$ the corresponding indeterminates. For each $i,j\in\{1,\ldots,n\}$, define the polynomials $g_{ij},f_{ij}\in\II_\Q(\G_{n,k})$ by $f_{ij}(\x):=\x_{ij}-\x_{ji}$ and $g_{ij}(\x):=\textstyle(\sum_{\ell=1}^n\x_{i\ell}\x_{\ell j})-\x_{ij}$. It follows that $\phi(x)=((f_{ij}(x))_{i,j},(g_{ij}(x))_{i,j})$ for every $x\in\R^{n^2}$. Define
\begin{align*}
S_1&:=\{(i,j)\in\{1,\ldots,n\}^2\,|\,i<j\},\\
S_2&:=\{(i,j)\in\{1,\ldots,n\}^2\,|\,i\leq j\leq k\},\\
S_3&:=\{(i,j)\in\{1,\ldots,n\}^2\,|\,k<i\leq j\}.
\end{align*}
Observe that the sum of the cardinalities of $S_1$, $S_2$ and $S_3$ is equal to $\frac{(n-1)n}{2}+\frac{k(k+1)}{2}+\frac{(n-k)(n-k+1)}{2}=n^2-k(n-k)$. By a direct computation, we see that
$$
\begin{array}{ll}
\nabla f_{ij}(D)=E_{ij}-E_{ji} & \text{ if $(i,j)\in S_1$,}\\
\nabla g_{ij}(D)=E_{ij} & \text{ if $(i,j)\in S_2$,}\\
\nabla g_{ij}(D)=-E_{ij} & \text{ if $(i,j)\in S_3$,}
\end{array}
$$
where $E_{ij}$ is the matrix in $\R^{n^2}$ whose $(i,j)$-coefficient is equal to $1$ and $0$ otherwise. It follows that the gradients listed above are linearly independent, so $\mr{rk}(J_\phi(D))\geq n^2-k(n-k)$, as claimed. Let $A\in\G_{n,k}$ and let $G\in O(n)$ be such that $D=G^TAG$. Define the linear automorphism $\psi:\R^{n^2}\to\R^{n^2}$ by $\psi(X):=G^TXG$. Since $\psi(A)=D$ and $(\psi\times\psi)\circ\phi=\phi\circ\psi$, we have that $J_{\psi\times\psi}(\phi(A))J_\phi(A)=J_\phi(D)J_\psi(A)$. Since both matrices $J_{\psi\times\psi}(\phi(A))$ and $J_\psi(A)$ are invertible, we have $\mr{rk}(J_\phi(A))=\mr{rk}(J_\phi(D))\geq n^2-k(n-k)$.

$(\mr{ii})$ We keep the notations of the previous part of the proof. Let $\Phi:\R^{n^2+n}\to\R^{2n^2+n}$ be the polynomial map $\Phi(X,y):=(X^T-X,X^2-X,Xy-y)$. As above, it suffices to show that $\mr{rk}(J_\Phi(A,b))\geq n^2+n-k(n-k+1)=n^2-k(n-k)+n-k$ for all $(A,b)\in\E_{n,k}$.

Define the polynomial $h_\ell\in\II_\Q(\E_{n,k})$ by $h_\ell(X,y):=(\sum_{j=1}^n\x_{\ell j}\y_j)-\y_\ell$ for every $\ell\in\{1,\dots,n\}$, so $\Phi(X,y)=(\phi(X),h_1(X,y),\ldots,h_n(X,y))$. Let $v=(v_1,\ldots,v_n)\in\R^n$ be such that $(D,v)\in\E_{n,k}$. By a simple computation, we see that $\nabla h_\ell(D,v)=\big(\sum_{j=1}^n v_j E_{\ell j},-e_\ell\big)\in\R^{n^2}\times\R^n$ for every $\ell\in\{k+1,\dots,n\}$, where $\{e_1,\ldots,e_n\}$ is the canonical vector basis of $\R^n$. Since $\mr{rk}(J_\phi(D))\geq n^2-k(n-k)$, we deduce that $\mr{rk}(J_\Phi(D,v))\geq n^2-k(n-k)+n-k$.

Pick any $(A,w)\in\E_{n,k}$. Let $G\in O(n)$ be such that $D=G^TAG$ and set $z:=G^Tw$. Observe that $Dz=G^TAGG^Tw=G^TAw=G^Tw=z$, i.e., $(D,z)\in\E_{n,k}$. Define the linear automorphism $\psi:\R^{n^2}\to\R^{n^2}$ as above $\psi(X):=G^TXG$, and the linear automorphism $\tau:\R^n\rightarrow\R^n$ by $\tau(y):=G^Ty$. Since $(\psi\times\tau)(A,w)=(D,z)$ and $(\psi\times\psi\times\tau)\circ\Phi=\Phi\circ(\psi\times\tau)$, we have that $J_{\psi\times\psi\times\tau}(\Phi(A,w))J_\Phi(A,w)=J_\Phi(D,z)J_{\psi\times\tau}(A,w)$. Bearing in mind that both matrices $J_{\psi\times\psi\times\tau}(\Phi(A,w))$ and $J_{\psi\times\tau}(A,w)$ are invertible, it follows that $\mr{rk}(J_\Phi(A,w))=\mr{rk}(J_\Phi(D,z))\geq n^2-k(n-k)+n-k$, as required.
\end{proof}

\begin{lem}\label{lem:gauss}
Let $X\subset\R^n$ be a $\Q$-nonsingular $\Q$-algebraic set of dimension $d$. We have:
\begin{itemize}
 \item[$(\mr{i})$] Let $\beta:X\rightarrow\G_{n,n-d}$ be the normal bundle map of $X$ in $\R^n$, i.e., $\beta(x)\in\R^{n^2}$ is the matrix associated to the orthogonal projection of $\R^n$ onto $T_x(X)^\perp$ with respect to the canonical basis of $\R^n$ for all $x\in X$. Then $\beta\in\reg^\Q(X,\G_{n,n-d})$.
 \item[$(\mr{ii})$] Let $N$ be the total space of the normal bundle of $X$ in $\R^n$, i.e., $N:=\{(x,y)\in X\times\R^n:y\in T_x(X)^\perp\}$. Then $N\subset\R^{2n}$ is a $\Q$-nonsingular $\Q$-algebraic set of dimension $n$.
\end{itemize}
\end{lem}
\begin{proof}
If $d=n$, the result is evident {as} $X=\R^n$. Suppose that $d<n$.

$(\mr{i})$ Let $a\in X$, let $f_1,\ldots,f_{n-d}\in\II_\Q(X)$ be such that $\nabla f_1(a),\ldots,\nabla f_{n-d}(a)$ are linearly independent and let $A:\R^n\to\R^{n^2}$ be the polynomial map $A(x):=\big(\frac{\partial f_i}{\partial\x_j}(x)\big)_{i=1,\ldots,n-d,\,j=1,\ldots,n}$.

First, we prove that there exists a $\Q$-Zariski open neighborhood $U$ of $a$ in $\R^n$ such that $X\cap U=\ZZ_\R(f_1,\ldots,f_{n-d})\cap U$ and $\mr{rk}(A(x))=n-d$ for all $x\in X\cap U$. By Lemma \ref{K-difference-preparation}, there exist a unique $\Q$-irreducible component $X'$ of $X$ such that $a\in X'$. Denote {by} $X''$ the union of all $\Q$-irreducible components of $X$ different from $X'$. Define $W:=\ZZ_\R(f_1,\ldots,f_{n-d})\subset\R^n$. By Propositions \ref{jacobian-criterion} and \ref{prop211}$(\mr{i})$, we deduce that there exist a unique $\Q$-irreducible component $W'$ of $W$ such that $a\in W'$ and $\dim_\R(W')=d=\dim_\R(X')$. Denote {by} $W''$ the union of all $\Q$-irreducible components of $W$ different from $W'$. By Lemma \ref{X} and \cite[Prop.3.3.11]{BCR}, there exists a semialgebraic open neighborhood $V$ of $a$ in $\R^n$ such that $X'\cap V=\ZZ_\R(f_1,\ldots,f_{n-d})\cap V=W'\cap V$ and $X'\cap V=W'\cap V$ is a Nash submanifold of $\R^n$ of dimension $d$. Observe that $(X'\cap W')\cap V=X'\cap V=W'\cap V$ and $X'\cap W'\subset X$. It follows that $\dim_\R(X'\cap W')=d$ {since} $d=\dim_\R((X'\cap W')\cap V)\leq\dim_\R(X'\cap W')\leq\dim_\R(X)=d$. By Theorem \ref{dimension}, we have $\dim_\Q(X')=\dim_\Q(W')=\dim_\Q(X'\cap W')=d$, so $X'=X'\cap W'=W'$ by Lemma \ref{dimirred}. The set $U:=\{x\in\R^n:\mr{rk}(A(x))=n-d\}\setminus(X''\cup W'')$ has the required properties.

By Definition \ref{def:Q-reg-function}, it suffices to show that the restriction $\beta|_{X\cap U}$ of $\beta$ to $X\cap U$ is $\Q$-regular. Let $x\in X\cap U$. Since $\mr{rk}(A(x))=n-d$, we have that the matrix $A(x)A(x)^T\in\R^{(n-d)^2}$ is invertible, so we can define the matrix $\eta(x)\in\G_{n,n-d}$ by $\eta(x):=A(x)^T(A(x)A(x)^T)^{-1}A(x)$. Since $\eta(x)v=0$ for all $v\in T_x(X)$, we deduce that $\eta(x)=\beta(x)$. Observe that the map $U\to\R^{n^2}$, $x\mapsto\eta(x)$ is $\Q$-regular since the matrix $(A(x)A(x)^T)^{-1}$ is computed in terms of the determinant and the minors of order $n-1$ of $A(x)A(x)^T$, which are polynomials with rational coefficients with respect to the entries of $A(x)A(x)^T$, thus the map $\beta|_{X\cap U}$ is also $\Q$-regular.

$(\mr{ii})$ Let $\widehat{\beta}:X\times\R^n\to\G_{n,n-d}\times\R^n$ be the map $\widehat{\beta}(x,y):=(\beta(x),y)$, which is $\Q$-regular by~$(\mr{i})$. Since $\widehat{\beta}$ is transverse to $\E_{n,n-d}$ in $\G_{n,n-d}\times\R^n$ and $\widehat{\beta}^{-1}(\E_{n,n-d})=N$, this item follows immediately from Proposition \ref{prop:Q-transverse}. 
\end{proof}

\begin{remark}\label{rem:gauss}
The proof of item $(\mr{i})$ of the previous corollary can be easily adapted to show that, given any $\Q$-algebraic set $X\subset\R^n$ of dimension $d$, the normal bundle map $\beta:\Reg^\Q(X)\to\G_{n,n-d}$ of $\Reg^\Q(X)$ in $\R^n$ is $\Q$-regular. $\sqbullet$
\end{remark}

Let $(x_0,x)=(x_0,x_1,\ldots,x_n)\in\R^{n+1}\setminus\{0\}$, let $|(x_0,x)|_{n+1}$ be the Euclidean norm of $(x_0,x)\in\R^{n+1}$ and let $[x_0,x]=[x_0,\ldots,x_n]$ be the point in $\PP^n(\R)$ corresponding to $(x_0,x)$. We consider $(x_0,x)\in\R^{n+1}$ as a column vector. Observe that $(x_0,x)(x_0,x)^T|(x_0,x)|_{n+1}^{-2}\in\G_{n+1,1}$ is the matrix associated to the orthogonal projection of $\R^{n+1}$ onto the vector line generated by $(x_0,x)$ with respect to the canonical vector basis of $\R^{n+1}$. Denote $\eta_n:\PP^n(\R)\to\G_{n+1,1}$ the biregular isomorphism defined by
\begin{equation}\label{eq:mu}
\textstyle
\eta_n([x_0,x]):=(x_0,x)(x_0,x)^T|(x_0,x)|_{n+1}^{-2}=\big(x_ix_j(\sum_{i=0}^nx_i^2)^{-1}\big)_{i,j\in\{0,\ldots,n\}}.
\end{equation}

\emph{Fix $m\in\N^*$ with $n\leq m$.} Let $P\in\R[\x_0,\x,\y_0,\y]:=\R[\x_0,\ldots,\x_n,\y_0,\ldots,\y_m]$, let $H_{n,m}$ be the nonsingular algebraic hypersurface of $\PP^n(\R)\times\PP^m(\R)$ and let $\HH_{n,m}$ be the nonsingular algebraic hypersurface of $\G_{n+1,1}\times\G_{m+1,1}$ defined by
\begin{align}
&\textstyle P(\x_0,\x,\y_0,\y):=\sum_{i=0}^n\x_i\y_i,\label{eq:milnor1}\\
&H_{n,m}:=\{([x_0,x],[y_0,y])\in\PP^n(\R) \times\PP^m(\R):P(x_0,x,y_0,y)=0\},\label{eq:milnor2}\\
&\HH_{n,m}:=(\eta_n\times\eta_m)(H_{n,m}).\label{eq:milnor3}
\end{align}
Evidently, $H_{n,m}\subset\PP^n(\R)\times\PP^m(\R)$ is $\Q$-algebraic {since} $P$ has rational coefficients.

Recall that $\theta^n_k:U^n_k\to\R^n$ is the affine chart $[x]\mapsto\big(\frac{x_0}{x_k},\ldots,\frac{x_{k-1}}{x_k},\frac{x_{k+1}}{x_k},\ldots,\frac{x_n}{x_k}\big)$, where $k\in\{0,\ldots,n\}$ and $U^n_k=\{[(x_0,x)]\in\PP^n(\R):x_k\neq0\}$, and $(\theta^n_k)^{-1}(x)=[x_1,\ldots,x_k,1,x_{k+1},\ldots,x_n]$.

\begin{lem} \label{lem:Q-Hnm}
$\HH_{n,m}\subset\R^{(n+1)^2+(m+1)^2}$ is a projectively $\Q$-closed $\Q$-nonsingular $\Q$-algebraic set. 
\end{lem}
\begin{proof}
Let $X=(x_{ij})_{i,j}$ and $Y=(y_{ij})_{i,j}$ be the coordinates of $\R^{(n+1)^2}$ and $\R^{(m+1)^2}$, respectively. Here $(i,j)$ varies in $\{0,\ldots,n\}\times\{0,\ldots,m\}$. Set $N:=(n+1)^2+(m+1)^2$. For every $k\in\{0,\ldots,n\}$, define $\Omega^n_k:=\{X\in\G_{n+1,1}:x_{kk}\neq0\}$, $\xi^n_k:\Omega^n_k\to U^n_k$ by $\xi^n_k(X):=[x_{0k},\ldots,x_{nk}]$ and $\eta^n_k:\Omega^n_k\to\R^n$ by $\eta^n_k:=\theta^n_k\circ\xi^n_k$. It is easy to verify that $\eta_n(U^n_k)=\Omega^n_k$, $\eta_n^{-1}(X)=\xi^n_k(X)$ for all $X\in\Omega^n_k$, and $\eta^n_k$ is a $\Q$-biregular isomorphism. In particular, the products $\Omega^n_k\times\Omega^m_h$ form a $\Q$-Zariski open cover of $\G_{n+1,1}\times\G_{m+1,1}$. For every $(k,h)\in\{0,\ldots,n\}\times\{0,\ldots,m\}$, let $f_{kh}\in\Q[\x,\y]:=\Q[\x_1,\ldots,\x_n,\y_1,\ldots,\y_m]$ be the polynomial defined by
$$
f_{kh}(\x,\y):=P(\x_1,\ldots,\x_k,1,\x_{k+1},\ldots,\x_n,\y_1,\ldots,\y_h,1,\y_{h+1},\ldots,\y_m)
$$
and let $V_{kh}\subset\R^{n+m}$ be the $\Q$-algebraic set $V_{kh}:=\ZZ_\R(f_{kh})$. Since $\eta^n_k\times\eta^m_h:\Omega^n_k\times\Omega^m_h\to\R^n\times\R^m=\R^{n+m}$ is $\Q$-regular and $(\eta^n_k\times\eta^m_h)^{-1}(V_{kh})=\HH_{n,m}\cap(\Omega^n_k\times\Omega^m_h)$, we deduce that $\HH_{n,m}\cap(\Omega^n_k\times\Omega^m_h)$ is $\Q$-Zariski closed in $\Omega^n_k\times\Omega^m_h$. Since this holds for every pair $(k,h)$, the whole set $\HH_{n,m}$ is $\Q$-Zariski closed in $\G_{n+1,1}\times\G_{m+1,1}$. By Lemma \ref{lem:projective}$(\mr{iv})$ and Lemma \ref{lem:Q-grassmannians}$(\mr{i})$, $\G_{n+1,1}\times\G_{m+1,1}\subset\R^N$ is projectively $\Q$-closed and $\Q$-nonsingular. Since $\HH_{n,k}\subset\G_{n+1,1}\times\G_{m+1,1}$, the set $\HH_{n,k}\subset\R^N$ is projectively $\Q$-closed by  Lemma \ref{lem:projective}$(\mr{ii})$.

It remains to prove that $\HH_{n,m}\subset\R^N$ is $\Q$-nonsingular. Let $a:=(X,Y)\in\HH_{n,m}$, let $(k,h)\in\{0,\ldots,n\}\times\{0,\ldots,m\}$ be such that $a\in\HH_{n,m}\cap(\Omega^n_k\times\Omega^m_h)$, and let $b:=(\eta^n_k(X),\eta^m_h(Y))\in V_{kh}$. By an easy computation, we see that $\nabla f_{kh}$ never vanishes on $V_{kh}$, so $\nabla f_{kh}(b)\neq0$. Since $\eta^n_k\times\eta^m_h$ is a $\Q$-biregular isomorphism and $(\eta^n_k\times\eta^m_h)^{-1}(V_{kh})=\HH_{n,m}\cap(\Omega^n_k\times\Omega^m_h)$, it follows that $T^\Q_a(\HH_{n,m})=d_a(\eta^n_k\times\eta^m_h)^{-1}(\mr{Span}(\nabla f_{kh}(b))^\perp)$, so $\dim(T^\Q_a(\HH_{n,m}))= n+m-1=\dim(T_a(\HH_{n,m}))$. By Proposition \ref{Y}, we deduce that $\HH_{n,m}\subset\R^N$ is $\Q$-nonsingular.
\end{proof}


\subsection{Projectively $\Q$-algebraic unoriented bordism and  homology}\label{subsec:Q-bordism-homology}
Let $d\in\N$ and let $Z$ be either the empty set or a $d$-dimensional compact smooth manifold. A $d$-dimensional compact smooth manifold $M$ is (unoriented) cobordant to $Z$ if there exists a $(d+1)$-dimensional compact manifold with boundary $T$ such that $\partial T$ is smoothly diffeomorphic to the disjoint union $M\sqcup Z$. Being cobordant defines an equivalence relation on $d$-dimensional compact smooth mani\-folds. The quotient set $\mk{N}_d$ of the resulting cobordism classes is an abelian group with the addition induced by the disjoint union, called $d$-dimensional unoriented cobordism group. The neutral element of $\mk{N}_d$ is the cobordism class of a $d$-dimensional compact smooth manifold cobordant to the empty set.  

Let $\mk{N}_*=\bigoplus_{d\in\N}\mk{N}_d$ be the unoriented cobordism group. The cartesian product induces a multiplication on $\mk{N}_*$ that makes it a graded algebra over $\Z/2$. In \cite{Th}, Thom proved that $\mk{N}_*$ is a polynomial algebra over $\Z/2$ generated by a $d$-dimensional compact smooth manifold $M_d$ for every $d\in\N$ not of the form $2^i-1$, where $M_d$ can be chosen to be $\PP^d(\R)$ if $d$ is even. In \cite{milnor}, Milnor proved that, if $d$ is odd, we can also assume that $M_d$ has the form $H_{n,m}$, see \eqref{eq:milnor2}. Recall that $\PP^d(\R)$ and $H_{n,m}$ are biregularly isomorphic to $\G_{d+1,1}$ and $\HH_{n,m}$ via $\eta_d$ and $\eta_n\times\eta_m$ respectively, see \eqref{eq:mu} and \eqref{eq:milnor3}.

\begin{lem}\label{lem:Q-cobordism}
For every $d\in\N$, every $d$-dimensional compact smooth manifold $M$ is cobordant to either the empty set or a $d$-dimensional projectively $\Q$-closed $\Q$-nonsingular $\Q$-algebraic subset $Z$ of some $\R^{N_d}$ of the following form: $Z$ is a finite disjoint union of subsets $Y+v$ of $\R^{N_d}$, where $v$ belongs to $\Q^{N_d}$ and
\begin{equation}\label{eq:Q-milnor}
Y=\G_{n_1+1,1}\times\ldots\times\G_{n_\alpha+1,1}\times\HH_{p_1,q_1}\times\ldots\times\HH_{p_\beta,q_\beta}
\end{equation}
for some $\alpha,\beta\in\N$ and $n_1,\ldots,n_\alpha,p_1,q_1,\ldots,p_\beta,q_\beta\in\N^*$, where $Y=\HH_{p_1,q_1}\times\cdots\times\HH_{p_\beta,q_\beta}$ if $\alpha=0$ and $\beta>0$, $Y=\G_{n_1+1,1}\times\ldots\times\G_{n_\alpha+1,1}$ if $\alpha>0$ and $\beta=0$, and $Y$ is the origin of $\R^{N_d}$ if $\alpha=\beta=0$.
\end{lem}
\begin{proof}
Let $\{Y_\ell\}_\ell$ be the finite family of all $d$-dimensional real algebraic sets $Y_\ell$ of the form~\eqref{eq:Q-milnor}, i.e.,  $Y_\ell=\G_{n_{\ell,1}+1,1}\times\cdots\times\G_{n_{\ell,\alpha_\ell}+1,1}\times\HH_{p_{\ell,1},q_{\ell,1}}\times\cdots\times \HH_{p_{\ell,\beta_\ell},q_{\ell,\beta_\ell}}$ for some $\alpha_\ell,\beta_\ell\in\N$ and $n_{\ell,1},\ldots,n_{\ell,\alpha_\ell},p_{\ell,1},q_{\ell,1},\ldots,p_{\ell,\beta_\ell},q_{\ell,\beta_\ell}\in\N^*$, where $\dim(Y_\ell)=d$. Let $N_d$ be a sufficiently large integer such that every $Y_\ell$ of such a family is contained in $\R^{N_d}$. By Lemmas \ref{lem:projective}$(\mr{iv})$, \ref{lem:Q-grassmannians}$(\mr{i})$ and \ref{lem:Q-Hnm}, we know that every $Y_\ell\subset\R^{N_d}$ is a projectively $\Q$-closed $\Q$-nonsingular $\Q$-algebraic sets. For every $\ell$, choose a vector $v_\ell\in\Q^{N_d}$ such that the sets $\{Y_\ell+v_\ell\}_\ell$ are pairwise disjoint. By Lemma \ref{lem:projective}$(\mr{i})(\mr{iii})$, every finite union of sets $Y_\ell+v_\ell\subset\R^{N_d}$ is a projectively $\Q$-closed $\Q$-nonsingular $\Q$-algebraic set. The above-mentioned results of Thom and Milnor \cite{Th,milnor} prove that the finite unions of the sets $\{Y_\ell+v_\ell\}_\ell$ form a system of generators of $\mk{N}_d$. 
\end{proof}

\begin{remark}\label{rmk:N_d}
In the statement of Lemma \ref{lem:Q-cobordism}, the natural number $N_d$ can be chosen to be $(2d+1)^2$. If $d=0$, this is evident. Suppose $d>0$. Let us keep the notations used in the proof of that lemma. Observe that $\dim(Y_\ell+v_\ell)=\dim(Y_\ell)=\sum_{i=1}^{\alpha_\ell}n_{\ell,i}+\sum_{j=1}^{\beta_\ell}(p_{\ell,j}+q_{\ell,j}-1)=d$ and $Y_\ell+v_\ell$ is contained in $\R^{N_{d,\ell}}$, where
$$
\textstyle
N_{d,\ell}:=\sum_{i=1}^{\alpha_\ell}(n_{\ell,i}+1)^2+\sum_{j=1}^{\beta_\ell}((p_{\ell,j}+1)^2+(q_{\ell,j}+1)^2).
$$
It is easy to prove by induction on $k\in\N^*$ that the following inequality holds:
\begin{equation*}\label{eq:estimate}
\textstyle
\text{$\sum_{s=1}^k(c_s+1)^2\leq\big(1+\sum_{s=1}^kc_s\big)^2\;$ for all $k\in\N^*$ and $(c_1,\ldots,c_k)\in(\N^*)^k$.}
\end{equation*}
Thus, we deduce
$$
\textstyle
N_{d,\ell}\leq\big(1+\sum_{i=1}^{\alpha_\ell}n_{\ell,i}+\sum_{j=1}^{\beta_\ell}(p_{\ell,j}+q_{\ell,j})\big)^2=(1+d+\beta_\ell)^2\leq(2d+1)^2.
$$
If we set $N_d:=(2d+1)^2$, then $N_{d,\ell}\leq N_d$ for all $\ell$, as required. $\sqbullet$
\end{remark}

In what follows, $W$ denotes a set, which is tacitly assumed to be non-empty.

Let $W$ be a finite CW complex, let $\mk{N}_*(W)=\bigoplus_{d\in\N}\mk{N}_d(W)$ be the unoriented bordism group of~$W$ and let $H_*(W,\Z/2)=\bigoplus_{d\in\N}H_d(W,\Z/2)$ be the singular homology group of $W$ over $\Z/2$. For every $d\in\N$, we denote {by} $\mu_d:\mk{N}_d(W)\to H_d(W,\Z/2)$ the natural homomorphism defined as follows: if $f:M\to W$ is a continuous map from a $d$-dimensional compact smooth manifold $M$ to $W$, if $[f:M\to W]$ is the corresponding element in $\mk{N}_d(W)$ and if $[M]\in H_d(M,\Z/2)$ is the fundamental class of $M$, then $\mu_d([f:M\to W]):=f_*([M])$. In \cite{Th}, Thom proved that every homomorphism $\mu_d$ is surjective. In \cite{CF},  Conner and Floyd proved the following {result}:

\begin{thm}[{\cite[Thm.(17.1), p.46]{CF}}]\label{thm:CF}
Let $W$ be a finite CW~complex. Choose a basis $\{\alpha_{d,i}\}_i$ of each $H_d(W,\Z/2)\neq\{0\}$ and, for every $\alpha_{d,i}$, a continuous map $f_{d,i}:M_{d,i}\to W$ from a $d$-dimen\-sional compact smooth manifold $M_{d,i}$ to $W$ such that $\mu_d([f_{d,i}:M_{d,i}\to W])=\alpha_{d,i}$ (which exists by the surjectivity of $\mu_d$). Define the $\mk{N}_*$-homomorphism $h:H_*(W,\Z/2)\otimes\mk{N}_*\to\mk{N}_*(W)$ by $h(\alpha_{d,i}\otimes 1):=[f_{d,i}:M_{d,i}\to W]$ for all $\alpha_{d,i}$. Then $h$ is a $\mk{N}_*$-isomorphism.
\end{thm}

Let us introduce the concepts of projectively $\Q$-algebraic homology and unoriented bordism classes. Recall that every subset of every $\R^k$ is assumed to be endowed with the Euclidean topology. 

\begin{defn}\label{def:Q-homology}
Let $W$ be a subset of $\R^k$, let $d\in\N$ and let $\alpha\in H_d(W,\Z/2)$. We say that $\alpha\in H_d(W,\Z/2)$ is \emph{projectively $\Q$-algebraic} if there exist a $d$-dimensional projectively $\Q$-closed $\Q$-nonsingular $\Q$-algebraic set $P\subset\R^m$ and a $\Q$-regular map $g:P\rightarrow W$ such that $g_*([P])=\alpha$, i.e., $\mu_d([g:P\rightarrow W])=\alpha$. If for all $d\in\N$ every element of $H_d(W,\Z/2)$ is projectively $\Q$-algebraic, we say that $W$ has \emph{projectively $\Q$-algebraic homology}. $\sqbullet$
\end{defn}

\begin{defn}\label{def:Q-bordism}
Let $W$ be a subset of $\R^k$, let $d\in\N$, let $f:M\to W$ be a smooth map from a $d$-dimensional compact smooth manifold $M$ to $W$, and let $\mk{f}=[f:M\to W]\in\mk{N}_d(W)$ be the unoriented bordism class of $f$. We say that $\mk{f}$ is \emph{projectively $\Q$-algebraic} if $\mk{f}$ is represented by (or, equivalently, if $f$ is bordant~to) a $\Q$-regular map from a $d$-dimensional projectively $\Q$-closed $\Q$-nonsingular $\Q$-algebraic set to~$W$. More explicitly, $\mk{f}$ is projectively $\Q$-algebraic if there exist a $d+1$-dimensional compact smooth manifold $T$ with boundary $\partial T$, a $d$-dimensional projectively $\Q$-closed $\Q$-nonsingular $\Q$-algebraic set $P\subset\R^m$, a $\Q$-regular map $g:P\to W$, a smooth diffeomorphism $\psi:M\sqcup P\to\partial T$ and a smooth map $F:T\rightarrow W$ such that $F(\psi(x))=f(x)$ for all $x\in M$, and $F(\psi(x))=g(x)$ for all $x\in P$. If for all $d\in\N$ every element of $\mk{N}_d(W)$ is projectively $\Q$-algebraic, we say that $W$ has \emph{projectively $\Q$-algebraic unoriented bordism}. $\sqbullet$
\end{defn}

If $S$ is any non-empty set, then $\mr{id}_S:S\to S$ denotes the identity map on $S$.

\begin{remark}\label{rem:244}
Let $M$ and $P$ be two disjoint cobordant $d$-dimensional compact smooth submani\-folds of $\R^n$. Identifying $\R^n$ with the subset $\R^n\times\{0\}$ of $\R^n\times\R^m=\R^{n+m}$ for some sufficiently large $m\in\N$ if necessary, we can increase the ambient dimension $n$ in such a way that $n+1\geq2(d+1)+1$. Let $T$ be a $(d+1)$-dimensional compact smooth manifold with boundary and let $\psi:M\sqcup P\to\partial T$ be a smooth diffeomorphism. We can embed and double the cobordism $T$ between $M$ and $P$ following the argument used in \cite[\S{\it b}), pp.176-177]{togn:algmodel}. By the collaring theorem (see \cite[{Thm.6.1, p.113}]{hirsch:difftop}), there exist an open neighborhood $U$ of $\partial T$ in $T$ and a smooth diffeomorphism $\Phi:U\to \partial T\times[0,1)$ such that $\Phi(p)=(p,0)$ for all $p\in\partial T$. Let $\phi:U\to(M\sqcup P)\times[0,1)$ be the smooth diffeomorphism $\phi:=(\psi^{-1}\times\mr{id}_{[0,1)})\circ\Phi$. Observe that $\phi(p)=(\psi^{-1}(p),0)$ for all $p\in\partial T$. Define $A:=T\setminus\partial T$, $B:=\phi^{-1}((M\sqcup P)\times(0,\frac{1}{2}])\subset A$, $N:=\R^n\times(0,+\infty)$ and the smooth map $\theta:B\to N$ by $\theta(x,x_{n+1}):=\phi(x,x_{n+1})$. Since $n+1\geq2(d+1)+1$ and {since} Tietze's theorem ensures the existence of a continuous extension of $\theta$ from $A$ to $N$, we can apply to $\theta$ the extension theorem \cite[Thm.5$(\mr{f})$]{whit:cobordism} obtaining a smooth embedding $\Theta:A\to N$ extending $\theta$. Let $R:\R^{n+1}\to\R^{n+1}$ be the reflection $R(x,x_{n+1}):=(x,-x_{n+1})$ and let $S'$ be the compact smooth submanifold $\Theta(A)\sqcup((M\sqcup P)\times\{0\})\sqcup R(\Theta(A))$ of $\R^{n+1}$. {By} compactness of $T$, {there} exists $\epsilon>0$ such that $S'\cap(\R^n\times(-\epsilon,\epsilon))=(M\sqcup P)\times(-\epsilon,\epsilon)$. Let $L:\R^{n+1}\to\R^{n+1}$ be the linear isomorphism $L(x,x_{n+1}):=(x,\epsilon^{-1}x_{n+1})$. The compact smooth submanifold $S:=L(S')$ of $\R^{n+1}$ has the following property:
$S\cap(\R^n\times(-1,1))=(M\sqcup P)\times(-1,1)$.

In addition, if $F:T\to W$ is a smooth map, i.e., a bordism between the smooth maps $M\to W$, $x\mapsto F(\psi(x))$ and $P\to W$, $x\mapsto F(\psi(x))$, then we can double such a bordism by means of the smooth map $G:S\to W$ defined by $G(x,x_{n+1}):=F(\Theta^{-1}(L^{-1}(x,|x_{n+1}|)))$ if $x_{n+1}\neq0$ and $G(x,0):=F(\psi(x))$. $\sqbullet$
\end{remark}

The following is a $\Q$-version of \cite[Lem.2.5]{ak1981}.

\begin{lem}\label{lem:Q_homology}
Let $W$ be a subset of $\R^k$. Suppose that $W$ admits a structure of finite CW complex (for example, when $W$ is a compact smooth submanifold of $\R^k$). Then the following assertions are equivalent:
\begin{itemize}
 \item[$(\mr{i})$] $W$ has projectively $\Q$-algebraic homology.
 \item[$(\mr{ii})$] $W$ has projectively $\Q$-algebraic unoriented bordism.
\end{itemize}
\end{lem}
\begin{proof}
$(\mr{i})\Longrightarrow(\mr{ii})$ Choose a basis $\{\alpha_{s,i}\}_i$ of each $H_s(W,\Z/2)\neq\{0\}$. For every $\alpha_{s,i}$, let $g_{s,i}:P_{s,i}\to W$ be a $\Q$-regular map from a $s$-dimensional projectively $\Q$-closed $\Q$-nonsingular $\Q$-algebraic set $P_{s,i}\subset\R^{M_s}$ to $W$ such that $\mu_s([g_{s,i}:P_{s,i}\to W])=\alpha_{s,i}$, where $M_s$ is a sufficiently large natural number. For every $t\in\N$, let $\{Y_{t,j}\}_j\subset\R^{N_t}$ be all the possible $t$-dimensional projectively $\Q$-closed $\Q$-nonsingular $\Q$-algebraic set of the form \eqref{eq:Q-milnor}, where $N_t$ is a sufficiently large natural number. By Theorem \ref{thm:generic-projection} and Remark \ref{rmk:N_d}, one can choose $M_s=2s+1$ and $N_t=(2t+1)^2$. For every $s,i,t,j$, define the map $g_{s,i}^{t,j}:V_{s,i}^{t,j}\to W$ by $V_{s,i}^{t,j}:=P_{s,i}\times Y_{t,j}\subset\R^{M_s+N_t}$ and $g_{s,i}^{t,j}(p,y):=g_{s,i}(p)$. By Lemmas \ref{lem:Q-basic}$(\mr{ii})$ and \ref{lem:projective}$(\mr{iv})$, we have that $g_{s,i}^{t,j}$ is $\Q$-regular and $V_{s,i}^{t,j}\subset\R^{M_s+N_t}$ is a $(s+t)$-dimensional projectively $\Q$-closed $\Q$-nonsingular $\Q$-algebraic set.

Let $d\in\N$ and let $L_d\in\N$ be the maximum of the sum $M_s+N_t$ with $s+t=d$. For every $s,i,t,j$ with $s+t=d$, the set $V_{s,i}^{t,j}$ is contained in $\R^{L_d}$ and therefore {it} is a $d$-dimensional projectively $\Q$-closed $\Q$-nonsingular $\Q$-algebraic subset of $\R^{L_d}$ again by Lemma \ref{lem:projective}$(\mr{iv})$. Since these subsets of $\R^{L_d}$ are compact, we can choose vectors $v_{s,i}^{t,j}\in\Q^{L_d}$ such that the sets $\widehat{V}_{s,i}^{t,j}:=V_{s,i}^{t,j}+v_{s,i}^{t,j}$ are pairwise disjoint. Define the map $\widehat{g}_{s,i}^{t,j}:\widehat{V}_{s,i}^{t,j}\to W$ by $\widehat{g}_{s,i}^{t,j}(p):=g_{s,i}^{t,j}(p-v_{s,i}^{t,j})$. Observe that the translation $\R^{L_d}\to\R^{L_d}$, $p\mapsto p-v_{s,i}^{t,j}$ is $\Q$-biregular map. By Lemmas \ref{lem:Q-basic}$(\mr{ii})$ and \ref{lem:projective}$(\mr{iii})$, we deduce that $\widehat{g}_{s,i}^{t,j}$ is $\Q$-regular and $\widehat{V}_{s,i}^{t,j}\subset\R^{L_d}$ is a $d$-dimensional projectively $\Q$-closed $\Q$-nonsingular $\Q$-algebraic set. By Lemma \ref{lem:Q-cobordism} and Theorem \ref{thm:CF}, we have that the set of all the unoriented bordism classes $[\widehat{g}_{s,i}^{t,j}:\widehat{V}_{s,i}^{t,j}\to W]$ with $s+t=d$ is a $\Z/2$-basis of $\mk{N}_d(W)$.

$(\mr{ii})\Longrightarrow(\mr{i})$ Let $\alpha\in H_d(W,\Z/2)$ for some $d\in\N$. By the surjectivity of $\mu_d$, there exist $\mk{f}=[f:M\to W]\in\mk{N}_d(W)$ such that $\mu_d(\mk{f})=\alpha$. Let $P$, $g$, $\psi$, $T$ and $F$ be as in Definition \ref{def:Q-bordism}. Then $\alpha=\mu_d([g:P\to W])$ is projectively $\Q$-algebraic.
\end{proof}

Let $d\in\{0,\ldots,n\}$, let $I_{d,n}:\PP^d(\R)\to\PP^n(\R)$ be the regular embedding $[x_0,\ldots,x_d]\mapsto[x_0,\ldots,x_d,0,\ldots,0]$ and let $J_{d,n}:\G_{d+1,1}\to\G_{n+1,1}$ be the $\Q$-regular embedding defined as follows: if $A=(a_{ij})_{i,j}\in\G_{d+1,1}$, then $J_{d,n}(A)=(b_{ij})_{i,j}$ with $b_{ij}:=a_{ij}$ for all $i,j\in\{0,\ldots,d\}$ and $b_{ij}:=0$ otherwise. It is immediate to verify that the following diagram commutes (see \eqref{eq:mu}):

\begin{center}
\begin{tikzpicture}

\node (A) at (0,0) [] {$\PP^d(\R)$}; 

\node (B) at (2.4,0) [] {$\PP^n(\R)$}; 

\node (A2) at (0,-1.6) [] {$\G_{d+1,1}$}; 

\node (B2) at (2.4,-1.6) [] {$\G_{n+1,1}$}; 

\draw [->] (A) --node[above, ]{\scriptsize $I_{d,n}$} (B);

\draw [->] (A) --node[right, ]{\scriptsize $\eta_d$} (A2);

\draw [->] (B) --node[right, ]{\scriptsize $\eta_n$} (B2);

\draw [->] (A2) --node[above, ]{\scriptsize $J_{d,n}$} (B2);

\node at (1.2,-.7) [] {$\circlearrowleft$};
\end{tikzpicture}
\end{center}
Since $(I_{d,n})_*([\PP^d(\R)])$ is a generator of $H_d(\PP^n(\R),\Z/2)$, we deduce that $(J_{d,n})_*([\G_{d+1,1}])$ is a generator of $H_d(\G_{n+1,1},\Z/2)$. This proves that $\G_{n+1,1}$ has a projectively $\Q$-algebraic homology.

\begin{lem}\label{lem:kunneth}
Let $\ell\in\N^*$, let $n_1,\ldots,n_\ell\in\N^*$ and let $W:=\G_{n_1+1,1}\times\ldots\times\G_{n_\ell+1,1}\subset\R^N$, where $N:=\sum_{i=1}^\ell(n_i^2+1)$. Then $W$ has projectively $\Q$-algebraic unoriented bordism.
\end{lem}
\begin{proof}
Let $d\in\{0,\ldots,N\}$. By Lemma \ref{lem:Q_homology}, it suffices to prove that $H_d(W,\Z/2)$ has a basis whose elements are projectively $\Q$-algebraic. Let $K_d$ be the subset of $\N^\ell$ of all $(d_1,\ldots,d_\ell)$ such that every $d_i$ belongs to $\{0,\ldots,n_i\}$ and $\sum_{i=1}^\ell d_i=d$. For every $D=(d_1,\ldots,d_\ell)\in K_d$, define $\G_D:=\G_{d_1+1,1}\times\ldots\times\G_{d_\ell+1,1}\subset\R^N$ and $g_D:\G_D\to W$ by $g_D:=J_{d_1,n_1}\times\ldots\times J_{d_\ell,n_\ell}$. By Lemmas \ref{lem:projective}$(\mr{iv})$ and \ref{lem:Q-grassmannians}$(\mr{i})$, we deduce that $\G_D\subset\R^N$ is a projectively $\Q$-closed $\Q$-nonsingular $\Q$-algebraic set. In addition, $g_D$ is $\Q$-regular {since} each of its factors $J_{d_i,n_i}$ is {so}. The K\"{u}nneth formula assures that the set $\{(g_D)_*([\G_D])\}_{D\in P_d}$ generates the $\Z/2$-vector space $H_d(W,\Z/2)$. Thus, $W$ has projectively $\Q$-algebraic unoriented bordism, as required.
\end{proof}


\section{$\Q$-algebraic approximations}\label{sec:Q-alg-approx}
In this section, we present $\Q$-versions of some fundamental approximation theorems originated in Nash's article \cite{nash} which were later developed by Wallace \cite{Wa}, Tognoli \cite{togn:algmodel}, Benedetti and Tognoli \cite{BT}, and Akbulut and King \cite{ak1981,akbking:tras}. Here we will use the results of real $\Q$-algebraic geometry obtained in Section \ref{sec:basic-Q-theory} to further develop some of these approximation theorems.

For recent developments in real algebraic approximation theory, we refer the reader to \cite{BK,Be2022,BeWi} and to the references therein.

Let us fix a notation that we will use throughout the rest of the article. Let $M$ be a smooth manifold with possibly non-empty boundary, let $N$ be a smooth manifold (without bounda\-ry) and let $\cinfty(M,N)$ be the set of smooth maps from $M$ to $N$ endowed with the weak $\cinfty$~topolo\-gy, see \cite[Ch.1\;\&\;2]{hirsch:difftop}. Let $f,g\in\cinfty(M,N)$. We say that $g$ is \textit{arbitrarily $\cinfty$ close to $f$} if an arbitrary neighborhood $\mc{U}$ of $f$ in $\cinfty(M,N)$ has been chosen and it holds that $g\in\mc{U}$.


\subsection{$\Q$-stable pairs}

Let $U$ be a non-empty open subset of $\R^n$, let $P$ be a subset of $U$, let $\cinfty(U)$ be the set of all real-valued smooth functions defined on $U$ endowed with the usual ring structure induced by the pointwise addition and multiplication, and let $\mc{I}_U^\infty(P)$ be the ideal of $\cinfty(U)$ of all smooth functions vanishing on $P$. Identify the polynomials in $\R[\x]$ with the corresponding functions in $\cinfty(\R^n)$. Given any set $F\subset\cinfty(\R^n)$, denote {by} $F\cinfty(U)$ the ideal of $\cinfty(U)$ generated by $\{f|_U\}_{f\in F}$.

Let us introduce the concept of $\Q$-stable pair, which is a variant of `variet\`a algebrica di $\R^n$ quasi regolare' introduced by Tognoli in \cite[p.168]{togn:algmodel} and `approximable pair' introduced by Akbulut and King in \cite[p.58]{akbking:tras} (see also the concept of `faithful subvariety' used in \cite[p.211]{BK}). In what follows, if $S$ and $T$ are sets, if $A$, $B$ and $C$ are subsets of $S$ with $\varnothing\neq C\subset A\cap B$ and if $f:A\to T$ and $g:B\to T$ are maps, then we say that $f=g$ on $C$ if $f(x)=g(x)$ for all $x\in C$.

\begin{defn} \label{def:Q-pair}
Let $P$ and $L$ be two subsets of $\R^n$ such that $L\subset P$ and $L$ is $\Q$-algebraic. We say that the pair $(P,L)$ is \emph{stable for $\Q$-approximations}, or \emph{$\Q$-stable} for short, if every function $f\in\II^\infty_{\R^n}(P)$ has the following property: for every $a\in L$, there exist an open neighborhood $U$ of $a$ in $\R^n$, $s\in\N^*$, $g_1,\dots,g_s\in\II_\Q(L)$ and $u_1,\ldots,u_s\in\cinfty(U)$ such that $f=\sum_{i=1}^su_ig_i$ on $U$.

We say that $L$ is \emph{$\Q$-stable} if $(L,L)$ is. $\sqbullet$
\end{defn}

Evidently, if $L=\varnothing$, then $(P,L)$ is always $\Q$-stable. If $L\neq\varnothing$, then in the preceding definition we can always assume that $g_1,\dots,g_s$ are fixed generators of $\II_\Q(L)$ in $\Q[\x]$.

Let $S$ be a subset of $\R^n$, let $r\in\R$, let $\alpha=(\alpha_1,\ldots,\alpha_n)$ be a multi-index in $\N^n$ and let $x\in\R^n$. We denote $\mr{int}_{\R^n}(S)$ the interior of $S$ in $\R^n$, $|r|$ the absolute value of $r$, $|\alpha|$ the order $\sum_{i=1}^n\alpha_i$ of~$\alpha$, $D_\alpha$ the partial derivative operator $\partial^{|\alpha|}/\partial x_1^{\alpha_1}\cdots\partial x_n^{\alpha_n}$, $|x|_n$ the Euclidean norm of $x$ in $\R^n$ and $\B_n(r)$ the open ball $\{x\in\R^n:|x|_n<r\}$ for every $r>0$.

The following result on smooth functions is known. However, we are not aware of a reference that completely proves its statement. Therefore, for completeness, we provide a proof in Appendix~\ref{appendix-A}.

\begin{lem}\label{lem:diff-diff}
Let $a\in\R^n$, let $d\in\{0,\ldots,n-1\}$, let $f_1,\ldots,f_{n-d}\in\cinfty(\R^n)$ be such that $f_1(a)=\ldots=f_{n-d}(a)=0$ and the gradients $\nabla f_1(a),\ldots,\nabla f_{n-d}(a)$ are linearly independent. Choose an open neighborhood $U$ of $a$ in $\R^n$ and set $V:=\{x\in U:f_1(x)=0,\ldots,f_{n-d}(x)=0\}$. We have:
\begin{itemize}
 \item[$(\mr{i})$] There exists an open neighborhood $U'$ of $a$ in~$U$ with the following property: for every $f\in\II^\infty_U(V)$, there exist $u_1,\ldots,u_{n-d}\in\cinfty(U')$ such that $f=\sum_{i=1}^{n-d}u_if_i$ on $U'$.
 \item[$(\mr{ii})$] Suppose that $d>0$. Let $e\in\N^*$ and let $g_1,\ldots,g_e\in\cinfty(\R^n)$ be such that $e\leq d$, $g_1(a)=\ldots=g_e(a)=0$ and the gradients $\nabla f_1(a),\ldots,\nabla f_{n-d}(a),\nabla g_1(a),\ldots,\nabla g_e(a)$ are linearly independent. Define the function $g\in\cinfty(\R^n)$ and the set $W\subset V$ by $g:=\prod_{i=1}^eg_i$ and $W:=\{x\in V:g(x)=0\}$. Then there exists an open neighborhood $U'$ of $a$ in~$U$ with the following property: for every $f\in\II^\infty_U(W)$, there exist $u_0,u_1,\ldots,u_{n-d}\in\cinfty(U')$ such that $f=u_0g+\sum_{i=1}^{n-d}u_if_i$ on $U'$.
\end{itemize}
\end{lem}

The next result collects four important properties of $\Q$-stable pairs. The second property justifies the nomenclature `$\Q$-stable': if $(P,L)$ is $\Q$-stable with $L$ compact, then every smooth function $f:\R^n\to\R$ vanishing on $P$ can be $\cinfty$ approximated by rational polynomials vanishing on $L$, so the vanishing of $f$ on $P$ is `stable for $\Q$-approximations' on $L$. 

\begin{lem}\label{lem:basic-Q-pair}
Let $P$ and $L$ be two subsets of $\R^n$ such that $\varnothing\neq L\subset P$ and $L$ is $\Q$-algebraic. We have:
\begin{itemize}
 \item[$(\mr{i})$] $(P,L)$ is $\Q$-stable if and only if $\II^\infty_{\R^n}(P)\subset\II_\Q(L)\cinfty(\R^n)$, i.e., if every $f\in\II^\infty_{\R^n}(P)$ can be expressed as $f=\sum_{i=1}^su_ig_i$ on the whole $\R^n$ for some $s\in\N^*$, $g_1,\dots,g_s\in\II_\Q(L)$ and $u_1,\ldots,u_s\in\cinfty(\R^n)$.
 \item[$(\mr{ii})$] If $L$ is compact in $\R^n$ and $(P,L)$ is $\Q$-stable, then for every $f\in\II^\infty_{\R^n}(P)$, for every $\varepsilon>0$, for every $h\in\N$ and for every compact neighborhood $K$ of $L$ in $\R^n$, there exists a polynomial $g\in\II_\Q(L)$ such that
$$
\textstyle
\text{$\max_{x\in K}|D_\alpha f(x)-D_\alpha g(x)|<\varepsilon\;$ for all $\alpha\in\N^n$ with $|\alpha|\leq h$.}
$$
 \item[$(\mr{iii})$] If $L\setminus\mr{int}_{\R^n}(P)\subset\Reg^\Q(L)$, then $(P,L)$ is $\Q$-stable. In particular, every disjoint union of finitely many $\Q$-nonsingular $\Q$-algebraic subsets of $\R^n$ of possibly different dimensions is $\Q$-stable.
 \item[$(\mr{iv})$] If $L\subset\R^n$ is a $\Q$-stable $\Q$-algebraic set and $h\in\N^*$, then $L\times\{0\}\subset\R^n\times\R^h=\R^{n+h}$ is also a $\Q$-stable $\Q$-algebraic set.
\end{itemize}
\end{lem}
\begin{proof}
Pick generators $g_1,\ldots,g_s$ of $\II_\Q(L)$ in $\Q[\x]$.

$(\mr{i})$ The `if' implication is evident. Let us prove the `only if' implication. Suppose that $(P,L)$ is $\Q$-stable. Let $f\in\II^\infty_{\R^n}(P)$. By Definition~\ref{def:Q-pair} and by the paracompactness of $\R^n$, there exist a locally finite  countable family $\{U_j\}_{j\in\N^*}$ of non-empty open subsets of $\R^n$ such that $L\subset\bigcup_{j\in\N^*}U_j$ and, for every $j\in\N^*$, smooth functions $u_{j1},\ldots,u_{js}\in\cinfty(U_j)$ such that $f=\sum_{i=1}^su_{ji}g_i$ on~$U_j$. Let $U_0:=\R^n\setminus L$ and let $u_{0i}:=fg_i(\sum_{i=1}^sg_i^2)^{-1}\in\cinfty(U_0)$ for every $i\in\{1,\ldots,s\}$ (so $f=\sum_{i=1}^su_{0i}g_i$ on $U_0$). Pick a smooth partition of the unity $\{\varphi_j\}_{j\in\N}$ on $\R^n$ subordinate to $\{U_j\}_{j\in\N}$ (i.e., $\mr{supp}(\varphi_j)\subset U_j$) and, for every $i\in\{1,\ldots,s\}$, define $u_i:=\sum_{j\in\N}\varphi_ju_{ji}\in\cinfty(\R^n)$, where $\varphi_ju_{ji}$ is assumed to be null outside $U_j$. It follows that $f=\sum_{i=1}^su_ig_i$ on $\R^n$.

$(\mr{ii})$ By $(\mr{i})$, there exist $u_1,\ldots,u_s\in\cinfty(\R^n)$ such that $f=\sum_{i=1}^su_ig_i$ on $\R^n$. By the Weierstrass approximation theorem and by the density of $\Q$ in $\R$, there exists a polynomial $v_i\in\Q[\x]$ arbitrarily $\cinfty$ close to $u_i$ on $\R^n$. It is now sufficient to set $g:=\sum_{i=1}^sv_ig_i$.

$(\mr{iii})$ Let $f\in\II^\infty_{\R^n}(P)$. Set $d:=\dim(L)$. If $d=n$, then $P=L=\R^n$ so it is evident that $(P,L)$ is $\Q$-stable. Suppose $d<n$. We have to prove that $f$ has the property stated in Definition~\ref{def:Q-pair}. Pick a point $a\in L$. If $a\in\mr{int}_{\R^n}(P)$, then it suffices to set $U:=\mr{int}_{\R^n}(P)$ and every $u_i:\equiv0$. Suppose that $a\in L\setminus\mr{int}_{\R^n}(P)$. Since $L\setminus\mr{int}_{\R^n}(P)\subset\Reg^\Q(L)$, there exist $j_1,\ldots,j_{n-d}\in\{1,\ldots,s\}$ and an open neighborhood $U$ of $a$ in $\R^n$ such that $\nabla g_{j_1}(a),\ldots,\nabla g_{j_{n-d}}(a)$ are linearly independent and $L\cap U=\ZZ_\R(g_{j_1},\ldots,g_{j_{n-d}})\cap U$ (so $s\geq n-d$). Rearranging the indices, we {may} assume $j_i=i$. By Lemma \ref{lem:diff-diff}$(\mr{i})$, shrinking $U$ around $a$ if necessary, there exist $u_1,\ldots,u_{n-d}\in\cinfty(U)$ such that $f=\sum_{i=1}^{n-d}u_ig_i$ on $U$, so $f=\sum_{i=1}^nu_ig_i$ on $U$ with $u_i:\equiv0$ for every $i\in\{n-d+1,\ldots,s\}$ if $s>n-d$. This proves that, if $L\setminus\mr{int}_{\R^n}(P)\subset\Reg^\Q(L)$, then $(P,L)$ is $\Q$-stable. In particular, if $L\subset\R^n$ is $\Q$-nonsingular, then $L$ is $\Q$-stable. Since the concept of $\Q$-stable pair is local, we immediately deduce that the disjoint union of finitely many $\Q$-nonsingular $\Q$-algebraic subsets of $\R^n$ of possibly different dimensions is also $\Q$-stable.

$(\mr{iv})$ Let $(x,y):=(x_1,\ldots,x_n,y_1,\ldots,y_h)$ be the coordinates of $\R^{n+h}$, let $f\in\II^\infty_{\R^{n+h}}(L\times\{0\})$, let $a\in L$, let $U$ be an open neighborhood of $a$ in $\R^n$ and let $u_1,\ldots,u_s\in\cinfty(U)$ be such that $f(x,0)=\sum_{i=1}^su_i(x)g_i(x)$ for all $x\in U$. By Lemma \ref{lem:diff-diff}$(\mr{i})$, shrinking $U$ around $a$ if necessary, there exist $\epsilon>0$ and $v_1,\ldots,v_h\in\cinfty(U\times\B_h(\epsilon))$ such that $f(x,y)-f(x,0)=\sum_{j=1}^hv_j(x,y)y_j$ for all $(x,y)\in U\times\B_h(\epsilon)$. It follows that $f(x,y)=\sum_{i=1}^su_i(x)g_i(x)+\sum_{j=1}^hv_j(x,y)y_j$ for all $(x,y)\in U\times\B_h(\epsilon)$. Since all polynomials $g_i(\x)$ and $\y_j$ belong to $\II_\Q(L\times\{0\})$, the proof is complete.
\end{proof}

Let us introduce the concept of $\Q$-nonsingular $\Q$-algebraic hypersurface.

\begin{defn}
Let $X\subset\R^n$ be a $\Q$-nonsingular $\Q$-algebraic set and let $Y$ be a subset of $X$. We say that $Y$ is a \emph{$\Q$-nonsingular $\Q$-algebraic hypersurface of $X$} if $Y\subset\R^n$ is a $\Q$-nonsingular $\Q$-algebraic set of dimension $\dim(X)-1$. $\sqbullet$
\end{defn}

The next result provides possibly singular examples of $\Q$-stable $\Q$-algebraic sets.

\begin{lem}\label{lem:35}
Let $X\subset\R^n$ be a $\Q$-nonsingular $\Q$-algebraic set, let $\ell\in\N^*$ and let $Y_1,\ldots,Y_\ell$ be $\Q$-nonsingular $\Q$-algebraic hypersurfaces of $X$. Suppose that $Y_1,\ldots,Y_\ell$ are in general position in $X$ (in the usual smooth sense). Then $\bigcup_{i=1}^\ell Y_i\subset\R^n$ is a $\Q$-stable $\Q$-algebraic set.
\end{lem}
\begin{proof}
Let $d:=\dim(X)$, let $L:=\bigcup_{i=1}^\ell Y_i$, let $f\in\II^\infty_{\R^n}(L)$ and let $a\in L$. Suppose that $d<n$. If $d=n$, the proof is similar. By hypothesis, there exist polyno\-mials $f_1,\ldots,f_{n-d},g_1,\ldots,g_e\in\Q[\x]$ for some $e\in\N^*$ with $e\leq\ell$, indices $j_1,\ldots,j_e\in\{1,\ldots,\ell\}$ with $j_1<\ldots<j_e$ and an open neighborhood $U$ of $a$ in $\R^n$ such that the polynomials $f_i$ vanish on~$X$, $g_1(a)=\ldots=g_e(a)=0$, the gradients $\nabla f_1(a),\ldots,\nabla f_{n-d}(a),\nabla g_1(a),\ldots,\nabla g_e(a)$ are linearly independent, $X\cap U=\ZZ_\R(f_1,\ldots,f_{n-d})\cap U$, $Y_{j_k}\cap U=X\cap U\cap\ZZ_\R(g_k)$ for all $k\in\{1,\ldots,e\}$ and $Y'\cap U=\varnothing$, where $Y':=\bigcup_{j\in\{1,\ldots,\ell\}\setminus\{j_1,\ldots,j_e\}}Y_j$. In particular, $L\cap U=\bigcup_{k=1}^e(X\cap U\cap\ZZ_\R(g_k))$. By Lemma \ref{lem:diff-diff}$(\mr{ii})$, shrinking $U$ around $a$ if necessary, there exist $u_0,u_1,\ldots,u_{n-d}\in\cinfty(U)$ such that $f=u_0g+\sum_{i=1}^{n-d}u_if_i$ on $U$, where $g:=\prod_{k=1}^eg_k$. Choose a polynomial $p\in\Q[\x]$ such that $\ZZ_\R(p)$ is equal to the union of $Y'$ and of all the $\Q$-irreducible components of $Y_{j_k}$ that do not intersect $U$ for all $k\in\{1,\ldots,e\}$. Observe that $p$ never vanishes on $U$, $f=\frac{u_0}{p}(pg)+\sum_{i=1}^{n-d}u_if_i$ on~$U$ and $f_1,\ldots,f_{n-d}\in\II_\Q(L)$. The polynomial $pg$ belongs to $\II_\Q(L)$ as well, because if a $\Q$-irreducible component $Y'_{j_k}$ of some $Y_{j_k}$ intersects $U$ then $\ZZ_\R(g)\supset Y'_{j_k}\cap U$ so $\ZZ_\R(g)\supset Y'_{j_k}$ by Lemma \ref{dimirred} and Theorem \ref{dimension}. This proves that $L\subset\R^n$ is $\Q$-stable, as required.
\end{proof}

Recall that $\R^n$ is identified with the vector subspace $\R^n\times\{0\}$ of $\R^{n+1}=\R^n\times\R$.

\begin{lem}\label{lem:L}
Let $L\subset\R^n$ be a compact $\Q$-stable $\Q$-algebraic set, let $K$ be a compact neighborhood of $L$ in $\R^n$ and let $f\in\II^\infty_{\R^n}(L\cup(\R^n\setminus K))$. Then, for every $\epsilon>0$, for every $h\in\N$ and for every non-empty compact subset $K'$ of $\R^n$, there exists a $\Q$-regular function $g\in\reg^\Q(\R^n)$ with the following properties:
\begin{itemize}
 \item[$(\mr{i})$] There exist $e\in\N$ and $p\in\II_\Q(L)$ such that $\deg(p)\leq 2e$ and $g(x)=p(x)(1+|x|_n^2)^{-e}$ for all $x\in\R^n$.
 \item[$(\mr{ii})$] $\sup_{x\in\R^n}|f(x)-g(x)|<\epsilon$.
 \item[$(\mr{iii})$] $\max_{x\in K'}|D_\alpha f(x)-D_\alpha g(x)|<\epsilon$ for all $\alpha\in\N^n$ with $|\alpha|\leq h$.
\end{itemize}
\end{lem}

\begin{proof}
We improve the strategy used in the proofs of \cite[Lem.2.1\;\&\;2.2]{ak1981} and \cite[Lem.2.8.1]{akbking:tras}, see also \cite[Lem.4]{togn:algmodel}. Suppose that $L\neq\varnothing$. If $L=\varnothing$, the proof is similar but simpler. Let $\sph^n$ be the standard unit sphere of $\R^{n+1}=\R^n\times\R$, let $N=(0,\ldots,0,1)$ be its north pole and let $\theta:\sph^n\setminus\{N\}\to\R^n$ be its stereographic projection from $N$. Recall that $\theta(x,x_{n+1})=x(1-x_{n+1})^{-1}$ and $\theta^{-1}(x)=(2x,-1+|x|_n^2)(1+|x|_n^2)^{-1}$. Define $d:=\dim(L)$. Observe that $d<n$ {since} $L$ is compact and {therefore} $L\neq\R^n$. Choose generators $g_1,\ldots,g_s$ of $\II_\Q(L)$ in $\Q[\x]$ and, for every $i\in\{1,\ldots,s\}$, write $g_i$ as follows: $g_i=\sum_{j=0}^{d_i}g_{ij}$, where $d_i:=\deg(g_i)$ and every $g_{ij}$ is a homogeneous polynomial in $\Q[\x]$ of degree $j$. Since $L\neq\varnothing$, we can assume that every $d_i$ is positive. Define the polynomial $G_i\in\Q[\x,\x_{n+1}]$ by $G_i(\x,\x_{n+1}):=\sum_{j=0}^{d_i}(1-\x_{n+1})^{d_i-j}g_{ij}(\x)$. We have:
\begin{equation}\label{eq:u_i}
G_i(N)=g_{id_i}(0)=0
\end{equation}
and
\begin{equation}\label{eq:p_i}
\text{$(g_i\circ\theta)(x,x_{n+1})=(1-x_{n+1})^{-d_i}G_i(x,x_{n+1})\;$ for all $(x,x_{n+1})\in\sph^n\setminus\{N\}$.}
\end{equation}

Let $H$ be a compact neighborhood of $N$ in $\R^{n+1}$ such that $H\cap\theta^{-1}(K)=\varnothing$, let $L':=\{N\}\cup\theta^{-1}(L)$, let $P':=\mr{int}_{\R^{n+1}}(H)\cup L'=\mr{int}_{\R^{n+1}}(H)\sqcup\theta^{-1}(L)$, let $G_{s+1}\in\Q[\x,\x_{n+1}]$ be the polynomial defined by $G_{s+1}(\x,\x_{n+1}):=-1+\sum_{i=1}^{n+1}\x_i^2$ and let $\rho:\R^{n+1}\setminus\{0\}\to\sph^n$ be the standard smooth retraction $\rho(x,x_{n+1}):=(x,x_{n+1})(|x|_n^2+x_{n+1}^2)^{-1/2}$. 

Let us show that $(P',L')$ is a $\Q$-stable pair of $\R^{n+1}$. By \eqref{eq:u_i} and \eqref{eq:p_i}, we know that $L'=\ZZ_\R(G_1,\ldots,G_{s+1})$, so $L'\subset\R^{n+1}$ is $\Q$-algebraic. Let $q'\in\II^\infty_{\R^{n+1}}(P')$ and let $a\in L'$. We have to prove that $q'$ has the property stated in Definition~\ref{def:Q-pair}. If $a\in\mr{int}_{\R^{n+1}}(H)$, this is evident. Suppose that $a\in\theta^{-1}(L)$. Let $b:=\theta(a)\in L$ and let $q\in\cinfty(\R^n)$ be the smooth function $q(x):=q'(\theta^{-1}(x))$. {Observe that $q\in\II^\infty_{\R^n}(L)$.} Since $L\subset\R^n$ is $\Q$-stable, there exist an open neighborhood $U$ of $b$ in $\R^n$ and $u_1,\ldots,u_s\in\cinfty(U)$ such that $q=\sum_{i=1}^su_ig_i$ on $U$, so $q'(x,x_{n+1})=q(\theta(x,x_{n+1}))=\sum_{i=1}^s u_i(\theta(x,x_{n+1}))(1-x_{n+1})^{-d_i}G_i(x,x_{n+1})$ for all $(x,x_{n+1})\in\theta^{-1}(U)$. Define the open neighborhood $U':=\{(x,x_{n+1})\in\R^{n+1}\setminus\{0\}:\rho(x,x_{n+1})\in\theta^{-1}(U),x_{n+1}\neq1\}$ of $a$ in $\R^{n+1}$, the function $u'_i\in\cinfty(U')$ by $u'_i(x,x_{n+1}):=u_i(\theta(\rho(x,x_{n+1})))(1-x_{n+1})^{-d_i}$ for every $i\in\{1,\ldots,s\}$ and the function $Q'\in\cinfty(U')$ by $Q':=\sum_{i=1}^su'_iG_i|_{U'}$. Since $q'=Q'$ on $\theta^{-1}(U)$, by Lemma \ref{lem:diff-diff}$(\mr{i})$, shrinking $U'$ around $a$, there exists $u'_{s+1}\in\cinfty(U')$ such that $q'=Q'+u'_{s+1}G_{s+1}$ on $U'$. It follows that $q'=\sum_{i=1}^{s+1}u'_iG_i$ on $U'$. This proves that the pair $(P',L')$ is $\Q$-stable.

Let us construct a smooth extension $F\in\cinfty(\R^{n+1})$ of the composition $f\circ\theta:\sph^n\setminus\{N\}\to\R$. Let $\psi\in\cinfty(\R)$ be such that $\psi\geq0$ on $\R$, $\mr{supp}(\psi)\subset\big[\frac{1}{2},\frac{3}{2}\big]$ and $\psi(1)=1$, and let $N_+$ be the half line $\{tN\in\R^{n+1}\,|\,t\geq0\}$. Define the extension $F:\R^{n+1}\to\R$ of $f\circ\theta$ by $F(x,x_{n+1}):=\psi(|x|_n^2+x_{n+1}^2)f(\theta(\rho(x,x_{n+1})))$ if $(x,x_{n+1})\not\in N_+$ and $F(x,x_{n+1}):=0$ if $(x,x_{n+1})\in N_+$. Since $\mr{supp}(f)$ is contained in the compact subset $K$ of $\R^n$ by hypothesis, we have that $F$ is smooth.

By Lemma \ref{lem:basic-Q-pair}$(\mr{ii})$, there exists $p'\in\II_\Q(L')$ arbitrarily $\cinfty$ close to $F$. In particular, we can assume that $|F-p'|<\veps$ on $\sph^n$, so $|f(x)-p'(\theta^{-1}(x))|=|F(\theta^{-1}(x))-p'(\theta^{-1}(x))|<\epsilon$ for all $x\in\R^n$. Write $p'$ as follows: $p'=\sum_{i=0}^ep'_i$, where $e:=\deg(p')$ and every $p'_i$ is a homogeneous polynomial in $\Q[\x,\x_{n+1}]$ of degree $i$. Define $p(x):=\sum_{i=0}^e(1+|x|^2_n)^{e-i}p'_i(2x,-1+|x|_n^2)\in\Q[\x]$ and $g\in\Reg^\Q(\R^n)$ by $g(x):=p(x)(1+|x|_n^2)^{-e}$. Observe that $\deg(p)\leq2e$. Since $g(x)=p'(\theta^{-1}(x))$ for all $x\in\R^n$, we have that $p\in\II_\Q(L)$ and $|f-g|<\epsilon$ on $\R^n$. Finally, consider the smooth map $\Psi:\R^n\to\R^{n+1}$ defined by $\Psi(x):=\theta^{-1}(x)$. Since the pullback map $\Psi^*:\cinfty(\R^{n+1})\to\cinfty(\R^n)$ is continuous, we can choose $\Psi^*(p')=p'\circ\Psi=g$ arbitrarily $\cinfty$ close to $\Psi^*(F)=F\circ\Psi=f$.
\end{proof}


\subsection{Relative $\Q$-algebraic approximations and proof of Theorem \ref{thm:NTQ}}

Let $X$ be a subset of $\R^n$. We denote $\cl_{\R^n}(X)$ the Euclidean closure of $X$ in $\R^n$. Suppose that $X$ is open in $\R^n$ and set $\overline{X}:=\cl_{\R^n}(X)$. Let $U$ be an open subset of $\R^n$. Suppose that $\overline{X}\subset U$. We say that $\overline{X}$ is a \emph{smooth submanifold with boundary of $U$} if the topological boundary $\overline{X}\setminus X$ of $\overline{X}$ in $U$ is a smooth hypersurface of~$\R^n$. Observe that, if $L$ is a compact subset of $U$, then there always exists an open neighborhood $V$ of $L$ in $U$ such that $\overline{V}:=\cl_{\R^n}(V)$ is a compact smooth submanifold with boundary of $U$: it suffices to consider a non-negative smooth function $f:\R^n\to\R$ such that $f^{-1}(0)=L$ and $f\equiv1$ out of a compact neighborhood of $L$ in $U$, and define $V:=f^{-1}([0,\epsilon))$, where $\epsilon$ is a regular value of $f$ belonging to $(0,1)$. As we mentioned at the beginning of Section~\ref{sec:Q-alg-approx}, if $\overline{V}$ is a compact smooth submanifold with boundary of $U$, then it makes sense to speak about smooth maps defined on $\overline{V}$ with values in any smooth manifold $N$ and therefore consider the corresponding set $\cinfty(\overline{V},N)$ of smooth maps equipped with the weak $\cinfty$ topology. This will help {to} simplify the statements and the proofs of the following approximation results.

The following is a `$\Q$-nonsingular' version of \cite[Thm.2.8.3, p.63]{akbking:tras}.

\begin{lem}\label{lem:283}
Let $L\subset\R^n$ be a non-empty compact $\Q$-stable $\Q$-algebraic set, let $U$ be an open neighborhood of $L$ in $\R^n$, let $W\subset\R^k$ be a $\Q$-nonsingular $\Q$-algebraic set and let $f:U \to W$ be a smooth map such that the restriction $f|_L:L\to W$ of $f$ to $L$ is $\Q$-regular. Choose an open neighborhood $V$ of $L$ in $\R^n$ such that $\overline{V}:=\cl_{\R^n}(V)$ is a compact smooth submanifold with boundary of $U$. Then there exist a smooth map $v:\overline{V}\to\R^k$, a $\Q$-algebraic subset $Z$ of $\R^{n+k}=\R^n\times\R^k$ of dimension~$n$ and a $\Q$-regular map $\eta:Z \to W$ with the following properties:
\begin{itemize}
 \item[$(\mr{i})$] $\widehat{V}:=\{(x,v(x))\in\R^{n+k}:x\in V\}$ is an open subset of $\Reg^\Q(Z)$.
 \item[$(\mr{ii})$] $v(x)=0$ and $\eta(x,0)=f(x)$ for all $x\in L$. In particular, we have $L\times\{0\}\subset\widehat{V}$.
 \item[$(\mr{iii})$] $v$ can be chosen arbitrarily $\cinfty$ close to the zero map $\overline{V}\to\R^k$, $x\mapsto0$, and $\eta$ can be chosen such that, if $\widehat{v}:\overline{V}\to Z$ denotes the smooth map $\widehat{v}(x):=(x,v(x))$, then $\eta\circ\widehat{v}$ is arbitrarily $\cinfty$ close to $f|_{\overline{V}}$.
\end{itemize}
\end{lem}
\begin{proof}
Let $V_1$ and $V_2$ be relatively compact open neighborhoods of $L$ in $\R^n$ such that $\overline{V}\subset V_1$, $\overline{V_1}:=\cl_{\R^n}(V_1)\subset V_2$ and $\cl_{\R^n}(V_2)\subset U$. We can assume that $\overline{V_1}$ is a smooth submanifold with boundary of $U$. Consider the compact subset $f(\overline{V_1})$ of $W$ and choose a relatively compact open neighborhood $W'$ of $f(\overline{V_1})$ in $W$. Let $N:=\{(w,y) \in W \times \R^k:y\in T_w(W)^\perp\}$ be the total space of the normal bundle of $W$ in~$\R^k$, let $\pi:N\to W$ be the bundle projection $(w,y)\mapsto w$ and let $\theta:N\to \R^k$ be the smooth map $\theta(w,y):=w+y$. By the Inverse Function Theorem, there exists $\delta>0$ such that $T:=\theta(N\cap(W'\times\B_k(\delta)))$ is an open neighborhood of $W'$ in $\R^k$ and the restriction $\theta':N\cap(W'\times\B_k(\delta))\to T$ of $\theta$ is a smooth diffeomorphism. Thus, $T$ is a tubular neighborhood of the smooth submanifold $W'$ in $\R^k$ and the smooth map $\rho:T\to W'$, $y\mapsto\pi((\theta')^{-1}(y))$ is the corresponding closest point map.

By Lemma~\ref{lem:Q-basic}$(\mr{iii})$, there exists a $\Q$-regular map $f^*:\R^n\to\R^k$ such that $f^*(x)=f(x)$ for all $x\in L$. Let $\psi:\R^n\to\R$ be a smooth function such that $\psi\equiv1$ on $V_2$ and $\mr{supp}(\psi)$ is a compact subset of $U$. Define the smooth map $\widetilde{f}:\R^n\to\R^k$ by $\widetilde{f}(x):=\psi(x)(f(x)-f^*(x))$ for all $x\in U$ and $\widetilde{f}(x):=0$ for all $x\in\R^n\setminus U$. Observe that $\widetilde{f}(x)=f(x)-f^*(x)$ for all $x\in V_2$ and $\widetilde{f}(x)=0$ for all $x\in L\cup(\R^n\setminus\mr{supp}(\psi))$. Applying Lemma \ref{lem:L} to every component of $\widetilde{f}$, we obtain a $\Q$-regular map $\widetilde{g}:\R^n\to\R^k$ such that $\widetilde{g}=0$ on~$L$ and $\widetilde{g}$ is arbitrarily $\cinfty$ close to $\widetilde{f}$. Define the $\Q$-regular map $g:\R^n\to\R^k$ by $g:=\widetilde{g}+f^*$. Observe that $g=f$ on $L$ and $g|_{V_2}$ is arbitrarily $\cinfty$ close to $f|_{V_2}$.

Let~$F:\overline{V_1}\times\B_k(\delta)\to\R^{2k}$ be the smooth map $F(x,y):=(f(x)+y,y)$. We claim that $F$ is transverse to $N$ in $\R^{2k}$. Indeed, if $(x,y)\in\overline{V_1}\times\B_k(\delta)$ and $(f(x)+y,y)\in N$, then $f(x)\in W'\subset T$, $f(x)+y\in W$ and $y\in T_{f(x)+y}(W)^\perp$, so $y=0$ by the construction of $T$. In addition, we have: $T_{(f(x),0)}(N)=T_{f(x)}(W)\times T_{f(x)}(W)^\perp$, $dF_{(x,0)}(\R^{n+k})$ contains the diagonal $\Delta$ of $\R^{2k}=\R^k\times\R^k$, $\Delta\cap(T_{f(x)}(W)\times T_{f(x)}^\perp(W))=\{(0,0)\}$ and therefore $dF_{(x,0)}(\R^{n+k})+T_{(f(x),0)}(N)=\R^{2k}$. This proves that $F$ is transverse to $N$ in $\R^{2k}$, as claimed.

Choose $g$ such that $g|_{V_2}$ is sufficiently $\cinfty$ close to $f|_{V_2}$ so that $g(V_1)\subset g(\overline{V_1})\subset T$, and define the smooth map $v_1:V_1\to\R^k$ by $v_1(x):=\rho(g(x))-g(x)$ and $\widehat{V}_1:=\{(x,v_1(x))\in\R^{n+k}:x\in V_1\}$. By \cite[Thm.2.1(b),~p.74]{hirsch:difftop}, if $G:\R^{n+k}\to\R^{2k}$ is the $\Q$-regular map $G(x,y):=(g(x)+y,y)$, we can assume that $G|_{\overline{V_1}\times\B_k(\delta)}$ is arbitrarily $\cinfty$ close to $F$ and so transverse to $N$ in $\R^{2k}$. Define the $\Q$-algebraic set $Z_1\subset\R^{n+k}$ by $Z_1:=G^{-1}(N)$. The construction of $T$ assures that $Z_1\cap(V_1\times\B_k(\delta))=\widehat{V}_1$. Indeed, if $x\in V_1$, then $g(x)\in T$ and there exists a unique $y\in\B_k(\delta)$ such that $g(x)+y\in W$ and $y\in T_{g(x)+y}(W)^\perp$, which is $y=v_1(x)$. Since $G|_{V_1\times\B_k(\delta)}$ is transverse to $N$ in $\R^{2k}$, Proposition \ref{prop:Q-transverse} assures that $\widehat{V}_1\subset\Reg^\Q(Z_1,n)$. Denote $Z\subset\R^{n+k}$ the union of all $\Q$-irreducible components of $Z_1$ of dimension $n$, so $Z\subset\R^{n+k}$ is $\Q$-algebraic and $\dim(Z)=n$. By Proposition \ref{prop211}$\mr{(i)}$, we have that $\widehat{V}_1=Z\cap(V_1\times\B_k(\delta))\subset\Reg^\Q(Z)$. Thus, $\widehat{V}_1=Z\cap(V_1\times\B_k(\delta))$ is an open subset of $\Reg^\Q(Z)$.

Let $v:V\to\R^k$ be the restriction $v:=v_1|_V$ of $v_1$ to $V$ and let $\widehat{V}:=\{(x,v(x))\in\R^{n+k}:x\in V\}$. Since $\widehat{V}=\widehat{V}_1\cap(V\times\R^k)$, we deduce that $\widehat{V}$ is an open subset of $\Reg^\Q(Z)$.

Define the $\Q$-regular map $\eta:Z\to W$ by $\eta(x,y):=g(x)+y$. Since $g=f$ on $L$, we have that $v|_L=0$ so $L\times\{0\}\subset \widehat{V}$ and $\eta(x,0)=g(x)=f(x)$ for all $x\in L$.

We have just proved that $(\mr{i})$ and $(\mr{ii})$ are satisfied. It remains to show $(\mr{iii})$. Consider the restrictions $f|_{V_1}$ and $g|_{V_1}$ as maps in $\cinfty(V_1,T)$, and $\rho$ as a map in $\cinfty(T,W)$ in the natural way. Since the push-forward map $\rho_*:\cinfty(V_1,T)\to\cinfty(V_1,W)$ is continuous and $g$ can be chosen so that $g|_{V_1}$ is arbitrarily $\cinfty$ close to $f|_{V_1}$, we can assume that $\rho_*(g|_{V_1})$ is arbitrarily $\cinfty$ close to $\rho_*(f|_{V_1})=f|_{V_1}$. Now, if $\widehat{v}:\overline{V}\to Z$ denotes the smooth map $\widehat{v}(x):=(x,v(x))$, then $\eta\circ\widehat{v}=\rho_*(g|_{V_1})|_{\overline{V}}$ {as} $\eta(x,v(x))=\rho(g(x))$ for all $x\in\overline{V}$. This proves $(\mr{iii})$.
\end{proof}

\begin{remark}\label{rem:vvarnothing}
If $L=\varnothing$, the statement of Lemma \ref{lem:283} remains valid: it suffices to cancel the condition `\emph{the restriction $f|_L:L\to W$ of $f$ to $L$ is $\Q$-regular}' and item $(\mr{ii})$. In this case the proof is similar but simpler. $\sqbullet$
\end{remark}

Let $M$ be a smooth submanifold of $\R^n$ of dimension $d$. Recall that  the normal bundle map $\beta_M:M\to\G_{n,n-d}$ of $M$ in $\R^n$ is defined by $\beta_M(x):=T_x(M)^\perp$. Here we are identifying the vector subspace $T_x(M)^\perp$ of $\R^n$ with the matrix in $\R^{n^2}$ associated to the orthogonal projection of $\R^n$ onto $T_x(X)^\perp$ with respect to the canonical basis of $\R^n$.

Recall that, for every $n,m\in\N$ with $m>n$, we identify $\R^n$ with the subset $\R^n\times\{0\}$ of $\R^n\times\R^{m-n}=\R^m$, so we can write $\R^n\subset\R^m$ and every subset of $\R^n$ is also a subset of $\R^m$.

The following proposition is a strong `$\Q$-nonsingular' version of \cite[Prop.2.8]{ak1981}, see also \cite[Thms.2.8.4, pp.65-66]{akbking:tras}.

\begin{thm}\label{thm:Q_tognoli}
Let $M$ be a compact smooth submanifold of $\R^n$ of dimension $d$, let $L$ be a non-empty projectively $\Q$-closed $\Q$-stable $\Q$-algebraic $\Q$-algebraic subset of $\R^n$ with $L\subset M$, let $\beta_M:M\to\G_{n,n-d}$ be the normal bundle map of $M$ in $\R^n$, let $W\subset\R^k$ be a $\Q$-nonsingular $\Q$-algebraic set and let $f:M \to W$ be a smooth map. Suppose that the unoriented bordism class of $f$ is projectively $\Q$-algebraic, and the restrictions $\beta_M|_L:L\to\G_{n,n-d}$ and $f|_L:L\to W$ are $\Q$-regular. Then there exist a projectively $\Q$-closed $\Q$-nonsingular $\Q$-algebraic set $M'\subset\R^{n+t}=\R^n\times\R^t$ for some $t\in\N^*$, a smooth diffeomorphism $\phi:M\to M'$ and a $\Q$-regular map $f':M'\to W$ such~that:
\begin{itemize}
 \item[$(\mr{i})$] $L\times\{0\}\subset M'$, $\phi(x)=(x,0)$ and $f'(x,0)=f(x)$ for all $x\in L$.
 \item[$(\mr{ii})$] The smooth embedding $M\hookrightarrow\R^{n+t}$, $x\mapsto\phi(x)$ is arbitrarily $\cinfty$ close to the inclusion map $M\hookrightarrow\R^{n+t}$, $x\mapsto(x,0)$.
 \item[$(\mr{iii})$] $f'\circ\phi$ is arbitrarily $\cinfty$ close to $f$.
\end{itemize}
\end{thm}
\begin{proof}
Choose $k\in\N$ so that $n+k+1\geq 2(d+1)+1$ (thus $n+k\geq2d+1$). Consider $M$ (and thus $L$) as a subset of $\R^{n+k}$. By Lemmas \ref{lem:projective}$(\mr{iv})$ and \ref{lem:basic-Q-pair}$(\mr{iv})$, $L$ is also a projectively $\Q$-closed $\Q$-stable $\Q$-algebraic subset of $\R^{n+k}$. Let $\beta'_M:M\to\G_{n+k,n+k-d}$ be the normal bundle map of $M$ in $\R^{n+k}$. Observe that, for all $a\in M$, $\beta'_M(a)$ is the $(n+k)\times(n+k)$-matrix with the leading $n\times n$-submatrix equal to $\beta_M(a)$ and with the $i^\mr{th}$ row and the $i^\mr{th}$ column equal to the $i^\mr{th}$ vector of the canonical basis of $\R^{n+k}$ for every $i\in\{n+1,\ldots,n+k\}$, when $k\geq1$. It follows that $\beta'_M|_L$ is $\Q$-regular. To simplify notations, we rename $n+k$ as $n$, and $\beta'_M$ as $\beta_M$.

The unoriented bordism class of $f$ is projectively $\Q$-algebraic by hypothesis, so Theo\-rem \ref{thm:generic-projection} and Remark \ref{rem:244} assure that there exist a $d$-dimensional projectively $\Q$-closed $\Q$-nonsingular $\Q$-algebraic set $P\subset\R^n$ with $M\cap P=\varnothing$, a $\Q$-regular map $g:P\to W$, $d+1$-dimensional compact smooth submanifold $S$ of $\R^{n+1}=\R^n\times\R$ and a smooth map $G:S\to W$ such that $S\cap(\R^n\times(-1,1))=(M\sqcup P)\times(-1,1)$, $G(x,0)=f(x)$ for all $x\in M$ and $G(x,0)=g(x)$ for all $x\in P$.

By Lemmas \ref{lem:projective}$(\mr{i})(\mr{iv})$ and \ref{lem:basic-Q-pair}$(\mr{iii})(\mr{iv})$, $(L\times\{0\})\sqcup(P\times\{0\})=(L\sqcup P)\times\{0\}\subset\R^{n+1}$ is a non-empty projectively $\Q$-closed $\Q$-stable $\Q$-algebraic set.

Let $\G:=\G_{n+1,n-d}\subset\R^{(n+1)^2}$, let $\E:=\E_{n+1,n-d}\subset\R^{(n+1)^2+n+1}$ be the total space of the universal vector bundle over $\G$, let $U$ be an open tubular neighborhood of $S$ in $\R^{n+1}$, let $B:S\to\G$ be the normal bundle map of $S$ in $\R^{n+1}$, let $\rho\in\cinfty(U,S)$ be the closest point map and let $\theta:U \to\E$ be the smooth map $\theta(x,x_{n+1}):=(B(\rho(x,x_{n+1})),(x,x_{n+1})-\rho(x,x_{n+1}))$. Evidently, $\theta^{-1}(\G\times\{0\})=S$ and $\theta$ is transverse to $\G\times\{0\}$ in $\E$. Let $\beta_P:P\to\G_{n,n-d}$ be the normal bundle map of $P$ in $\R^n$, which is $\Q$-regular by Lemma \ref{lem:gauss}$(\mr{i})$. Observe that $\theta(x,0)=(B(x,0),0)$ for all $x\in L\sqcup P$. Since $S\cap(\R^n\times(-1,1))=(M\sqcup P)\times(-1,1)$, we have that $B(x,0)$ is the $(n+1)\times(n+1)$-matrix with the last row and the last column equal to zero and with the leading $n\times n$-submatrix equal to $\beta_M(x)$ if $x\in M$ and $\beta_P(x)$ if $x\in P$. Since $\beta_M|_L$ and $\beta_P$ are $\Q$-regular, Lemma \ref{lem:Q-basic}$(\mr{iv})$ assures that the restriction $\theta|_{(L\sqcup P)\times\{0\}}$ is $\Q$-regular.  

Let $F: U \to W$ be the smooth map defined by $F:=G\circ\rho$ and let $\theta\times F:U\to \E\times W$ be the product map. Since $F(x,0)=G(x,0)=f|_L(x)$ for all $x\in L$ and $F(x,0)=G(x,0)=g(x)$ for all $x\in P$, Lemma \ref{lem:Q-basic}$(\mr{iv})$ also assures that the restriction $F|_{(L\sqcup P)\times\{0\}}$ is $\Q$-regular. It follows that the restriction $(\theta\times F)|_{(L\sqcup P)\times\{0\}}=\theta|_{(L\sqcup P)\times\{0\}}\times F|_{(L\sqcup P)\times\{0\}}$ is $\Q$-regular as well. Observe that $\E\times W$ is a $\Q$-nonsingular $\Q$-algebraic subset of $\R^{(n+1)^2+n+1+k}$ by Lemma \ref{lem:Q-grassmannians}$(\mr{ii})$.

Let $V$ be an open neighborhood of $S$ in $\R^{n+1}$ such that $\overline{V}:=\cl_{\R^{n+1}}(V)$ is a compact smooth submanifold with boundary of $U$. Apply Lemma \ref{lem:283} to $\theta\times F:U\to\E\times W$. We obtain a smooth map $v:\overline{V}\to\R^s$ with $s:=(n+1)^2+n+1+k$, a $\Q$-algebraic subset $Z$ of $\R^{n+1+s}=\R^n\times\R\times\R^s$ of dimension $n+1$ and a $\Q$-regular map $\eta:Z\to\E\times W$ such that:
\begin{itemize}
 \item $\widehat{V}:=\{(x,x_{n+1},v(x,x_{n+1}))\in\R^{n+1+s}:(x,x_{n+1})\in V\}$ is an open subset of $\Reg^\Q(Z)$.
 \item $(L\sqcup P)\times\{0\}\times\{0\}\subset\widehat{V}$, and $v(x,0)=0$ and $\eta(x,0,0)=(\theta\times F)(x,0)$ for all $x\in L\sqcup P$.
 \item $v$ can be chosen arbitrarily $\cinfty$ close to the zero map $\overline{V}\to\R^s$, $(x,x_{n+1})\mapsto0$. Moreover, if $\widehat{v}:\overline{V}\to Z$ denotes the smooth map $\widehat{v}(x,x_{n+1}):=(x,x_{n+1},v(x,x_{n+1}))$, then $\eta$ can be chosen such that $\eta\circ\widehat{v}$ is arbitrarily $\cinfty$ close to $(\theta\times F)|_{\overline{V}}$.
\end{itemize}
In addition, using again Lemma \ref{lem:projective}$(\mr{iv})$, we deduce that
\begin{itemize}
 \item $(L\sqcup P)\times\{0\}\times\{0\}\subset\R^{n+1+s}$ is a projectively $\Q$-closed $\Q$-algebraic set.
\end{itemize} 

Since $(\theta\times F)|_{\overline{V}}$ is transverse to $\G\times\{0\}\times W$ in $\E\times W$, $((\theta\times F)|_{\overline{V}})^{-1}(\G\times\{0\}\times W))=\overline{V}\cap\theta^{-1}(\G\times\{0\})=S$ and $S\cap\partial\overline{V}=\varnothing$, by \cite[Thm.14.1.1]{BCR}, if $\eta\circ\widehat{v}$ is sufficiently $\cinfty$ close to $(\theta\times F)|_{\overline{V}}$, then $\eta\circ\widehat{v}$ is transverse to $\G\times\{0\}\times W$ in $\E\times W$, $S_1:=(\eta\circ\widehat{v})^{-1}(\G\times\{0\}\times W)$ is a compact smooth submanifold of $\overline{V}$ with $S_1\cap\partial\overline{V}=\varnothing$, and there exists a smooth diffeomorphism $\xi:\overline{V}\to\overline{V}$ arbitrarily $\cinfty$ close to $\mr{id}_{\overline{V}}$ such that $\xi(S)=S_1$ and $\xi$ fixes $\partial\overline{V}$, i.e., $\xi(p)=p$ for all $p\in\partial\overline{V}$. In addition, since $(\eta\circ\widehat{v})(x,0)=(\theta\times F)(x,0)$ for all $x\in L\sqcup P\subset S$, we have that $(L\sqcup P)\times\{0\}\subset S_1$ and we can also assume that $\xi$ fixes $(L\sqcup P)\times\{0\}$.

Since $\widehat{v}|_V:V\to Z$ is a smooth embedding, it follows that $\eta|_{\widehat{V}}$ is transverse to $\G\times\{0\}\times W$ in $\E\times W$. Let $Z_1'\subset\R^{n+1+s}$ be the $\Q$-algebraic set $Z_1':=\eta^{-1}(\G\times\{0\}\times W)$ and let $\widehat{S}_1:=\widehat{V}\cap Z_1'$. By Proposition \ref{prop:Q-transverse}, we have $\widehat{S}_1\subset\Reg^{\Q}(Z_1',d+1)$, so $\widehat{S}_1$ is an open subset of $\Reg^{\Q}(Z_1',d+1)$. Denote $Z_1\subset\R^{n+1+s}$ the union of all $\Q$-irreducible components of $Z_1'$ of dimension $d+1$, so $\dim(Z_1)=d+1$. By Proposition \ref{prop211}$(\mr{i})$, we deduce that $\widehat{S}_1=\widehat{V}\cap Z_1$ is an open subset of $\Reg^\Q(Z_1)$ (and so of $Z_1$). Since $\widehat{S}_1=\widehat{v}(S_1)=\widehat{v}(\xi(S))$, $\widehat{S}_1$ is also compact in $Z_1$ (and so in $\R^{n+1+s}$). It follows that $\widehat{S}_1$ is a union of some connected components of $Z_1$. In addition, the map $\widehat{\xi}:S\to\R^{n+1+s}$, sending $(x,x_{n+1})$ to $\widehat{v}(\xi(x,x_{n+1}))$, is a well-defined smooth embedding arbitrarily $\cinfty$ close to the inclusion map $J:S\hookrightarrow\R^{n+1+s}$, $(x,x_{n+1})\mapsto(x,x_{n+1},0)$ such that $\widehat{\xi}(S)=\widehat{S}_1$ and $\widehat{\xi}(x,0)=(x,0,0)$ for all $x\in L\sqcup P$.

Let $(x,x_{n+1},y)$ be the coordinates of $\R^{n+1+s}=\R^n\times\R\times\R^s$, let $\pi_{n+1}:\R^{n+1+s}\to\R$ be the projection $\pi_{n+1}(x,x_{n+1},y):=x_{n+1}$ and let $H$ be the hyperplane $(\pi_{n+1})^{-1}(0)$ of $\R^{n+1+s}$. Since $\widehat{\xi}$ is arbitrarily $\cinfty$ close to $J$ and $(L\sqcup P)\times\{0\}\times\{0\}\subset\widehat{S}_1\cap H$, we have that $\widehat{S}_1$ is transverse to $H$ in $\R^{n+1+s}$, $\widehat{S}_1\cap H=M_1\sqcup(P\times\{0\}\times\{0\})$ for some compact smooth submanifold $M_1$ of $\R^{n+1+s}$ containing $L\times\{0\}\times\{0\}$ and there exists a smooth diffeomorphism $\phi_1:M\to M_1$ such that $\phi_1(x)=(x,0,0)$ for all $x\in L$ and the smooth embedding $M\hookrightarrow\R^{n+1+s}$, $x\mapsto\phi_1(x)$ is arbitrarily $\cinfty$ close to the inclusion map $M\hookrightarrow\R^{n+1+s}$, $x\mapsto(x,0,0)$.

Since $\widehat{S}_1$ is compact in $\R^{n+1+s}$ and $Z_1\setminus\widehat{S}_1$ is closed in $Z_1$ (and so in $\R^{n+1+s}$), there exists a compact neighborhood $K$ of $\widehat{S}_1$ in $\R^{n+1+s}$ such that $K\cap(Z_1\setminus\widehat{S}_1)=\varnothing$. Let $\y=(\y_1,\ldots,\y_s)$ and $(\x,\x_{n+1},\y)=(\x_1,\ldots,\y_s)$ be the indeterminates corresponding to the coordinates of $\R^s$ and $\R^{n+1+s}$, respectively. Since $(L\sqcup P)\times\{0\}\times\{0\}$ is a non-empty projectively $\Q$-closed $\Q$-algebraic subset of $\R^{n+1+s}$, by Lemma \ref{lem:overt}, there exists a {non-constant} overt polynomial $q\in\Q[\x,\x_{n+1},\y]$ such that $\ZZ_\R(q)=(L\sqcup P)\times\{0\}\times\{0\}$. Since the function $q:\R^{n+1+s}\to\R$ is proper, replacing $q$ with $cq^2$ for some rational number $c>0$ if necessary, we can assume that: $q\in\Q[\x,\x_{n+1},\y]$ is overt, $\ZZ_\R(q)=(L\sqcup P)\times\{0\}\times\{0\}$, $q\geq0$ on $\R^{n+1+s}$ and $q\geq2$ on $\R^{n+1+s}\setminus K$. Let $K_1$ be a compact neighborhood of $\widehat{S}_1$ in $\R^{n+1+s}$ such that $K_1$ is a compact smooth submanifold with boundary of $\mr{int}_{\R^{n+1+s}}(K)$. Using a smooth partition of unity subordinate to the open cover $\{\mr{int}_{\R^{n+1+s}}(K),\R^{n+1+s}\setminus K_1\}$ of $\R^{n+1+s}$, we can define a smooth function $h:\R^{n+1+s}\to\R$ such that $h=\pi_{n+1}$ on $K_1$ and $h=q$ on $\R^{n+1+s}\setminus K$. Apply Lemma \ref{lem:L} to $h-q$ obtaining a $\Q$-regular function $u_1:\R^{n+1+s}\to\R$, a natural number $e$ and a polynomial $p\in\Q[\x,\x_{n+1},\y]$ of degree $\leq 2e$ such that $u_1(z)=p(z)(1+|z|_{n+1+s}^2)^{-e}$ for all $z\in\R^{n+1+s}$, $(L\sqcup P)\times\{0\}\times\{0\}\subset\ZZ_\R(p)$, $\sup_{x\in\R^{n+1+s}}|h(x)-q(x)-u_1(x)|<1$ and $u_1|_{K_1}$ is arbitrarily $\cinfty$ close to $(h-q)|_{K_1}=(\pi_{n+1}-q)|_{K_1}$. Let $u:\R^{n+1+s}\to\R$ be the $\Q$-regular map $u:=u_1+q$ and let $v\in\Q[\x,\x_{n+1},\y]$ be the unique polynomial such that $v(z)=p(z)+q(z)(1+|z|_{n+1+s}^2)^e$ and so $u(z)=v(z)(1+|z|_{n+1+s}^2)^{-e}$ for all $z\in\R^{n+1+s}$. Since $q$ is a {non-constant} overt polynomial and $\deg(p)\leq2e$, we immediately deduce that the polynomial $v$ is overt. In addition, we have that $u(x,0,0)=0$ for all $x\in L\sqcup P$, $u>1$ on $\R^{n+1+s}\setminus K$ and $u|_{K_1}$ is arbitrarily $\cinfty$ close to $\pi_{n+1}|_{K_1}$.

Since $u|_{K_1}$ can be chosen arbitrarily $\cinfty$ close to $\pi_{n+1}|_{K_1}$, we can also assume that $u|_{\widehat{S}_1}$ is arbitrarily $\cinfty$ close to $\pi_{n+1}|_{\widehat{S}_1}$. Observe that $0$ is a regular value of $\pi_{n+1}|_{\widehat{S}_1}$ and $(\pi_{n+1}|_{\widehat{S}_1})^{-1}(0)=M_1\sqcup(P\times\{0\}\times\{0\})$. By \cite[Thm.14.1.1]{BCR}, if $u|_{\widehat{S}_1}$ is sufficiently $\cinfty$ close to $\pi_{n+1}|_{\widehat{S}_1}$, then $0$ is a regular value of $u|_{\widehat{S}_1}$, $Q:=(u|_{\widehat{S}_1})^{-1}(0)$ is a compact smooth submanifold of $\widehat{S}_1$, and there exists a smooth diffeomorphism $\widehat{\phi}_1:\widehat{S}_1\to\widehat{S}_1$ arbitrarily $\cinfty$ close to $\mr{id}_{\widehat{S}_1}$ such that $\widehat{\phi}_1(M_1\sqcup(P\times\{0\}\times\{0\}))=Q$. Since $u|_{\widehat{S}_1}(x,0,0)=0=\pi_{n+1}|_{\widehat{S}_1}(x,0,0)$ for all $x\in L\sqcup P$, we have $(L\sqcup P)\times\{0\}\times\{0\}\subset Q$ and we can assume that $\widehat{\phi}_1$ fixes $(L\sqcup P)\times\{0\}\times\{0\}$. Therefore, if we define $M':=\widehat{\phi}_1(M_1)$, then $L\times\{0\}\times\{0\}\subset M'$ and $Q=M'\sqcup(P\times\{0\}\times\{0\})$. In addition, if we denote $\phi:M\to M'$ the smooth diffeomorphism $\phi(x):=\widehat{\phi}_1(\phi_1(x))$, then we can assume that the smooth embedding $M\hookrightarrow\R^{n+1+s}$, $x\mapsto\phi(x)$ is arbitrarily $\cinfty$ close to the inclusion map $M\hookrightarrow\R^{n+1+s}$, $x\mapsto(x,0,0)$. Observe that $\phi(x)=(x,0,0)$ for all $x\in L$.

Consider the restriction $u|_{Z_1}:Z_1\to\R$, which is a $\Q$-regular map. Since $u>1$ on $\R^{n+1+s}\setminus K\supset Z_1\setminus\widehat{S}_1$, we have $Q=(u|_{Z_1})^{-1}(0)$. In addition, $u|_{Z_1}$ coincides with $u|_{\widehat{S}_1}$ on the open subset $\widehat{S}_1$ of $\Reg^\Q(Z_1)$ and $0$ is a regular value of $u|_{\widehat{S}_1}$, so we deduce that $Q\subset\Reg^\Q((u|_{Z_1})^{-1}(0))=\Reg^\Q(Q)$ by Proposition \ref{prop:Q-transverse}. This proves that $Q=M'\sqcup(P\times\{0\}\times\{0\})$ is a $\Q$-nonsingular $\Q$-algebraic subset of $\R^{n+1+s}$ of dimension $d$. Since $P\times\{0\}\times\{0\}$ is also a $\Q$-nonsingular $\Q$-algebraic subset of $\R^{n+1+s}$ of dimension $d$, Proposition \ref{prop:Q_setminus} assures that $M'$ is a $\Q$-nonsingular $\Q$-algebraic subset of $\R^{n+1+s}$ of dimension $d$ as well.

Since $u(z)=v(z)(1+|z|_{n+1+s}^2)^{-e}$ for all $z\in\R^{n+1+s}$ and $v\in\Q[\x,\x_{n+1},\y]$ is overt, we have $M'\subset Q=(u|_{Z_1})^{-1}(0)\subset\ZZ_\R(v)$ and so the $\Q$-algebraic set $M'\subset\R^{n+1+s}$ is projectively $\Q$-closed by Lemma \ref{lem:projective}$(\mr{ii})$. 

Let $\Pi:\E\times W\to W$ be the projection onto the second factor and let $f':M'\to W$ be the $\Q$-regular map $f'(z):=\Pi(\eta(z))$. Since $\eta$ can be chosen such that $\eta\circ\widehat{v}$ is arbitrarily $\cinfty$ close to $(\theta\times F)|_{\overline{V}}$, $M\times\{0\}\subset\overline{V}$ and the smooth embedding $M\hookrightarrow\R^{n+1+s}$, $x\mapsto\phi(x)$ can be chosen arbitrarily $\cinfty$ close to the smooth embedding $M\hookrightarrow\R^{n+1+s}$, $x\mapsto\widehat{v}(x,0)$, then we can assume that the smooth map $M\times\{0\}\to W$, $(x,0)\mapsto\Pi(\eta(\phi(x)))=f'(\phi(x))$ is arbitrarily $\cinfty$ close to the smooth map $M\times\{0\}\to W$, $(x,0)\mapsto F(x,0)=f(x)$. In other words, $f'\circ\phi$ is arbitrarily $\cinfty$ close to $f$. Observe that $f'(x)=\Pi(\eta(x,0,0))=F(x,0)=f|_L(x)=f(x)$ for all $x\in L$.

To complete the proof, it is now sufficient to set $t:=1+s$.
\end{proof}

\begin{remark}\label{rem:varnothing}
If $L=\varnothing$, the statement of Theorem \ref{thm:Q_tognoli} remains valid: it suffices to cancel the conditions `\emph{the restrictions $\beta_M|_L:L\to\G_{n,n-d}$ and $f|_L:L\to W$ are $\Q$-regular}' and item $(\mr{i})$. In this case the proof is similar but simpler. $\sqbullet$
\end{remark}

We are ready to prove the Nash-Tognoli theorem over the rationals.

\begin{proof}[Proof of Theorem \ref{thm:NTQ}]
Let $M$ be a compact smooth manifold of dimension $d$ and let $\psi:M\to\R^{2d+1}$ be a smooth embedding. Set $n:=2d+1$ for short. Identifying $M$ with $\psi(M)$, we can assume that $M$ is a compact smooth submanifold of $\R^n$ and $\psi$ is the inclusion map $M\hookrightarrow\R^n$. 

Let $W$ be the singleton $\{0\}$ of $\R$ and let $f:M\to W$ be the constant map. Apply Theorem~\ref{thm:Q_tognoli} to $f$ with $L=\varnothing$ (see Remark \ref{rem:varnothing}). We obtain a projectively $\Q$-closed $\Q$-nonsingular $\Q$-algebraic set $M_1\subset\R^{n+t}=\R^n\times\R^t$ for some $t\in\N^*$ and a smooth diffeomorphism $\phi_1:M\to M_1$ such that the smooth embedding $M\hookrightarrow\R^{n+t}$, $x\mapsto\phi_1(x)$ is arbitrarily $\cinfty$ close to the inclusion map $M\hookrightarrow\R^{n+t}$, $x\mapsto(x,0)$.

Let $y=(y_1,\ldots,y_t)$ be the coordinates of $\R^t$ and let $(x,y)$ be the coordinates of $\R^{n+t}=\R^n\times\R^t$. Endow the set $\mc{M}_{n,t}(\Q)$ of rational $n\times t$-matrices with the topology induced by the Euclidean one of $\mc{M}_{n,t}(\R)=\R^{nt}$. Combining Theorem \ref{thm:generic-projection} and the density of $\Q$ in $\R$, there exists a matrix $A\in\mc{M}_{n,t}(\Q)$ arbitrarily close to the zero matrix $O$ such that the corresponding projection $\pi_A:\R^{n+t}\to\R^n$, $(x,y)\mapsto x-Ay$ (here $x$ and $y$ are interpreted as column vectors) has the following properties: $M':=\pi_A(M_1)\subset\R^n$ is a projectively $\Q$-closed $\Q$-nonsingular $\Q$-algebraic set and the restriction $\pi'_A:M_1\to M'$ of $\pi_A$ is a $\Q$-biregular isomorphism. Let $\phi:M\to\R^n$ be the smooth embedding $\phi(x):=\pi_A(\phi_1(x))$. If we choose $A$ sufficiently close to $O$ in $\mc{M}_{n,t}(\Q)$, then we can assume that $\phi$ is arbitrarily $\cinfty$ close to the inclusion map $\psi:M\hookrightarrow\R^n$. This completes the proof.
\end{proof}

\begin{remark}\label{rem112}
If we replace $n=2d+1$ with a larger natural number $n$, the previous proof continues to work {as} Theorem \ref{thm:generic-projection} remains valid, see Remark \ref{rem111}. $\sqbullet$
\end{remark}

Another consequence of Theorems \ref{thm:generic-projection} and \ref{thm:Q_tognoli} is the following $\Q$-version of \cite[Thm.2.10]{ak1981}, which is a relative version of Theorem \ref{thm:NTQ}.

\begin{cor}\label{thm:Q_tico_approx}
Let $M$ be a compact smooth submanifold of $\R^{2d+1}$ of dimension $d$, let $\ell\in\N^*$ and let $\{M_i\}_{i=1}^\ell$ be a family of smooth hypersurfaces of $M$ in general position. Then there exist a projectively $\Q$-closed $\Q$-nonsingular $\Q$-algebraic set $M'\subset\R^{2d+1}$, a family $\{M'_i\}_{i=1}^\ell$ of $\Q$-nonsingular $\Q$-algebraic hypersurfaces of $M'$ and a smooth diffeomorphism $\psi:M\to M'$ such that $\psi(M_i)=M'_i$ for every $i\in\{1,\ldots,\ell\}$, and the smooth embedding $M\hookrightarrow\R^{2d+1}$, $x\mapsto\psi(x)$ is arbitrarily $\cinfty$ close to the inclusion map $M\hookrightarrow\R^{2d+1}$.
\end{cor}

In order to prove this corollary, we adapt the strategy of the proof of \cite[Thm.2.10]{ak1981} to the present $\Q$-algebraic setting. First, we need to recall a known result on `small' automorphisms, the proof of which will be outlined here for the sake of completeness.

\begin{lem}\label{lem:small-embeddings}
Let $M$ be a compact smooth manifold, let $\ell\in\N^*$ and let $\{M_i\}_{i=1}^\ell$ be a family of smooth hypersurfaces of $M$ in general position. Then, for every neighborhood $\mc{U}$ of $\mr{id}_M$ in $\cinfty(M,M)$, there exists a neighborhood $\mc{V}$ of $\mr{id}_M$ in $\cinfty(M,M)$ with the following property: if $\{\xi_i:M\to M\}_{i=1}^\ell$ is a family of smooth diffeomorphisms belonging to $\mc{V}$, then there exists a smooth diffeomorphism $\xi:M\to M$ belonging to $\mc{U}$ such that $\xi(M_i)=\xi_i(M_i)$ for every $i\in\{1,\ldots,\ell\}$.
\end{lem}
\begin{proof}
By Whitney's embedding theorem, we can assume that $M$ is a compact smooth submanifold of some $\R^n$. Let $\rho:U\to M$ be an open tubular neighborhood of $M$ in $\R^n$. If the neighborhood $\mc{V}$ of $\mr{id}_M$ in $\cinfty(M,M)$ is sufficiently small, then the point $(1-t)x+t\xi_i(x)$ belongs to $U$ for every $i\in\{1\ldots,\ell\}$, $t\in[0,1]$ and $x\in M_i$, so we can define the smooth map $\xi^*_{it}:M_i\to M$ by $\xi^*_{it}(x):=\rho((1-t)x+t\xi_i(x))$. Since the set of smooth embeddings from $M_i$ to $M$ is open in $\cinfty(M_i,M)$ (see \cite[Thm.1.4]{hirsch:difftop}), shrinking $\mc{V}$ around $\mr{id}_M$ if necessary, we can assume that $\{\xi^*_{it}\}_{t\in[0,1]}$ is an isotopy from the inclusion $M_i\hookrightarrow M$ to $\xi^*_{i1}=\xi_i|_{M_i}$. Now we can apply  \cite[Lem.2.9]{ak1981} to the isotopies $\{\xi^*_{it}\}_{t\in[0,1]}$ for $i=1,\ldots,\ell$, obtaining an isotopy $\{\xi_t:M\to M\}_{t\in[0,1]}$ such that $\xi_t(M_i)=\xi^*_{it}(M_i)$ and $\xi_t\in\mc{U}$, provided $\mc{V}$ is sufficiently small, for all $i\in\{1,\ldots,\ell\}$ and $t\in[0,1]$. Finally, it is sufficient to set $\xi:=\xi_1$ and observe that $\xi(M_i)=\xi^*_{i1}(M_i)=\xi_i(M_i)$.   
\end{proof}

\begin{proof}[Proof of Corollary \ref{thm:Q_tico_approx}]
Let us start by recalling a well-known construction in {the theory of vector bundles}. Consider the ring $\cinfty(M)$ of real-valued smooth functions on $M$. {Since $M$ is compact,} let $\{f_{i1},\ldots,f_{im_i}\}$ be a finite system of generators of the vanishing ideal $I_i$ of $M_i$ in $\cinfty(M)$ {for every $i\in\{1,\ldots,\ell\}$}. Setting $m:=\max_{i\in\{1,\ldots,\ell\}}m_i$ and adding generators if necessary, we can assume that $m=m_1=\ldots=m_\ell$. Since $I_i$ is locally principal, there exists a smooth map $\sigma_i:M\to\G_{m,1}$ such that, for every $x\in M\setminus M_i$, $\sigma_i(x)$ is the matrix associated to the orthogonal projection of $\R^m$ onto the line generated by the vector $(f_{i1}(x),\ldots,f_{im}(x))$ with respect to the canonical basis of $\R^m$. Consider the universal vector bundle $\gamma_{m,1}:=(\E_{m,1},p_{m,1},\G_{m,1})$ over $\G_{m,1}$ and define the smooth map $f_i:M\to\E_{m,1}$ by $f_i(x):=(\sigma_i(x),(f_{i1}(x),\ldots,f_{im}(x)))$. Observe that $f_i$ is transverse in $\E_{m,1}$ to its zero section $\G_{m,1}\times\{0\}$. Since the Thom space of $\gamma_{m,1}$ is $\G_{m+1,1}$, we can assume that $f_i$ is a smooth map from $M$ to $\G_{m+1,1}$ transverse to $\G_{m,1}$ and such that $f_i^{-1}(\G_{m,1})=M_i$.

Define the subsets $W,W_1,\ldots,W_\ell$ of $\R^{\ell(m+1)^2}=(\R^{(m+1)^2})^\ell$ and the smooth map $f:M\to W$ as follows: $W:=(\G_{m+1,1})^\ell$, $W_i:=(\G_{m+1,1})^{i-1}\times\G_{m,1}\times(\G_{m+1,1})^{\ell-i}$ for $i=1,\ldots,\ell$, and $f:=f_1\times\ldots\times f_\ell$. Observe that $W\subset\R^{\ell(m+1)^2}$ is a $\Q$-nonsingular $\Q$-algebraic set, each $W_i$ is a $\Q$-nonsingular $\Q$-algebraic hypersurface of $W$, $f$ is transverse to each $W_i$ in $W$ and $f^{-1}(W_i)=M_i$ for all $i\in\{1,\ldots,\ell\}$. In addition, $W$ has projectively $\Q$-algebraic unoriented bordism by Lemma~\ref{lem:kunneth}. Set $n:=2d+1$ for short. Apply Theorem \ref{thm:Q_tognoli} to $f$ with $L=\varnothing$ (see Remark~\ref{rem:varnothing}). We obtain a projectively $\Q$-closed $\Q$-nonsingular $\Q$-algebraic set $M'\subset\R^{n+t}=\R^n\times\R^t$ for some $t\in\N^*$, a smooth diffeomorphism $\phi:M\to M'$ and a $\Q$-regular map $f':M'\to W$ such that the smooth embedding $M\hookrightarrow\R^{n+t}$, $x\mapsto\phi(x)$ is arbitrarily $\cinfty$ close to the inclusion map $M\hookrightarrow\R^{n+t}$, $x\mapsto(x,0)$, and $f'\circ\phi$ is arbitrarily $\cinfty$ close to $f$. By \cite[Thm.14.1.1]{BCR}, if $f'\circ\phi$ is sufficiently $\cinfty$ close to $f$, then $f'\circ\phi$ is transverse to each $W_i$ in $W$, $M_{i1}:=(f'\circ\phi)^{-1}(W_i)$ is a compact smooth hypersurface of $M$ and there exists a smooth diffeomorphism $\xi_i:M\to M$ arbitrarily $\cinfty$ close to $\mr{id}_M$ such that $\xi_i(M_i)=M_{i1}$. By Lemma \ref{lem:small-embeddings}, if each $\xi_i$ is sufficiently $\cinfty$ close to $\mr{id}_M$, then there exists a smooth diffeomorphism $\xi:M\to M$ arbitrarily $\cinfty$ close to $\mr{id}_M$ such that $\xi(M_i)=M_{i1}$ for all $i\in\{1,\ldots,\ell\}$. Observe that each $M'_i:=\phi(M_{i1})$ is equal to $(f')^{-1}(W_i)$. Since $\phi$ is a smooth diffeomorphism, $f'$ is also transverse to $W_i$ in $W$. By Proposition \ref{prop:Q-transverse}, $M'_i$ is a $\Q$-nonsingular $\Q$-algebraic hypersurface of $M'$. Define the smooth diffeomorphism $\psi:M\to M'$ by $\psi:=\phi\circ\xi$. Observe that $M'_i=\psi(M_i)$ for all $i\in\{1,\ldots,\ell\}$ and the smooth embedding $\psi':M\hookrightarrow\R^{n+t}$, $x\mapsto\psi(x)$ is arbitrarily $\cinfty$ close to the inclusion map $M\hookrightarrow\R^{n+t}$, $x\mapsto(x,0)$. Composing $\psi'$ with a suitable linear projection $\pi_A:\R^{n+t}\to\R^n$ as in the last part of the proof of Theorem \ref{thm:NTQ}, we can replace the ambient space $\R^{n+t}$ with $\R^n$.
\end{proof}


\subsection{Real $\Q$-algebraic blowing down with approximation}
Here we give a $\Q$-version of the blowing down lemma \cite[Lem.2.6.1]{akbking:tras} integrated with an approximation part.

Let $V$ and $W$ be two topological spaces, let $E$ be a non-empty topological subspace of $V$ and let $f:E\to W$ be a continuous map. Consider the disjoint union $V\sqcup W$ of $V$ and $W$ endowed with the natural topology, and the equivalence relation $\sim$ on $V\sqcup W$ that identifies each point $y$ of $E$ with $f(y)\in W$ or, equivalently, that has the following equivalence classes: $[y]_\sim=\{y\}$ if $y\in V\setminus E$ and $[x]_\sim=f^{-1}(x)\cup\{x\}$ if $x\in W$. The corresponding quotient topological space $V\cup_fW$ of $V\sqcup W$ modulo $\sim$ is called adjunction space of the system $V\supset E\stackrel{_f}{\to}W$. 

\begin{lem}\label{lem:blowing_down}
Let $V\subset\R^m$ be a projectively $\Q$-closed $\Q$-nonsingular $\Q$-algebraic set, let $W\subset\R^n$ be a $\Q$-algebraic set, let $F:V\to\R^n$ be a $\Q$-regular map, let $E:=F^{-1}(W)$ be such that $E\neq\varnothing$, let $f:E\to W$ be the restriction of $F$ from $E$ to $W$ and let $\Pi:V\sqcup W\to V\cup_f W$ be the quotient map from $V\sqcup W$ to the adjunction space $V\cup_f W$ of the system $V\supset E\stackrel{_f}{\to}W$. Then there exist a $\Q$-algebraic set $T\subset\R^{n+m}$, $\Q$-regular maps $\xi:V\to T$ and $\eta:W\to T$, and a homeomorphism $h :V\cup_f W\to T$ with the following properties:
\begin{itemize}
 \item[$(\mr{i})$] $h\circ\Pi|_V=\xi$ and $h\circ\Pi|_W=\eta$, i.e., the following diagram commutes:
\[
\centering
\begin{tikzcd}
V\sqcup W \arrow[d,"\Pi"] \arrow[r,"\xi\sqcup\eta"] & T \\
V\cup_f W \arrow[ur, bend right=20, "h"] & 
\end{tikzcd}
\]
 \item[$(\mr{ii})$] $\eta(x)=(x,0)\in\R^n\times\R^m=\R^{n+m}$ for all $x\in W$.
 \item[$(\mr{iii})$] $\xi(V)\cup(W\times\{0\})=T$, $\xi(V\setminus E)=T\setminus(W\times\{0\})$ and the restriction of $\xi$ from $V\setminus E$ onto $T\setminus(W\times\{0\})$ is a $\Q$-biregular isomorphism.
 \item[$(\mr{iv})$] The map $V\to\R^{n+m}$, $y\mapsto\xi(y)$ is arbitrarily $\cinfty$ close to $V\to\R^{n+m}$, $y\mapsto(F(y),0)$.
 \item[$(\mr{v})$] If $\max\{\dim(E),\dim(W)\}\leq\dim(V\setminus E)$, then $T\setminus(W\times\{0\})\subset\Reg^\Q(T)$.
\end{itemize}
\end{lem}
\begin{proof}
Let $x=(x_1,\ldots,x_n)$ be the coordinates of $\R^n$, let $y=(y_1,\ldots,y_m)$ be the coordinates of $\R^m$ and let $\x=(\x_1,\ldots,\x_n)$ and $\y=(\y_1,\ldots,\y_m)$ be the corresponding indeterminates. Let $w\in\Q[\x]$ and $v\in\Q[\y]$ be such that $\ZZ_\R(w)=W$ and $\ZZ_\R(v)=V$. Replacing $w$ with $cw$ for some sufficiently small constant $c\in\Q\setminus\{0\}$, we can assume that the corresponding polynomial function $w:\R^n\to\R$ is arbitrarily $\cinfty$ close to the zero function on $\R^n$. Since $V\subset\R^m$ is non-empty and projectively $\Q$-closed, the degree $d$ of $v$ is positive and, by Lemma \ref{lem:overt}, we can assume that $v$ is overt. Denote $v_d$ the leading homogeneous term of $v$. Let $p_1,\ldots,p_n,q\in\Q[\y]$ be such that $q(y)\neq0$ and $F(y)=\frac{p(y)}{q(y)}$ for all $y\in V$, where $p:\R^m\to\R^n$ is the polynomial map $p(y):=(p_1(y),\ldots,p_n(y))$. Choose $k\in\N^*$ such that $kd>\max\{\deg(q),\deg(p_1q),\ldots,\deg(p_nq)\}$, where we are assuming that the degree of the zero polynomial is equal to $0$. Observe that: $F(y)=\frac{p(y)q(y)}{q(y)^2+v(y)^{2k}}$ for all $y\in V$, $q^2+v^{2k}$ is an overt polynomial in $\Q[\y]$, whose degree is $2kd>2\max\{\deg(p_1q),\ldots,\deg(p_nq)\}$, whose leading homogeneous term is $v_d^{2k}$ and whose zero set is empty. Thus, replacing $p$ with $pq$ and $q$ with $q^2+v^{2k}$ if necessary, we can assume that $F(y)=\frac{p(y)}{q(y)}$ for all $y\in V$ with $p=(p_1,\ldots,p_n)$, where $p_1,\ldots,p_n,q$ are polynomials in $\Q[\y]$ such that $q$ is overt, $\deg(q)=2kd>2\max\{\deg(p_1),\ldots,\deg(p_n)\}$, the leading homogeneous term of $q$ is $v_d^{2k}$ and $\ZZ_\R(q)=\varnothing$.

Let $\Gamma=\{(F(y),y)\in\R^{n+m}:y\in V\}$ be the graph of $F$ and let $s$ be the unique polynomial in $\Q[\x,\y]$ such that $s(x,y)=|q(y)x-p(y)|_n^2+q(y)^2v(y)^2$ for all $(x,y)\in\R^{n+m}$, where $|\cdot|_n$ denotes the Euclidean norm of $\R^n$. Since $q$ never vanishes, a point $(x,y)\in\R^{n+m}$ belongs to $\ZZ_\R(s)$ if and only if $y\in\ZZ_\R(v)=V$ and $x=\frac{p(y)}{q(y)}=F(y)$, i.e., if and only if $(x,y)\in\Gamma$. This proves that $\ZZ_\R(s)=\Gamma$. In particular, for every $x\in W$, the set $\{y\in\R^m:s(x,y)=0\}$ coincides with $f^{-1}(x)$.

Given any polynomial $P\in\Q[\x,\y]=\Q[\x][\y]$, we denote $\deg_\y(P)$ the degree of $P$ with respect to $\y$. Observe that $s(\x,\y)=q(\y)^2(\sum_{i=1}^n\x_i^2)+(\sum_{i=1}^np_i(\y)^2)-2(\sum_{i=1}^n\x_iq(\y)p_i(\y))+q(\y)^2v(\y)^2$, where $\deg_\y(q(\y)^2(\sum_{i=1}^n\x_i^2))=4kd$, $\deg_\y(\sum_{i=1}^np_i(\y)^2)<2kd$, $\deg_\y(\sum_{i=1}^n\x_iq(\y)p_i(\y))\leq3kd$ and $\deg_\y(q(\y)^2v(\y)^2)=4kd+2d$. In addition, the leading homogeneous term of $q(\y)^2v(\y)^2$ is $v_d(\y)^{4k+2}$. It follows that there exists a polynomial $Q\in\Q[\x,\y]$ such that $s(\x,\y)=v_d(\y)^{4k+2}+Q(\x,\y)$ and $Q(\x,\y)=\sum_{\alpha\in D}Q_\alpha(\x)\y^\alpha$, where $D:=\{\alpha\in\N^m:|\alpha|\leq 4kd+2d-1\}$ and the $Q_\alpha$ are polynomials in $\Q[\x]$.

Consider the $\Q$-biregular isomorphism $\psi:(\R^n\setminus W)\times\R^m\to(\R^n\setminus W)\times\R^m$ defined by $\psi(x,y):=\big(x,\frac{y}{w(x)}\big)$, whose inverse is $\psi^{-1}(x,y)=(x,w(x)y)$. For every $(x,y)\in(\R^n\setminus W)\times\R^m$, it holds:
\begin{equation*}\label{equa}
\textstyle
s(\psi(x,y))=w(x)^{-4kd-2d}\big(v_d(y)^{4k+2}+w(x)\sum_{\alpha\in D}Q_\alpha(x)w(x)^{4kd+2d-1-|\alpha|}y^\alpha\big).
\end{equation*}
Define the polynomial $t\in\Q[\x,\y]$ by
\begin{equation*}\label{defa}
\textstyle
t(\x,\y):=v_d(\y)^{4k+2}+w(\x)\sum_{\alpha\in D}Q_\alpha(\x)w(\x)^{4kd+2d-1-|\alpha|}\y^\alpha
\end{equation*}
and the $\Q$-algebraic set $T:=\ZZ_\R(t)\subset\R^{n+m}$. Since $\ZZ_\R(s)=\Gamma$ and $s(\psi(x,y))=w(x)^{-4kd-2d}t(x,y)$ for all $(x,y)\in(\R^n\setminus W)\times\R^m$, we have that $\psi^{-1}(\Gamma\cap((\R^n\setminus W)\times\R^m))=T\cap((\R^n\setminus W)\times\R^m)$. Since $v$ is an overt polynomial of positive degree and  $t(x,y)=v_d(y)^{4k+2}$ for all $(x,y)\in W\times\R^m$, we also have that $\ZZ_\R(v_d)=\{0\}$ and so $T\cap(W\times\R^m)=W\times\{0\}$. It follows that the maps $\xi:V\to T$, $y\mapsto(F(y),w(F(y))y)$ and $\eta:W\to T$, $x\mapsto(x,0)$ are well-defined $\Q$-regular maps. Since the function $w:\R^n\to\R$ is arbitrarily $\cinfty$ close to the zero function on $\R^n$, we can assume that $\xi$ is arbitrarily $\cinfty$ close to the map $V\to\R^{n+m}$, $y\mapsto (F(y),0)$. Observe that $\eta(W)=W\times\{0\}$ and $\xi^{-1}(W\times\{0\})=F^{-1}(W)=E$. Moreover, the restri\-ction $\xi':V\setminus E\to T\setminus(W\times\{0\})$ of $\xi$ is a $\Q$-biregular isomorphism {as} $T\setminus(W\times\{0\})=T\cap((\R^n\setminus W)\times\R^m)=\psi^{-1}(\Gamma\cap((\R^n\setminus W)\times\R^m))$ and so the $\Q$-regular map $T\setminus(W\times\{0\})\to V\setminus E$, $(x,y)\mapsto\frac{y}{w(x)}$ is the well-defined inverse of $\xi'$. In particular, the semialgebraic sets $V\setminus E\subset\R^m$ and $T\setminus(W\times\{0\})\subset\R^{n+m}$ are semialgebraically homeo\-morphic, so $\dim(V\setminus E)=\dim(T\setminus(W\times\{0\}))$. It follows that $\dim(T)=\max\{\dim(T\setminus(W\times\{0\})),\dim(W\times\{0\})\}=\max\{\dim(V\setminus E),\dim(W)\}$. If in addition $\max\{\dim(E),\dim(W)\}\leq\dim(V\setminus E)$, then $\dim(V)=\dim(V\setminus E)=\dim(T)$ and thus $T\setminus(W\times\{0\})\subset\Reg^\Q(T)$ by Lemma \ref{861}.

Consider the disjoint union map $\xi\sqcup\eta:V\sqcup W\to T$ of $\xi$ and $\eta$. By construction, we have: $(\xi\sqcup\eta)^{-1}((x,y))=\big\{\frac{y}{w(x)}\big\}\subset V\setminus E$ if $(x,y)\in T\setminus(W\times\{0\})$, and $(\xi\sqcup\eta)^{-1}((x,0))=f^{-1}(x)\sqcup\{x\}\subset E\sqcup W$ if $(x,0)\in W\times\{0\}$. This proves that the family of the fibers of $\xi\sqcup\eta$ coincides with the family of the $\sim$-equivalence classes, so $\xi\sqcup\eta$ descends to the topological quotient, inducing a bijective continuous map $h:V\cup_fW\to T$ such that $h\circ\Pi=\xi\sqcup\eta$. Since $\xi\sqcup\eta$ is a closed map, we deduce that $h$ is a homeomorphism, as required.
\end{proof}


\section{Proofs of Theorems \ref{thm:main}, \ref{thm:main-2} and \ref{thm:main-germs}}\label{sec:proofs}
This section is divided into three parts. In the first two parts, we present some additional preparatory results that we use in the third part to prove our theorems.

\subsection{A relative Nash approximation theorem}
\label{BFR}
Let $N\subset\R^n$ be a Nash manifold, i.e., a locally closed semialgebraic subset of $\R^n$ that is also a smooth submanifold of $\R^n$, let $X$ be a semialgebraic subset of $N$ and let $Y$ be a semialgebraic subset of $\R^m$. A map $f:X\to Y$ is called semialgebraic if its graph is a semialgebraic subset of $\R^{n+m}=\R^n\times\R^m$. The map $f:X\to Y$ is said to Nash if there exist an open semialgebraic neighborhood $U$ of $X$ in $N$ and a map $F:U\to\R^m$ such that $F$ is semialgebraic and smooth, and $F(x)=f(x)$ for all~$x\in X$. A~set $X\subset N$ is said to be Nash if it is the zero set of a finite family of Nash functions defined on~$N$. Each Nash set $X\subset N$ decomposes into the finite union of its Nash irreducible components, see \cite[Cor.8.6.8]{BCR}. A Nash set $X\subset N$ is a Nash submanifold of $N$ if $X\subset\R^n$ is a Nash manifold.

\begin{defn}[{\cite[Def.1.1\,\&\,1.3 and p.63]{BaFeRu2014}}]\label{def:Nmsing}
A set $X\subset N$ is called \emph{Nash set with monomial singularities} if it is a Nash set and, for every $x\in X$, there exist a semialgebraic open neighborhood $U$ of $x$ in $N$ and a Nash diffeomorphism $u:U\to\R^d$, where $d=\dim(N)$, such that $u(x)=0$ and $u(X\cap U)$ is equal to a union of coordinate linear subspaces of $\R^d$. If in addition the Nash irreducible components of $X$ are Nash submanifolds of $N$, then $X\subset N$ is called \emph{Nash monomial crossings}. $\sqbullet$
\end{defn}

The finite unions of Nash submanifolds of $N$ in general position are Nash monomial crossings.

Our next result is a smooth variant of \cite[Thm.1.7\;\&\;Prop.8.2]{BaFeRu2014}.

\begin{thm}\label{thm:1.7}
Let $N\subset\R^n$ and $M\subset\R^m$ be Nash manifolds, let $X\subset N$ and $Y\subset M$ be Nash monomial crossings and let $f:N\to M$ be a smooth map such that $f(X)\subset Y$ and the restriction of $f$ from $X$ to $Y$ preserves Nash irreducible components, i.e., the image of each Nash irreducible component of $X$ under $f$ is contained in some Nash irreducible component of~$Y$. Then there exists a Nash map $g:N\to M$ arbitrarily $\cinfty$ close to $f$ such that $g(X)\subset Y$ and the restriction of $g$ from $X$ to $Y$ preserves Nash irreducible components.
\end{thm}

Since the set of smooth diffeomorphisms between compact smooth manifolds is open with respect to the weak $\cinfty$ topology (see \cite[Thm.1.7, p.38]{hirsch:difftop}) and the Inverse Function Theorem holds true in the Nash category (see \cite[Prop.2.9.7]{BCR}), an immediate consequence reads as follows.

\begin{cor}\label{cor:1.7-weak}
Let $N\subset\R^n$ and $M\subset\R^m$ be compact Nash manifolds, let $I$ be a finite set, let $\{N_i\}_{i\in I}$ be a family of Nash submanifolds of $N$ in general position, let $\{M_i\}_{i\in I}$ be a family of Nash submanifolds of $M$ in general position and let $\phi:N\to M$ be a smooth diffeomorphism such that $\phi(N_i)=M_i$ for all $i\in I$. Then there exists a Nash diffeomorphism $\psi:N\to M$ arbitrarily $\cinfty$~close to $\phi$ such that $\psi(N_i)=M_i$ for all $i\in I$.
\end{cor}

The proof of Theorem \ref{thm:1.7} cannot be obtained by an automatic reformulation of the proofs of \cite[Thm.1.7\;\&\;Prop.8.2]{BaFeRu2014}. The reason is that the results contained in \cite{BaFeRu2014} deal with Nash approximations of semialgebraic $\cnu$ maps between Nash sets with monomial singularities, where $\nu$ is a positive natural number. The case $\nu=\infty$ is not considered in \cite{BaFeRu2014} because a smooth semialgebraic map defined on a Nash manifold is precisely a Nash map by definition. In Appendix~\ref{appendix-B}, we will provide the proof of Theorem \ref{thm:1.7}.


\subsection{Small Lagrange-type interpolations}\label{lagrange} In this subsection, $x$ denotes the coordinate of~$\R$ and $\x$ the corresponding indeterminate. Let $\beta\in\N^*$ and let $c=(c_1,\ldots,c_\beta)\in\R^\beta$ be such that $c_j\neq c_{j'}$ for all $j,j'\in\{1,\ldots,\beta\}$ with $j\neq j'$. We denote $\Ll_{c,1},\ldots,\Ll_{c,\beta}\in\R[\x]$ the \emph{Lagrange basis polynomials associated to $c$}, i.e.,
\begin{equation}\label{eq:lagrange}
\Ll_{c,j}(\x):={\textstyle{\displaystyle\prod}_{s\in\{1,\ldots,\beta\}\setminus\{j\}}}\left(\frac{\x-c_s}{c_j-c_s}\right).
\end{equation}
If $d=(d_1,\ldots,d_\beta)\in\R^\beta$, then $\sum_{j=1}^\beta d_j\Ll_{c,j}$ is the Lagrange interpolation polynomial associated to $c$ and $d$. This polynomial coincides with the unique polynomial $\Ll\in\R[\x]$ of degree $<\beta$ such that $\Ll(c_j)=d_j$ for all $j\in\{1,\ldots,\beta\}$. Furthermore, if each $d_j$ is sufficiently small, then $\mc{L}$ is arbitrarily $\cinfty$ small on compact subsets of $\R$. We need a variant of this result in which we permit that the interpolating function $\Ll$ is regular (not only polynomial), it vanishes on a finite set disjoint from $\{c_1,\ldots,c_\beta\}$ and its smallness is also controlled at infinity.

Set $\R^+:=\{t\in\R:t>0\}$. The mentioned variant is as follows.

\begin{prop}\label{prop:small-lagrange}
Let $A$ be a finite subset of $\R$, let $\beta\in\N^*$ and let $b_1,\ldots,b_\beta\in\R$ be real numbers such that $b_j\neq b_{j'}$ for all $j,j'\in\{1,\ldots,\beta\}$ with $j\neq j'$, and $A\cap\{b_1,\ldots,b_\beta\}=\varnothing$. Let $k,m\in\N$ and let $\epsilon\in\R^+$. Then there exists $\delta\in\R^+$ with the following property: for each $c=(c_1,\ldots,c_\beta)\in\R^\beta$ such that $|b_j-c_j|<\delta$ for all $j\in\{1,\ldots,\beta\}$, there exists a regular function $L_c:\R\to\R$ such that:
\begin{itemize}
\item[$(\mr{i})$] $L_c(a)=0$ for all $a\in A$.
\item[$(\mr{ii})$] $L_c(c_j)=b_j-c_j$ for all $j\in\{1,\ldots,\beta\}$.
\item[$(\mr{iii})$] $|D_hL_c(x)|<\epsilon(1+x^2)^{-k}$ for all $h\in\{0,\ldots,m\}$ and for all $x\in\R$, where $D_h$ denotes the $h^{\mr{th}}$ derivative operator.
\end{itemize} 
\end{prop}
\begin{proof}
Let us assume that $A\neq\varnothing$. If $A=\varnothing$, the proof we present below (suitably simplified) {continues} to work. Let $a_1,\ldots,a_\alpha$ be the elements of $A$. Choose a natural number $\ell$ such that
\[\textstyle
\alpha+\beta-1-2\ell\leq-2k
\]
and define:
\begin{align*}
\delta_1&\textstyle:=\frac{1}{3}\min_{i,j\in\{1,\ldots,\beta\},i\neq j}|b_i-b_j|>0,\\
\delta_2&\textstyle:=\frac{1}{2}\min_{i\in\{1,\ldots,\alpha\},j\in\{1,\ldots,\beta\}}|a_i-b_j|>0,\\
\delta_3&\textstyle:=\min\{\delta_1,\delta_2\}>0,\\
K&\textstyle:=\bigcup_{j=1}^\beta[b_j-\delta_3,b_j+\delta_3],\\
H&\textstyle:=\prod_{j=1}^\beta[b_j-\delta_3,b_j+\delta_3],\\
p(\x)&\textstyle:=\prod_{i=1}^\alpha(\x-a_i)\in\R[\x].
\end{align*}
Observe that $\mr{dist}(A,K)\geq2\delta_2-\delta_3\geq\delta_2>0$. Thus, $p$ never vanishes on the compact subset $K$ of $\R$ so $M:=\sup_{x\in K}|p(x)|^{-1}(1+x^2)^{\ell}$ is finite (and positive). It follows that
\begin{equation}\label{eq:l1}
|p(x)|^{-1}(1+x^2)^{\ell}\leq M\; \text{ for all $x\in K$.}
\end{equation}
For each $j\in\{1,\ldots,\beta\}$ and $w\in\{1,\ldots,\alpha+\beta-1\}$, let $q_{j,w}(\y)\in\R[\y]=\R[\y_1,\ldots,\y_\beta]$ be the unique polynomials such that
\begin{equation}\label{eq:l1'}
\textstyle
p(x)\prod_{s\in\{1,\ldots,\beta\}\setminus\{j\}}(x-y_s)=\sum_{w=0}^{\alpha+\beta-1}q_{j,w}(y)x^w
\end{equation}
for all $(x,y)\in\R\times\R^\beta$. Observe that, for each $j,j'\in\{1,\ldots,\beta\}$ with $j\neq j'$, the intervals $[b_j-\delta_3,b_j+\delta_3]$ and $[b_{j'}-\delta_3,b_{j'}+\delta_3]$ are disjoint since their distance is $\geq3\delta_1-2\delta_3\geq\delta_1>0$. As a consequence, we have that
\[\textstyle
N_{j,w}:=\sup_{(y_1,\ldots,y_n)\in H}\big|q_{j,w}(y)\prod_{s\in\{1,\ldots,\beta\}\setminus\{j\}}(y_j-y_s)^{-1}\big|
\]
is finite (and positive) for all $j\in\{1,\ldots,\beta\}$ and $w\in\{1,\ldots,\alpha+\beta-1\}$. Set
\[\textstyle
N:=\max_{j\in\{1,\ldots,\beta\},w\in\{1,\ldots,\alpha+\beta-1\}}N_{j,w}>0.
\]
It follows that
\begin{equation}\label{eq:l2}
\textstyle
\big|q_{j,w}(y)\prod_{s\in\{1,\ldots,\beta\}\setminus\{j\}}(y_j-y_s)^{-1}\big|\leq N 
\end{equation}
for all $j\in\{1,\ldots,\beta\}$, $w\in\{1,\ldots,\alpha+\beta-1\}$ and $y\in H$. Given $w\in\{1,\ldots,\alpha+\beta-1\}$, let $g_w:\R\to\R$ be the smooth function defined by $g_w(x):=x^w(1+x^2)^{-\ell}$. By elementary considerations, we see that, for each $h\in\{0,\ldots,m\}$, there exists a constant $L_{w,h}>0$ such that $|D_hg_w(x)|\leq L_{w,h}(1+x^2)^{(w-2\ell-h)/2}$ for all $x\in\R$. Set
\[\textstyle
L:=\max_{w\in\{1,\ldots,\alpha+\beta-1\},h\in\{0,\ldots,m\}}L_{w,h}>0.
\]
Since $\alpha+\beta-1-2\ell\leq-2k$, we have that
\[\textstyle
|D_hg_w(x)|\leq L_{w,h}(1+x^2)^{(w-2\ell-h)/2}\leq L(1+x^2)^{(\alpha+\beta-1-2\ell)/2}\leq L(1+x^2)^{-k}
\]
for all $x\in\R$. It follows that
\begin{equation}\label{eq:l3}
|D_hg_w(x)|\leq  L(1+x^2)^{-k}
\end{equation}
for all $w\in\{1,\ldots,\alpha+\beta-1\}$, $h\in\{0,\ldots,m\}$ and $x\in\R$. Set
\[\textstyle
\delta:=\min\{\delta_3,(2\beta(\alpha+\beta-1)MNL)^{-1}\epsilon\}>0.
\]

Let $c_j\in\R$ such that
\begin{equation}\label{eq:l4}
\text{$|b_j-c_j|<\delta$ for each $j\in\{1,\ldots,\beta\}$.}
\end{equation}
Observe that
\begin{equation}\label{eq:l5}
\text{$c_j\in K$ for all $j\in\{1,\ldots,\beta\}$, and  $c:=(c_1,\ldots,c_\beta)$ belongs to $H$,}
\end{equation}
{as} $\delta\leq\delta_3$. Define the regular function $L_c:\R\to\R$ by setting
\[\textstyle
L_c(x):=\frac{p(x)}{(1+x^2)^\ell}\left(\sum_{j=1}^\beta\frac{(1+c_j^2)^\ell}{p(c_j)}(b_j-c_j)\Ll_{c,j}(x)
\right).
\]
Evidently, $L_c$ satisfies points $(\mr{i})$ and $(\mr{ii})$. Let us prove point $(\mr{iii})$. By \eqref{eq:lagrange} and \eqref{eq:l1'}, for every $x\in\R$, we have
\begin{align*}
L_c(x)&\textstyle=\sum_{j=1}^\beta(b_j-c_j)\frac{(1+c_j^2)^\ell}{p(c_j)}\frac{1}{\prod_{s\in\{1,\ldots,\beta\}\setminus\{j\}}(c_s-c_j)}\left(p(x)\prod_{s\in\{1,\ldots,\beta\}\setminus\{j\}}(x-c_j)\right)\frac{1}{(1+x^2)^\ell}=\\
&\textstyle=\sum_{j=1}^\beta\sum_{w=1}^{\alpha+\beta-1}(b_j-c_j)\frac{(1+c_j^2)^\ell}{p(c_j)}\frac{q_{j,w}(c)}{\prod_{s\in\{1,\ldots,\beta\}\setminus\{j\}}(c_s-c_j)}x^w(1+x^2)^{-\ell}.
\end{align*}
As a consequence, by \eqref{eq:l1}, \eqref{eq:l2}, \eqref{eq:l3}, \eqref{eq:l4} and \eqref{eq:l5}, it follows that
\begin{align*}
|D_hL_c(x)|&\textstyle\leq\sum_{j=1}^\beta\sum_{w=1}^{\alpha+\beta-1}\delta MNL(1+x^2)^{-k}=\\
&\textstyle=\beta(\alpha+\beta-1)MNL\delta(1+x^2)^{-k}\leq\frac{\epsilon}{2}(1+x^2)^{-k}<\epsilon(1+x^2)^{-k}
\end{align*}
for all $h\in\{0,\ldots,m\}$ and for all $x\in\R$, as required.
\end{proof}

Let $\Nn(\R)$ be the set of real-valued Nash functions on $\R$ and let $\tau$ be the topology of $\Nn(\R)$ for which a fundamental system of neighborhoods of $f\in\Nn(\R)$ is given by the sets
\[
\textstyle
\mr{U}_{k,m,\epsilon}(f):=\big\{g\in\Nn(\R):\big|D_h(g-f)(x)\big|<\epsilon(1+x^2)^{-k}\;\; \forall h\in\{0,\ldots,m\}, \forall x\in\R\big\},
\]
where $k$ and $m$ vary in $\N$ and $\epsilon$ in $\R^+$. The topology $\tau$ on $\Nn(\R)$ coincides with the `$\cinfty$ topology' on $\mr{N}^\omega(\R)=\Nn(\R)$ defined in \cite[\S II.1]{Sh}. Endow $\Nn(\R)$ with such a topology $\tau$. Let $\EuScript{D}$ be the subset of $\Nn(\R)$ of all Nash diffeomorphisms from $\R$ to $\R$. By \cite[Lem.II.1.7]{Sh}, $\EuScript{D}$ is open in $\Nn(\R)$ and the map $\mr{Inv}:\EuScript{D}\to\EuScript{D}$, sending $f$ into $f^{-1}$, is continuous with respect to the relative topology induced by $\tau$ on~$\EuScript{D}$.

Two consequences of Proposition \ref{prop:small-lagrange} are as follows.

\begin{cor}\label{cor:lagrange1}
Let $A$ and $B$ be two finite subsets of $\R$ such that $A\subset\Q$ and $B\subset\R\setminus\Q$. Then, for every neighborhood $\mc{U}$ of $\mr{id}_\R$ in $\Nn(\R)$, there exists a Nash diffeomorphism $\varphi:\R\to\R$ such that $\varphi\in\mc{U}$, $\varphi(a)=a$ for all $a\in A$, $\varphi(B)\subset\Q$ and $\varphi^{-1}:\R\to\R$ is a regular map.
\end{cor}
\begin{proof}
If $B=\varnothing$, then it suffices to set $\varphi:=\mr{id}_\R$. Suppose that $B\neq\varnothing$. Let $b_1,\ldots,b_\beta$ be the elements of $B$. Let $\mc{U}$ be a neighborhood of $\mr{id}_\R$ in $\Nn(\R)$. Choose a neighborhood $\mc{V}$ of $\mr{\id}_\R$ in $\Nn(\R)$ such that $\mc{V}\subset\EuScript{D}$ and $\mr{Inv}(\mc{V})\subset\mc{U}$. Shrinking $\mc{V}$ around $\mr{\id}_\R$ if necessary, we can assume that there exist $k,m\in\N$ and $\epsilon\in\R^+$ such that $\mc{V}=\mr{U}_{k,m,\epsilon}(\id_\R)$. Let $\delta\in\R^+$ be a positive real number with the properties $(\mr{i})$, $(\mr{ii})$ and $(\mr{iii})$ described in Proposition \ref{prop:small-lagrange}. Choose $c=(c_1,\ldots,c_\beta)\in\Q^\beta$ in such a way that $|b_j-c_j|<\delta$ for all $j\in\{1,\ldots,\beta\}$, and consider the regular function $L_c:\R\to\R$ given by the mentioned Proposition \ref{prop:small-lagrange}. Define the regular function $\psi:\R\to\R$ by setting $\psi:=\mr{id}_\R+L_c$. We know that $L_c(a)=0$ for all $a\in A$, $L_c(c_j)=b_j-c_j$ for all $j\in\{1,\ldots,\beta\}$ and $\psi\in\mc{V}$. It follows that $\psi(a)=a$ for all $a\in A$ and $\psi(c_j)=b_j$ for all $j\in\{1,\ldots,\beta\}$. Moreover, since $\mc{V}\subset\EuScript{D}$ and $\mr{Inv}(\mc{V})\subset\mc{U}$, we have that $\psi$ is a Nash diffeomorphism such that $\varphi:=\psi^{-1}\in\mc{U}$. The Nash diffeomorphism $\varphi$ has all the required properties. 
\end{proof}

\begin{cor}\label{cor:lagrange2}
Let $n\in\N^*$, and let $A$ and $B$ be two finite subsets of $\R^n$ such that $A\subset\Q^n$ and $B\subset\R^n\setminus\Q^n$. Then there exists a Nash diffeomorphism $\varphi:\R^n\to\R^n$ arbitrarily $\cinfty$ close to $\mr{id}_{\R^n}$ such that $\varphi(a)=a$ for all $a\in A$, $\varphi(B)\subset\Q^n$ and $\varphi^{-1}:\R^n\to\R^n$ is a regular map. More precisely, for every neighborhood $\mc{U}$ of $\mr{id}_{\R^n}$ in $\cinfty(\R^n,\R^n)$, there exists a Nash diffeomorphism $\varphi:\R^n\to\R^n$ such that $\varphi\in\mc{U}$, $\varphi(a)=a$ for all $a\in A$, $\varphi(B)\subset\Q^n$ and $\varphi^{-1}:\R^n\to\R^n$ is a regular map.
\end{cor}
\begin{proof}
For each $i\in\{1,\ldots,n\}$, let $\pi_i:\R^n\to\R$ be the projection $\pi_i(x_1,\ldots,x_n):=x_i$, let $A_i:=\pi_i(A\cup B)\cap\Q$ and let $B_i:=\pi_i(A\cup B)\setminus\Q$. By Corollary \ref{cor:lagrange1}, there exists a Nash diffeomorphism $\varphi_i:\R\to\R$ arbitrarily close to $\mr{id}_\R$ in $\Nn(\R)$ such that $\varphi_i(a)=a$ for all $a\in A_i$, $\varphi_i(B_i)\subset\Q$ and $\varphi_i^{-1}$ is a regular function. Since $\tau$ is finer than the relative topology induced by $\cinfty(\R)$ on $\Nn(\R)$, we can assume that each $\varphi_i$ is arbitrarily $\cinfty$ close to $\mr{id}_\R$. Define the Nash diffeomorphism $\varphi:\R^n\to\R^n$ by setting $\varphi(x_1,\ldots,x_n):=(\varphi_1(x_1),\ldots,\varphi_n(x_n))$. Observe that $\varphi(a)=a$ for all $a\in A$ and $\varphi(B)\subset\Q^n$. Furthermore, the $i^{\mr{th}}$ component of $\varphi$ equals $(\pi_i)^*(\varphi_i)=\varphi_i\circ\pi_i$, where $(\pi_i)^*:\cinfty(\R)\to\cinfty(\R^n)$ is the pullback map associated to $\pi_i$. Since each pullback map $(\pi_i)^*$ is continuous with respect to the weak $\cinfty$ topology, we may assume that each component of $\varphi$ is arbitrarily $\cinfty$ close to $\id_\R$, which is equivalent to assume that $\varphi$ is arbitrarily $\cinfty$ close to $\id_{\R^n}$.
\end{proof}

\subsection{Proofs of the theorems}

Recall that, for every $n,m\in\N$ with $m>n$, we identify $\R^n$ with the subset $\R^n\times\{0\}$ of $\R^n\times\R^{m-n}=\R^m$, so we can write $\R^n\subset\R^m$ and every subset of $\R^n$ is also a subset of $\R^m$. If $M$ and $N$ are smooth manifolds, then $\cinfty(M,N)$ is the set of smooth maps from $M$ to $N$ endowed with the weak $\cinfty$ topology. If $S$ is a topological space, then $\czero(S,\R^m)$ is the set of continuous maps from $S$ to $\R^m$ endowed with the compact-open topology. A Nash embedding is a smooth embedding that is also a Nash map, and a continuous embedding is a homeomorphism onto its image.

We need the last preliminary result.

\begin{lem}\label{lem:degen}
Let $Z\subset\R^n$ be a $\Q$-algebraic set of positive dimension and let $H$ be a finite subset of $Z$ such that $\Sing^\Q(Z)\subset H\subset\Q^n$. Denote $J:\R^n\to\R^{n+1}=\R^n\times\R$ the inclusion map $x\mapsto(x,0)$ and choose a neighborhood $\mc{U}$ of $J$ in $\czero(\R^n,\R^{n+1})$ and a neighborhood $\mc{V}$ of $J|_{\R^n\setminus H}$ in $\cinfty(\R^n\setminus H,\R^{n+1})$. Then there exist a semialgebraic continuous embedding $\nu:\R^n\to\R^{n+1}$ with the following properties:
\begin{itemize}
 \item[$(\mr{i})$] $\nu(x)=(x,0)$ for all $x\in H$.
 \item[$(\mr{ii})$] $Z':=\nu(Z)$ is a $\Q$-determined $\Q$-algebraic subset of $\R^{n+1}$ such that $\nu(Z\setminus H)=\Reg(Z')$.
 \item[$(\mr{iii})$] The restriction $\nu|_{\R^n\setminus H}:\R^n\setminus H\to\R^{n+1}$ is a Nash embedding.
 \item[$(\mr{iv})$] $\nu\in\mc{U}$ and $\nu|_{\R^n\setminus H}\in\mc{V}$.
\end{itemize}
\end{lem}
\begin{proof}
Let $I$ be the set of isolated points of $Z$. We have: $I\subset\Sing^\Q(Z)\subset H$. If $H=I$, then $I=\Sing(Z)=\Sing^\Q(Z)$ so the $\Q$-algebraic set $Z\subset\R^n$ is already $\Q$-determined and the result is trivial: it suffices to set $\nu:=J$. Suppose that $H\neq I$. Let $c$ be a positive small rational number, let $P\in\Q[\x]$ be the unique polynomial such that $P(x)=c\prod_{p\in H}|x-p|_n^2$ for all $x\in\R^n$, let $Q:=\x_{n+1}^3-P(\x)\in\Q[\x,\x_{n+1}]$, let $Y\subset\R^{n+1}$ be the $\Q$-algebraic set $\ZZ_\R(Q)$, let $\pi:Y\to\R^n$ be the projection $\pi(x,x_{n+1}):=x$, let $\nu:\R^n\to\R^{n+1}$ be the map $\nu(x):=(x,\sqrt[3]{P(x)})$ and let $Z':=\nu(Z)$. Observe that $\ZZ_\R(P)=H$ and $\nu|_{\R^n\setminus H}$ is a Nash embedding. If $c$ is sufficiently small, then $\nu$ is arbitrarily $\czero$ close to $J$ and $\nu|_{\R^n\setminus H}$ is arbitrarily $\cinfty$ close to $J|_{\R^n\setminus H}$. The set $Z'\subset\R^{n+1}$ is $\Q$-algebraic {as} $Z'=\pi^{-1}(Z)$. Since $\frac{\partial Q}{\partial x_{n+1}}(x,x_{n+1})\neq0$ for all $(x,x_{n+1})\in Y\setminus(H\times\{0\})$, $\pi$ is submersive at every point of $Y\setminus(H\times\{0\})$ and {therefore} $\nu(Z\setminus H)=\pi^{-1}(Z\setminus H)\subset\Reg^\Q(Z')$ by Proposition \ref{prop:Q-transverse}. In order to complete the proof, it remains to show that $H\times\{0\}\subset\Sing(Z')$ so $Z'\subset\R^{n+1}$ is $\Q$-determined. Choose $q\in H$. If $q\in I$, then $(q,0)$ is an isolated point of $Z'$ so it belongs to $\Sing(Z')$ {since} $Z'$ has positive dimension. Let $q\in H\setminus I$. Suppose that $(q,0)\in\Reg(Z')$. Since $(q,0)$ is not an isolated point of $Z'$, by the curve selection lemma, there exists a Nash curve $(\alpha,\beta):(-1,1)\to\R^{n+1}=\R^n\times\R$ whose image is contained in $Z'$ (and thus in $Y$) such that $(\alpha(0),\beta(0))=(q,0)$ and $(\alpha'(0),\beta'(0))\neq0$. Observe that $\beta(t)=\sqrt[3]{P(\alpha(t))}\geq0$ for all $t\in(-1,1)$, so $\beta'(0)=0$ and $\alpha'(0)\neq0$. It follows that $\lim_{t\to0}\beta(t)^3t^{-2}=0$ and $\lim_{t\to0}P(\alpha(t))t^{-2}=|\alpha'(0)|_n^2\prod_{p\in H\setminus\{q\}}|q-p|_n^2\neq0$. This is a contradiction because $\beta(t)^3=P(\alpha(t))$ for all $t\in(-1,1)$ so $0=\lim_{t\to0}\beta(t)^3t^{-2}=\lim_{t\to0}P(\alpha(t))t^{-2}\neq0$. 
\end{proof}

Let us give the proof of Theorem \ref{thm:main}. First, we deal with the compact case.

\subsubsection{Proof of Theorem \ref{thm:main}: the compact case.} \label{431} Let $X\subset\R^n$ be a compact algebraic set of dimension $d$ such that $\Sing(X)$ is finite. If $d=0$, then $X$ is finite and it suffices to apply Corollary \ref{cor:lagrange2} with $A:=X\cap\Q^n$ and $B:=X\setminus\Q^n$.

Suppose that $d\geq1$. We divide the rest of the proof into six steps.

\textsc{Step I.} Let $I$ be the set of isolated points of $X$. Since $d\geq1$, $I$ is a subset of $\Sing(X)$. Define $A:=\Sing(X)\cap\Q^n$ and $B:=\Sing(X)\setminus\Q^n$. By Corollary~\ref{cor:lagrange2}, there exists a Nash diffeomorphism $\varphi:\R^n\to\R^n$ arbitrarily $\cinfty$ close to $\mr{id}_{\R^n}$ such that $\varphi(a)=a$ for all $a\in A$, $\varphi(B)\subset\Q^n$ and $\varphi^{-1}$ is a regular map. Define $X_1:=\varphi(X)$. Since $X_1=(\varphi^{-1})^{-1}(X)$, it follows that $X_1\subset\R^n$ is an algebraic set and, by transversality (see \cite[Lem.2.2.13]{akbking:tras}), $\varphi(\Reg(X))\subset\Reg(X_1)$ or, equivalently, $\Sing(X_1)\subset\varphi(\Sing(X))$. Let $I_1:=\varphi(I)$. Observe that $I_1$ is the set of isolated points of $X_1$ and $I_1\subset\Sing(X_1)\subset\varphi(\Sing(X))\subset\Q^n$.

\textsc{Step II.}  If $\Sing(X_1)=I_1$, then $\Reg(X_1)=X_1\setminus\Sing(X_1)$ is a compact Nash submanifold of $\R^n$ and, to complete the proof, it suffices to apply Theorem~\ref{thm:NTQ} and Remarks~\ref{rem17}$(\mr{i})(\mr{ii})$ to~$\Reg(X_1)$.

Suppose that $\Sing(X_1)\setminus I_1=\{p_1,\ldots,p_s\}$ for some $s\in\N^*$. By Hironaka's resolution of singularities, there exist a compact nonsingular algebraic set $Y\subset\R^e$ of dimension $d$ for some $e\in\N^*$, a regular map $\pi:Y\to X_1$ and a family $\{Y_i\}_{i=1}^s$ of algebraic subsets of $Y$ such that $\pi(Y)=X_1\setminus I_1$, $\pi^{-1}(p_i)=Y_i$ for all $i\in\{1,\ldots,s\}$, the restriction of $\pi$ from $Y\setminus\bigcup_{i=1}^sY_i$ to $\Reg(X_1)=X_1\setminus\Sing(X_1)$ is a biregular isomorphism and each $Y_i$ is a finite union of nonsingular algebraic hypersurfaces of $Y$ in general position. In addition, using the standard generic proje\-ction theorem, we can assume that $e=2d+1$, i.e., $Y$ is a compact nonsingular algebraic subset of $\R^{2d+1}$.

\textsc{Step III.} By Corollary \ref{thm:Q_tico_approx}, there exist a projectively $\Q$-closed $\Q$-nonsingular $\Q$-algebraic set $V\subset\R^{2d+1}$, a family $\{V_i\}_{i=1}^s$ of $\Q$-algebraic subsets of $V$ and a smooth diffeomorphism $\psi:Y\to V$ such that each $V_i$ is a finite union of $\Q$-nonsingular $\Q$-algebraic hypersurfaces of $V$ in general position and $\psi(Y_i)=V_i$ for every $i\in\{1,\ldots,s\}$. By Corollary \ref{cor:1.7-weak}, we can also assume that $\psi$ is a Nash diffeomorphism.

Let $g:V\to\R^n$ be the Nash map $g(y):=\pi(\psi^{-1}(y))$ and let $E:=\bigcup_{i=1}^sV_i$. Observe that:
\begin{itemize}
 \item $g(V)=X_1\setminus I_1$ and $g^{-1}(p_i)=V_i$ for all $i\in\{1,\ldots,s\}$,
 \item the restriction $g|_E:E\to\R^n$ is $\Q$-regular {as} $V_i\subset\R^{2d+1}$ is $\Q$-algebraic, $g(V_i)=\{p_i\}$ and $p_i\in\Q$ for all $i\in\{1,\ldots,s\}$,
 \item the map $\widehat{g}:V\setminus E\to(X_1\setminus I_1)\setminus\{p_1,\ldots,p_s\}$, defined by $\widehat{g}(y):=g(y)$, is a Nash diffeomorphism.
\end{itemize}

\textsc{Step IV.} Let $y=(y_1,\ldots,y_{2d+1})$ be the coordinates of $\R^{2d+1}$, let $\y=(\y_1,\ldots,\y_{2d+1})$ be the corresponding indeterminates and let $q\in\Q[\y]$ be a polynomial such that $\ZZ_\R(q)=E$. Choosing a sufficiently small non-zero rational number $c$ and replacing $q$ with $cq$ if necessary, we can assume that the polynomial function $q:\R^{2d+1}\to\R$ corresponding to $q\in\Q[\y]$ is arbitrarily $\cinfty$ close to the zero function. By Lemma \ref{lem:Q-basic}$(\mr{iii})$, there exists a $\Q$-regular map $g^*:\R^{2d+1}\to\R^n$ such that $g^*|_E=g|_E$. Using a tubular neighborhood of $V$ in $\R^{2d+1}$ and a suitable smooth partition of unity, we can construct a smooth map $g^{**}:\R^{2d+1}\to\R^n$ such that $g^{**}=g$ on $V$ and $g^{**}=g^*$ outside of a compact neighborhood of $V$ in $\R^{2d+1}$. By Lemma \ref{lem:35}, $E\subset\R^{2d+1}$ is a $\Q$-stable $\Q$-algebraic set so we can apply Lemma \ref{lem:L} to the components of $g^{**}-g^*$ with $L:=E$, obtaining a $\Q$-regular map $h:\R^{2d+1}\to\R^n$ arbitrarily $\cinfty$ close to $g^{**}-g^*$ that vanishes on $E$. Define the $\Q$-regular map $G:V\to\R^n$ by $G(y):=h(y)+g^*(y)$. Observe that $G$ is arbitrarily $\cinfty$ close to $g$ and $G|_E=g|_E$. Define the $\Q$-regular map $F:V\to\R^{n+1}=\R^n\times\R$ and the $\Q$-algebraic set $W\subset\R^{n+1}=\R^n\times\R$ as follows:
\begin{itemize}
 \item $F:=G\times q$, i.e., $F(y):=(G(y),q(y))$ for every $y\in V$,
 \item $W:=\{(p_1,0),\ldots,(p_s,0)\}\subset\Q^{n+1}$.
\end{itemize}
We have $F^{-1}((p_i,0))=G^{-1}(p_i)\cap E=(G|_E)^{-1}(p_i)=(g|_E)^{-1}(p_i)=V_i$ for all $i\in\{1,\ldots,s\}$ so $F^{-1}(W)=E$. Choosing $G$ sufficiently $\cinfty$ close to $g$ and $q|_V$ sufficiently $\cinfty$ close to the zero function $c_0:V\to\R$, we can assume that $F$ is arbitrarily $\cinfty$ close to $g\times c_0:V\to\R^{n+1}$.

\textsc{Step V.} Let $f:E\to W$, $g':V\to X_1\setminus I_1$ and $w:W\to X_1\setminus I_1$ be the maps $f(y):=F(y)=(g(y),0)$ for all $y\in E$, $g'(y):=g(y)$ for all $y\in V$, and $w(p_i,0):=p_i$ for all $i\in\{1,\ldots,s\}$. Consider the adjunction topological space $V\cup_f W$ and the corresponding natural projection map $\Pi:V\sqcup W\to V\cup_f W$.
Observe that the map $g'\sqcup w:V\sqcup W\to X_1\setminus I_1$ is continuous and surjective, and the family of its fibers coincides with the family of the fibers of $\Pi$. Thus, $g'\sqcup w$ passes to the  quotient as a homeomorphism, i.e., there exists a homeomorphism $k:V\cup_f W\to X_1\setminus I_1$ such that $k\circ \Pi=g'\sqcup w$, i.e., the following diagram commutes:
\[
\centering
\begin{tikzcd}
V\sqcup W \arrow[d,"\Pi"] \arrow[r,"g'\sqcup w"] & X_1\setminus I_1\\
V\cup_f W \arrow[ur, bend right=20, "k"] & 
\end{tikzcd}
\]

Let $x=(x_1,\ldots,x_n)$ be the coordinates of $\R^n$ and let $(x,t,y)$ be the coordinates of $\R^{n+2d+2}=\R^n\times\R\times\R^{2d+1}$. By Lemma \ref{lem:blowing_down}, there exist a $\Q$-algebraic set $T\subset\R^{n+2d+2}$, $\Q$-regular maps $\xi:V\to T$ and $\eta:W\to T$, and a homeomorphism $h :V\cup_f W\to T$ such that $h\circ \Pi=\xi\sqcup\eta$, i.e., the following diagram commutes:
\[
\centering
\begin{tikzcd}
V\sqcup W \arrow[d,"\Pi"] \arrow[r,"\xi\sqcup\eta"] & T \\
V\cup_f W \arrow[ur, bend right=20, "h"] & 
\end{tikzcd}
\]
In addition, we have: $\eta(p_i,0)=(p_i,0,0)\in\R^{n+2d+2}$ for all $(p_i,0)\in W$, $W\times\{0\}\subset\xi(V)=T$, $\xi(V\setminus E)=T\setminus(W\times\{0\})$, the restriction of $\xi$ from $V\setminus E$ onto $T\setminus(W\times\{0\})$ is a Nash diffeomorphism (actually, a $\Q$-biregular isomorphism), $T\setminus(W\times\{0\})\subset\Reg^\Q(T)$ (or, equivalently, $\Sing^\Q(T)\subset W\times\{0\}$) and the map $\overline{\xi}:V\to\R^{n+2d+2}$, $y\mapsto\xi(y)$ is arbitrarily $\cinfty$ close to $V\to\R^{n+2d+2}$, $y\mapsto(F(y),0)$. Since $F$ can be chosen arbitrarily $\cinfty$ close to $g\times c_0$, we can also assume that $\overline{\xi}$ is arbitrarily $\cinfty$ close to  $V\to\R^{n+2d+2}$, $y\mapsto(g(y),0,0)$.

The homeomorphism $\Phi:=h\circ k^{-1}:X_1\setminus I_1\to T$ has the following explicit form:
\begin{itemize}
 \item $\Phi(x)=(x,0,0)$ for all $x\in\{p_1,\ldots,p_s\}$,
 \item $\Phi(x)=\xi(\widehat{g}^{-1}(x))$ for all $x\in(X_1\setminus I_1)\setminus\{p_1,\ldots,p_s\}$.
\end{itemize}
Evidently, $\Phi((X_1\setminus I_1)\setminus\{p_1,\ldots,p_s\})=T\setminus(W\times\{0\})$ and the map $(X_1\setminus I_1)\setminus\{p_1,\ldots,p_s\}\to\R^{n+2d+2}$, $x\mapsto\Phi(x)$ is a Nash embedding. Since $(X_1\setminus I_1)\setminus\{p_1,\ldots,p_s\}$ is dense in $X_1\setminus I_1$, the closure of the graph of the latter Nash embedding in $\R^{2n+2d+2}=\R^n\times\R^{n+2d+2}$ coincides with the graph of~$\Phi$. This proves that $\Phi$ is semialgebraic.

Define the $\Q$-algebraic subsets $I_2$ and $X_2$ of $\R^{n+2d+2}$ by $I_2:=I_1\times\{0\}\times\{0\}$ and $X_2:=T\cup I_2$. 

Let $\overline{\Phi}:X_1\setminus I_1\to\R^{n+2d+2}$ be the map $\overline{\Phi}(x):=\Phi(x)$. Since $\overline{\Phi}(x)=\overline{\xi}(\widehat{g}^{-1}(x))$ for all $x\in(X_1\setminus I_1)\setminus\{p_1,\ldots,p_s\}$, if $\overline{\xi}$ is sufficiently $\cinfty$ close to  $V\to\R^{n+2d+2}$, $y\mapsto(g(y),0,0)$, then we can assume that $\overline{\Phi}|_{(X_1\setminus I_1)\setminus\{p_1,\ldots,p_s\}}$ is arbitrarily $\cinfty$ close to $(X_1\setminus I_1)\setminus\{p_1,\ldots,p_s\}\to\R^{n+2d+2}$, $x\mapsto(g(\widehat{g}^{-1}(x)),0,0)=(x,0,0)$. Since the domain $V$ of $g$ is compact, we can also assume that $\overline{\Phi}$ is arbitrarily $\czero$ close to $X_1\setminus I_1\to\R^{n+2d+2}$, $x\mapsto(x,0,0)$ and thus $T\cap I_2=\varnothing$. It follows that $I_2$ is the set of isolated points of $X_2$ and $I_2\subset\Sing(X_2)\subset\Sing^\Q(X_2)\subset I_2\sqcup(W\times\{0\})$.

Define the semialgebraic continuous embedding $\phi_2:X\to\R^{n+2d+2}$ by $\phi_2(x):=\overline{\Phi}(\varphi(x))$ if $x\in X\setminus I$ and $\phi_2(x):=(\varphi(x),0,0)$ if $x\in I$. Observe that $\phi_2(X)$ is equal to the $\Q$-algebraic set $X_2\subset\R^{n+2d+2}$, $\Sing^\Q(X_2)\subset I_2\sqcup(W\times\{0\})\subset\phi_2(\Sing(X))$, and the restriction $\phi_2|_{\Reg(X)}:\Reg(X)\to\R^{n+2d+2}$ is a Nash embedding. In addition, if $\varphi$ is sufficiently $\cinfty$ close to $\mr{id}_{\R^n}$, then $\phi_2$ is arbitrarily $\czero$ close to the inclusion $X\hookrightarrow\R^{n+2d+2}$, $x\mapsto(x,0,0)$ and $\phi_2|_{\Reg(X)}$ is arbitrarily $\cinfty$ close to the inclusion $\Reg(X)\hookrightarrow\R^{n+2d+2}$, $x\mapsto(x,0,0)$. Observe that the map $\phi_2$ does not yet have all the required properties. The problem is that $X_2=\phi_2(X)$ may not be $\Q$-determined: we cannot in fact exclude that the subset $\Sing^\Q(X_2)\setminus\Sing(X_2)$ of $W\times\{0\}$ is non-empty.

\textsc{Step VI.} To overcome this problem, we {may} apply Lemma \ref{lem:degen} with $Z:=X_2\subset\R^{n+2d+2}$ and $H:=I_2\sqcup(W\times\{0\})\subset\Q^{n+2d+2}$, obtaining a semialgebraic continuous embedding $\nu:\R^{n+2d+2}\to\R^{n+2d+3}=\R^{n+2d+2}\times\R$ such that $\nu(p)=(p,0)$ for all $p\in I_2\sqcup(W\times\{0\})$, $X':=\nu(X_2)$ is a $\Q$-determined $\Q$-algebraic subset of $\R^{n+2d+3}$, $\nu(X_2\setminus(I_2\sqcup(W\times\{0\})))=\Reg(X')$, $\nu$ is arbitrarily $\czero$ close to the inclusion $J:\R^{n+2d+2}\hookrightarrow\R^{n+2d+3}$, $p\mapsto(p,0)$ and $\nu|_{\R^{n+2d+2}\setminus (W\times\{0\})}$ is arbitrarily $\cinfty$ close to $J|_{\R^{n+2d+2}\setminus (W\times\{0\})}$. The semialgebraic continuous embedding $\phi:=\nu\circ\phi_2:X\to\R^{n+2d+3}$ has all the required properties, and the proof of Theorem \ref{thm:main} in the compact case is complete.

\vspace{.5em}

Two remarks concerning the latter proof are in order.

\begin{remark}\label{remark48}
Let $X\subset\R^n$ be a compact algebraic set of dimension $d\geq1$ such that $\Sing(X)$ is finite. Observe that the preceding proof actually gives a compact $\Q$-determined $\Q$-algebraic set $X'\subset\R^m$ satisfying Theorem \ref{thm:main}$(\mr{i})(\mr{ii})(\mr{iii})$ with $m:=n+2d+3$, one additional dimension will be needed just if $X$ is noncompact. $\sqbullet$
\end{remark}

\begin{remark} \label{rem:degen}
Consider again a compact algebraic set $X\subset\R^n$ of dimension $d\geq1$ such that $\Sing(X)$ is finite, and choose any (possibly empty) finite subset $P$ of $X\cap\Q^n$. In the statement of Theorem \ref{thm:main}, we can replace properties $(\mr{i})$ and $(\mr{ii})$ with the following ones:
\begin{itemize}
 \item[$(\mr{i}')$] \textit{$X':=\phi(X)$ is a $\Q$-determined $\Q$-algebraic subset of $\R^{n+2d+3}$ such that $\phi(x)=(x,0)\in\R^n\times\R^{2d+3}=\R^{n+2d+3}$ for every $x\in P$, and $\Sing(X')=(P\times\{0\})\cup\phi(\Sing(X))$. In particular, if $P=\varnothing$, then $\Sing(X')=\phi(\Sing(X))$ and {therefore} $X$ and $X'$ {have} the same number of singularities.}
 \item[$(\mr{ii}')$] \textit{The restriction $\phi|_{\Reg(X)\setminus P}:\Reg(X)\setminus P\to\R^{n+2d+3}$ is a Nash embedding. In particular, if $P=\varnothing$, then $\phi(\Reg(X))=\Reg(X')$ and the restriction of $\phi$ from $\Reg(X)$ onto $\Reg(X')$ is a Nash diffeomorphism.}
\end{itemize}
To obtain these two properties, we can modify the previous proof as follows:
\begin{itemize}
 \item In \textsc{Step I}, define $A:=P\cup(\Sing(X)\cap\Q^n)$ instead of $A:=\Sing(X)\cap\Q^n$.
 \item In \textsc{Step II}, if $P\cup\varphi(\Sing(X))=I_1$, then we can again complete the proof using Theo\-rem~\ref{thm:NTQ} and Remarks~\ref{rem17}$(\mr{i})(\mr{ii})$. Suppose that $(P\cup\varphi(\Sing(X)))\setminus I_1=\{p_1,\ldots,p_s\}$ for some $s\in\N^*$. Using Hironaka's resolution of singularities, blowing up the strict transform of the finite set $(P\cup\varphi(\Sing(X)))\setminus\Sing(X_1)$ and using the generic projection theorem, we {may} again assume that there exist a compact nonsingular algebraic set $Y\subset\R^{2d+1}$ of dimension~$d$, a regular map $\pi:Y\to X_1$ and a family $\{Y_i\}_{i=1}^s$ of algebraic subsets of $Y$ such that $\pi(Y)=X_1\setminus I_1$, $\pi^{-1}(p_i)=Y_i$ for all $i\in\{1,\ldots,s\}$, the restriction of $\pi$ from $Y\setminus\bigcup_{i=1}^sY_i$ to $X_1\setminus(P\cup\varphi(\Sing(X)))$ is a biregular isomorphism and each $Y_i$ is a finite union of nonsingular algebraic hypersurfaces of $Y$ in general position.
 \item Repeat \textsc{Steps III, IV, V} and \textsc{VI} (with minor obvious changes), obtaining a semialgebraic continuous embedding $\phi:X\to\R^{n+2d+3}$ with the previous properties $(\mr{i}')$ and $(\mr{ii}')$, and with the property $(\mr{iii})$ stated in Theorem \ref{thm:main}.~$\sqbullet$
\end{itemize}
\end{remark}

\subsubsection{Proof of Theorem \ref{thm:main}: the noncompact case} \label{432} Let $X\subset\R^n$ be a noncompact algebraic set of dimension $d$. Since $X$ is infinite, $d$ is positive. If $d=n$, then $X=\R^n$ and there is nothing to prove. 

Suppose that $d<n$. Using a translation of $\R^n$ along a suitable small vector, we can assume that the origin $0$ of $\R^n$ does not belong to $X$. Let $P:=\{0\}$, let $\theta:\R^n\setminus P\to\R^n\setminus P$ be the inversion $\theta(x):=x|x|_n^{-2}$, and let $X_\bullet\subset\R^n$ be the algebraic set defined by $X_\bullet:=P\sqcup\theta(X)$, which is the (algebraic) Alexandrov compactification of $X$. Since $\theta$ is a biregular isomorphism (with $\theta^{-1}=\theta$), it follows that $\Sing(X_\bullet)\subset P\sqcup\theta(\Sing(X))$. Apply the compact case of Theorem~\ref{thm:main} and Remarks~\ref{remark48}\,\&\,\ref{rem:degen} proven above to $X_\bullet$ and $P$, obtaining a semialgebraic continuous embedding $\phi_\bullet:X_\bullet\to\R^{n+2d+3}$ such that $X'_\bullet:=\phi_\bullet(X_\bullet)\subset\R^{n+2d+3}$ is a $\Q$-determined $\Q$-algebraic set, $\phi_\bullet(0)=(0,0)\in\R^n\times\R^{2d+3}=\R^{n+2d+3}$, $\Sing(X'_\bullet)=\{(0,0)\}\cup\phi_\bullet(\Sing(X_\bullet))=\{(0,0)\}\sqcup\phi_\bullet(\theta(\Sing(X)))$, the restriction $\phi_\bullet|_{\theta(\Reg(X))}:\theta(\Reg(X))\to\R^{n+2d+3}$ is a Nash embedding, $\phi_\bullet$ is arbitrarily $\czero$ close to the inclusion $J_\bullet:X_\bullet\to\R^{n+2d+3}$, $x\mapsto(x,0)$ and $\phi_\bullet|_{\theta(\Reg(X))}$ is arbitrarily $\cinfty$ close to the inclusion $J:\theta(\Reg(X))\to\R^{n+2d+3}$, $x\mapsto(x,0)$.

Let $(x,z)=(x_1,\ldots,x_n,z_1,\ldots,z_{2d+3})$ be the coordinates of $\R^{n+2d+3}$, let $Q$ be the singleton $\{(0,0)\}$ of $\R^{n+2d+3}$, let $\Theta:\R^{n+2d+3}\setminus Q\to\R^{n+2d+3}\setminus Q$ be the inversion $\Theta(x,z):=(x,z)|(x,z)|_{n+2d+3}^{-2}$, and let $\phi_1:X\to\R^{n+2d+3}$ be the semialgebraic continuous embedding  $\phi_1(x):=\Theta(\phi_\bullet(\theta(x)))$. Since $\Theta$ is a $\Q$-biregular isomorphism, $X_1:=Q\sqcup\phi_1(X)=Q\sqcup\Theta(X'_\bullet\setminus Q)$ is a $\Q$-algebraic subset of $\R^{n+2d+3}$ such that $\Sing^\Q(X_1)\setminus Q=\Sing(X_1)\setminus Q=\Theta(\phi_\bullet(\theta(\Sing(X))))=\phi_1(\Sing(X))$ and the restriction $\phi_1|_{\Reg(X)}:\Reg(X)\to\R^{n+2d+3}$, $x\mapsto\Theta(\phi_\bullet(\theta(x)))$, is a Nash embedding. Since $\Theta(x,0)=(\theta(x),0)$ for all $x\in\R^n$, if $\phi_\bullet$ is sufficiently $\czero$ close to the inclusion $J_\bullet$ and $\phi_\bullet|_{\theta(\Reg(X))}$ is sufficiently $\cinfty$ close to the inclusion $J$, then $\phi_1$ is arbitrarily $\czero$ close to the inclusion $X\to\R^{n+2d+3}$, $x\mapsto\Theta(J_\bullet(\theta(x)))=(\theta(\theta(x)),0)=(x,0)$ and $\phi_1|_{\Reg(X)}$ is arbitrarily $\cinfty$ close to the inclusion $\Reg(X)\to\R^{n+2d+3}$, $x\mapsto\Theta(J(\theta(x)))=(\theta(\theta(x)),0)=(x,0)$.

Let $c$ be a small non-zero rational number, let $\Psi:\R^{n+2d+3}\setminus Q\to\R$ be the $\Q$-regular function $\Psi(x,z):=c|(x,z)|_{n+2d+3}^{-2}$, let $(x,z,t)$ be the coordinates of $\R^{n+2d+4}=\R^{n+2d+3}\times\R$ and let $S$ be the $\Q$-nonsingular $\Q$-algebraic set $S:=\{(x,y,t)\in\R^{n+2d+4}:t|(x,z)|_{n+2d+3}^2-c=0\}$. Observe that $S$ coincides with the graph of $\Psi$ and $\Psi$ is arbitrarily $\cinfty$ close to the zero function $\R^{n+2d+3}\setminus Q\to\R$, $(x,z)\mapsto0$ if $c$ is sufficiently small. Define the semialgebraic continuous embedding $\phi:X\to\R^{n+2d+4}$ and the $\Q$-algebraic set $X'\subset\R^{n+d+4}$ by $\phi(x):=(\phi_1(x),\Psi(\phi_1(x)))$ and $X':=\phi(X)=S\cap((X_1\setminus Q)\times\R)$. Since the map $X_1\setminus Q\to X'$, $(x,t)\mapsto(x,t,\Psi(x,y))$ is $\Q$-biregular and $\Sing^\Q(X_1)\setminus Q=\Sing(X_1)\setminus Q$, Lemma \ref{861} implies that the $\Q$-algebraic set $X'\subset\R^{n+2d+4}$ is $\Q$-determined. The map $\phi$ has all the required properties. This completes the proof of Theorem \ref{thm:main} in the noncompact case.

\begin{remark}\label{rem:degen2}
The map $\phi:X\to\R^{n+2d+4}$, constructed in the previous proof when $X\subset\R^n$ is noncompact, has the property $\Sing(\phi(X))=\phi(\Sing(X))$, in addition to the three properties stated in Theorem \ref{thm:main}. $\sqbullet$
\end{remark}

\subsubsection{Proof of Theorem \ref{thm:main-2}} {To prove the compact case, we can modify the proof of the compact case of Theorem \ref{thm:main} as follows:}
\begin{itemize}
 \item Repeat the first three steps until you define the projectively $\Q$-closed $\Q$-nonsingular $\Q$-algebraic set $V\subset\R^{2d+1}$, the corresponding family $\{V_i\}_{i=1}^s$ of $\Q$-algebraic subsets of $V$ and $E:=\bigcup_{i=1}^sV_i$. Now denote $g:E\to\R$ the $\Q$-regular map such that $g(V_i)=\{i\}$ for all $i\in\{1,\ldots,s\}$, and $G:V\to\R$ a $\Q$-regular function that extends $g$ (see Lemma~\ref{lem:Q-basic}$(\mr{iii})$). Consider a polynomial $q\in\Q[\y_1,\ldots,\y_{2d+1}]$ such that $\ZZ_\R(q)=E$. Define the $\Q$-regular map $F:V\to\R^2$ by $F(y):=(G(y),q(y))$. Since $F^{-1}((i,0))=V_i$ for all $i\in\{1,\ldots,s\}$, the adjunction space $V\cup_{F|_E}\{(1,0),\ldots,(s,0)\}$ is homeomorphic to $X_1\setminus I_1$.
 \item Repeat \textsc{Steps V} and \textsc{VI} replacing each $p_i$ with $i$, and $n$ with $1$ (and making minor obvious changes). We obtain a $\Q$-determined $\Q$-algebraic subset $X'$ of $\R^{1+2d+3}=\R^{2d+4}$ and a semialgebraic homeomorphism $\eta:X\to X'$ such that $\eta(\Reg(X))\subset\Reg(X')$ and the restriction of $\eta$ from $\Reg(X)$ to $\eta(\Reg(X))$ is a Nash diffeomorphism.
 \item By modifying the previous part of this proof as indicated in Remark \ref{rem:degen} (with $P:=\varnothing$), we can also assume that $\eta(\Reg(X))=\Reg(X')$.
\end{itemize}

To prove the noncompact case of Theorem \ref{thm:main-2}, we can modify the proof of the noncompact case of Theorem \ref{thm:main} using the same strategy employed above for the compact case. We obtain a $\Q$-determined $\Q$-algebraic subset $X'$ of $\R^{1+2d+4}=\R^{2d+5}$ and a semialgebraic homeomorphism $\eta:X\to X'$ with the required properties.

\begin{remark} \label{rem:degen-2}
Let $X\subset\R^n$ be a compact algebraic set of dimension $d\geq1$ such that $\Sing(X)$ is finite, and let $P=\{p_1,\ldots,p_r\}$ be any non-empty finite subset of $X\cap\Q^n$. In the statement of Theorem \ref{thm:main-2}, we can assume that the {$\Q$-determined $\Q$-algebraic} set $X'\subset\R^{2d+4}$ and the map $\eta:X\to X'$ have the following additional properties:
\begin{itemize}
 \item \textit{$X'$ contains the finite subset $\{(1,0),\ldots,(r,0)\}$ of $\R^{2d+4}=\R\times\R^{2d+3}$ and $\eta(p_i)=(i,0)$ for all $i\in\{1,\ldots,r\}$.}
\end{itemize}
To obtain these additional properties, it suffices to combine the strategy described in Remark~\ref{rem:degen} with the previous proof of Theorem \ref{thm:main-2}. $\sqbullet$
\end{remark}

\subsubsection{Proof of Theorem \ref{thm:main-germs}} If $d=0$, the result is evident {since} $O$ is an isolated point of $X$. Suppose that $d\geq1$. Choose $\epsilon>0$ sufficiently small in such a way that the closed ball $\overline{\B}_n(\epsilon)$ of $\R^n$ centered at $O$ of radius $\epsilon$ intersects $\Sing(X)$ only at $O$, its boundary $\partial\overline{\B}_n(\epsilon)$ intersects transversally $\Reg(X)$ in $\R^n$ and $\dim(X\cap\overline{\B}_n(\epsilon))=d$. Let $\pi:\R^{n+1}=\R^n\times\R\to\R^n$ be the projection $\pi(x,x_{n+1}):=x$, let $\sph^n(\epsilon)$ be the $n$-sphere of $\R^{n+1}$ centered at the origin $O_1=(O,0)$ of radius $\epsilon$ and let $X_1:=\pi^{-1}(X)\cap\sph^n(\epsilon)$. Observe that $\sph^n(\epsilon)$ intersects transversally $\Reg(\pi^{-1}(X))=\Reg(X)\times\R$ in $\R^{n+1}$ since $\partial\overline{\B}_n(\epsilon)$ intersects transversally $\Reg(X)$ in $\R^n$. Hence, $X_1\subset\R^{n+1}$ is a compact algebraic set of dimension $d$ such that $\Sing(X_1)=\{N,S\}$, where $N:=(0,\ldots,0,\epsilon)$ and $S:=(0,\ldots,0,-\epsilon)$. Consider the translation $\tau:\R^{n+1}\to\R^{n+1}$ associated to the vector $-N$, i.e., $\tau(x,x_{n+1}):=(x,x_{n+1}-\epsilon)$. Observe that $X_2:=\tau(X_1)\subset\R^{n+1}$ is also a compact algebraic set of dimension~$d$ such that $\Sing(X_2)=\{O_1,S-N\}$. In addition, if $U:=X\cap\B_n(\epsilon)$ and $U_2:=X_2\cap\{x_{n+1}>-\epsilon\}$, then the map $\phi_2:U\to U_2$, $x\mapsto(x,\sqrt{\epsilon^2-|x|_n^2}-\epsilon)$ is a semialgebraic homeomorphism from the open semialgebraic neighborhood $U$ of $O$ in $X$ onto {the} open semialgebraic neighborhood $U_2$ of $O_1$ in~$X_2$. Apply Theorem~\ref{thm:main-2} and Remark \ref{rem:degen-2} with $P:=\{O_1\}$ to $X_2$, obtaining a $\Q$-determined $\Q$-algebraic subset $X'_2$ of $\R^{2d+4}$ and a semialgebraic homeomorphism $\eta_2:X_2\to X'_2$ such that $\eta_2(\Reg(X_2))=\Reg(X'_2)$, the restriction of $\eta_2$ from $\Reg(X_2)$ to $\Reg(X'_2)$ is a Nash diffeomorphism and $\eta_2(O_1)=(1,0)\in\R^{2d+4}=\R\times\R^{2d+3}$. Let $\tau':\R^{2d+4}\to\R^{2d+4}$ be the translation associated to the vector $-(1,0)$. Since $(1,0)\in\Q^{2d+4}$, $\tau'$ is a $\Q$-biregular isomorphism. Thus, the set $X':=\tau'(X'_2)$ is also a $\Q$-determined $\Q$-algebraic subset of $\R^{2d+4}$ such that the origin $O'$ of $\R^{2d+4}$ is an isolated point of $\Sing(X')$. Define the open semialgebraic neighborhood $U'$ of $O'$ in $X'$ and the map $\phi:U\to U'$ by $U':=\tau'(\eta_2(U_2))$ and $\phi(x):=\tau'(\eta_2(\phi_2(x)))$. By construction, $\phi$ has all the required properties.


\addtocontents{toc}{\protect\setcounter{tocdepth}{1}}

\appendix

\section{Proofs of some preparatory lemmas}\label{appendix-A}
We collect here the proofs of some results presented above.

\begin{proof}[Proof of Lemma \ref{X}]
By \cite[Prop.3.3.10]{BCR}, $a$ belongs to only one irreducible component $X'$ of $X$ such that $\dim(X')=e$ and $a\in\Reg(X')$. Let $X''$ be the union of all the irreducible components of $X$ different from $X'$. Choose $p\in\R[\x]$ such that $\ZZ_\R(p)=X'$. By \cite[Prop.2.2.11]{akbking:tras}, there exist polynomials $h_0,h_1,\ldots,h_{n-e}\in\R[\x]$ and an open neighborhood $U$ of $a$ in $\R^n$ such that $h_0(x)\neq0$ and $h_0(x)p(x)=\sum_{i=1}^{n-e}h_i(x)f_i(x)$ for all $x\in U$. It follows that $X'\cap U\subset\ZZ_\R(f_1,\ldots,f_{n-e})\cap U\subset\ZZ_\R(p)\cap U=X'\cap U$, so $X'\cap U=\ZZ_\R(f_1,\ldots,f_{n-e})\cap U$. Shrinking $U$ around $a$ if necessary, we can assume that $X''\cap U=\varnothing$ so $X\cap U=X'\cap U=\ZZ_\R(f_1,\ldots,f_{n-e})\cap U$, as required.
\end{proof}

\begin{proof}[Proof of Lemma \ref{lem:Q-basic}]
$(\mr{i})$ The `if' implication is evident. Let us prove the `only if' implication. Working component by component, we can assume that  $m=1$ and $Y=\R$. Let $a\in X$, let $U_a$ be a $K$-Zariski open neighborhood of $a$ in $\R^n$ and let $p_a,q_a\in\Q[\x]$ such that $q_a(x)\neq0$ and $f(x)=\frac{p_a(x)}{q_a(x)}$ for each $x\in X\cap U_a$. Choose $u_a\in\Q[\x]$ such that $\ZZ_\R(u_a)=\R^n\setminus U_a$. Replacing $p_a$ with $p_au_a$, $q_a$ with $q_au_a$ and $U_a$ with $U_a\setminus\ZZ_\R(q_a)$, we can assume that $U_a=\R^n\setminus\ZZ_\R(q_a)$. Since the $\Q$-Zariski topology on $X$ is Noetherian (and thus compact), there exist finitely many points $a_1,\ldots,a_\ell$ of $X$ such that $X=\bigcup_{i=1}^\ell(X\setminus\ZZ_\R(q_{a_i}))$. Define $p,q\in \Q[\x]$ by $p:=\sum_{i=1}^\ell q_{a_i}p_{a_i}$ and $q:=\sum_{i=1}^\ell q_{a_i}^2$. Observe that $X\cap\ZZ_\R(q)=\varnothing$ and $f(x)=\frac{p(x)}{q(x)}$ for all~$x\in X$.

$(\mr{ii})$ By $(\mr{i})$, there exist $p_1,\ldots,p_m,q\in\Q[\x]$ and $\xi_1,\ldots,\xi_{\ell+1}\in\Q[\y]:=\Q[\y_1,\ldots,\y_m]$ such that $q(x)\neq0$ and $f(x)=\big(\frac{p_1(x)}{q(x)},\ldots,\frac{p_m(x)}{q(x)}\big)$ for all $x\in X$, and $\xi_{\ell+1}(y)\neq0$ and $g(y)=\big(\frac{\xi_1(y)}{\xi_{\ell+1}(y)},\ldots,\frac{\xi_\ell(y)}{\xi_{\ell+1}(y)}\big)$ for all $y\in Y$. For each $k\in\{1,\ldots,\ell+1\}$, write $\xi_k$ as follows: $\xi_k=\sum_{j=0}^{c_k}\xi_{k,j}$, where $c_k$ is the degree of $\xi_k$ and each $\xi_{k,j}$ is a homogeneous polynomial in $\Q[\y]$ of degree $j$. Let $c:=\max\{c_1,\ldots,c_{\ell+1}\}$ and, for each $k\in\{1,\ldots,\ell+1\}$, let $\eta_k\in\Q[\x]$ be the polynomial $\eta_k:=\sum_{j=0}^{c_k}q^{c-j}\xi_{k,j}(p_1,\ldots,p_m)$, so $\eta_k(x):=q(x)^c\xi_k(f(x))$ for all $x\in X$. It follows that $g\circ f\in\reg^\Q(X,Z)$ since $\eta_{\ell+1}(x)\neq0$ and $(g\circ f)(x)=\big(\frac{\eta_1(x)}{\eta_{\ell+1}(x)},\ldots,\frac{\eta_\ell(x)}{\eta_{\ell+1}(x)}\big)$ for all $x\in X$.

$(\mr{iii})$ Consider again $p_1,\ldots,p_m,q\in\Q[\x]$ such that $q(x)\neq0$ and $f(x)=\big(\frac{p_1(x)}{q(x)},\ldots,\frac{p_m(x)}{q(x)}\big)$ for all $x\in X$. Consider also $p\in\Q[\x]$ such that $\ZZ_\R(p)=X$. Since $\ZZ_\R(p^2+q^2)=\varnothing$, we can define $F\in\reg^\Q(\R^n,\R^m)$ by $F(x):=\big(\frac{p_1(x)q(x)}{p(x)^2+q(x)^2},\ldots,\frac{p_m(x)q(x)}{p(x)^2+q(x)^2}\big)$ so that $F(x)=f(x)$ for all $x\in X$.

$(\mr{iv})$ This item follows immediately from Definition \ref{def:Q-reg-function}.

$(\mr{v})$ Let $q\in\Q[\y]$ be such that $W=Y\cap\ZZ_\R(q)$ and let $h:X\to Y$ be the composition $h(x):=q(f(x))$. By $(\mr{i})$ and $(\mr{ii})$, there exist $h_1,\ldots,h_{m+1}\in\Q[\x]$ such that $h_{m+1}(x)\neq0$ and $h=\big(\frac{h_1(x)}{h_{m+1}(x)},\ldots,\frac{h_m(x)}{h_{m+1}(x)}\big)$ for all $x\in X$. Consider again $p\in\Q[\x]$ such that $\ZZ_\R(p)=X$. We have $f^{-1}(W)=\ZZ_\R(p,h_1,\ldots,h_m)$, so $f^{-1}(W)\subset\R^n$ is $\Q$-algebraic.
\end{proof}

\begin{proof}[Proof of Lemma \ref{lem:diff-diff}]
$(\mr{i})$ Let $p_1,\ldots,p_d$ be affine functions on $\R^n$ vanishing at $a$ such that $\nabla f_1(a)$, $\ldots,\nabla f_{n-d}(a)$, $\nabla p_1(a),\ldots,\nabla p_d(a)$ are linearly independent. By the Inverse Function Theorem, there exist a sufficiently small $\epsilon>0$ and an open neighborhood $U'$ of $a$ in $U$ such that the map $\psi:U'\to\B_n(\epsilon)$, $x\mapsto(f_1(x),\ldots,f_{n-d}(x),p_1(x),\ldots,p_d(x))$ is a well-defined smooth diffeomorphism. Define $F\in\cinfty(\B_n(\epsilon))$ by $F(x):=f(\psi^{-1}(x))$. Observe that $\psi^{-1}(\B_n(\epsilon)\cap(\{0\}\times\R^d))=V\cap U'$ so $F$ vanishes on $\{0\}\times\B_d(\epsilon)$. We have
$$\textstyle
F(y)=\int_0^1\frac{d}{dt}F(ty',y'')\,dt=\sum_{i=1}^{n-d}y_i\int_0^1\frac{\partial F}{\partial y_i}(ty',y'')\,dt,
$$
where $y':=(y_1,\ldots,y_{n-d})$, $y'':=(y_{n-d+1},\ldots,y_n)$ and $(y',y'')\in\B_n(\epsilon)$. For all $i\in\{1,\ldots,n-d\}$, define $u_i\in\cinfty(U')$ by $u_i(x):=\int_0^1\frac{\partial F}{\partial y_i}(t\psi'(x),\psi''(x))\,dt$, where $\psi'(x):=(f_1(x),\ldots,f_{n-d}(x))$ and $\psi''(x):=(p_1(x),\ldots,p_d(x))$. It follows that $f=\sum_{i=1}^{n-d}u_if_i$ on $U'$.

$(\mr{ii})$ Suppose that $e<d$. If $e=d$, the proof is similar. We start as above. Choose $d-e$ affine functions $p_1,\ldots,p_{d-e}$ on $\R^n$ vanishing at $a$ such that $\nabla f_1(a),\ldots,\nabla f_{n-d}(a)$, $\nabla g_1(a),\ldots,\nabla g_e(a)$, $\nabla p_1(a),\ldots,\nabla p_{d-e}(a)$ are linearly independent. Therefore, there exist an open neighborhood $U'$ of $a$ in $U$ and a sufficiently small $\epsilon>0$ such that the map $\psi:U'\to\B_n(\epsilon)$, $x\mapsto(f_1(x),\ldots,f_{n-d}(x),g_1(x),\ldots,g_e(x),p_1(x),\ldots,p_{d-e}(x))$ is a well-defined smooth diffeomorphism. Let $x':=(x_1,\ldots,x_{n-d})$, let $x'':=(x_{n-d+1},\ldots,x_n)$ and let $x=(x',x'')$ be the coordinates of $\R^n$. Replacing the functions $f_i$ with $f_i\circ\psi^{-1}$ and $g_j$ with $g_j\circ\psi^{-1}$, we can assume that $a$ is the origin of~$\R^n$, $U=\B_n(\epsilon)$, $f_i(x)=x_i$ for all $i\in\{1,\ldots,n-d\}$ and $g_j(x)=x_{n-d+j}$ for all $j\in\{1,\ldots,e\}$, so $g(x)=\prod_{j=1}^ex_{n-d+j}$, $V=\{0\}\times\B_d(\epsilon)$ and $W=\bigcup_{j=1}^e\{x\in V:x_{n-d+j}=0\}$. By the preceding part of the proof, we have that $f(x',x'')-f(0,x'')=\sum_{i=1}^{n-d}u_i(x)x_i$ for all $x=(x',x'')\in\B_n(\epsilon)$ and for some $u_1,\ldots,u_{n-d}\in\cinfty(\B_n(\epsilon))$. Define $F\in\cinfty(\B_d(\epsilon))$ by $F(x''):=f(0,x'')$. Identify $V=\{0\}\times\B_d(\epsilon)$ with $\B_d(\epsilon)$ and $W$ with $\bigcup_{j=1}^e\{x''\in\B_d(\epsilon):x_{n-d+j}=0\}$. It remains to prove that $F(x'')=u_0\prod_{j=1}^ex_{n-d+j}$ for some $u_0\in\cinfty(\B_d(\epsilon))$. Since $F$ vanishes on~$W$ (and therefore on $\{x''\in\B_d(\epsilon):x_{n-d+1}=0\}$), using again the preceding part of the proof, we have that $F(x'')=u_{0,1}(x'')x_{n-d+1}$ on $\B_d(\epsilon)$ for some $u_{0,1}\in\cinfty(\B_d(\epsilon))$. If $e=1$, we are done by setting $u_0:=u_{0,1}$. If $e>1$, then $F$ also vanishes on $\{x''\in\B_d(\epsilon):x_{n-d+2}=0\}$, so $u_{0,1}$ vanishes on $\{x''\in\B_d(\epsilon):x_{n-d+2}=0,x_{n-d+1}\neq0\}$. By density, $u_{0,1}$ vanishes on the whole $\{x''\in\B_d(\epsilon):x_{n-d+2}=0\}$, so $u_{0,1}(x'')=u_{0,2}(x'')x_{n-d+2}$ on $\B_d(\epsilon)$ for some $u_{0,2}\in\cinfty(\B_d(\epsilon))$. It follows that $F(x'')=u_{0,2}(x'')x_{n-d+1}x_{n-d+2}$ on $\B_d(\epsilon)$. Proceeding inductively in this way, we obtain a function $u_{0,e}\in\cinfty(\B_d(\epsilon))$ such that $F(x'')=u_{0,e}(x'')\prod_{j=1}^ex_{n-d+j}$ on $\B_d(\epsilon)$. It is now sufficient to set $u_0:=u_{0,e}$, completing the proof.
\end{proof}

\section{Proof of Theorem \ref{thm:1.7}} \label{appendix-B}
Here we prove Theorem \ref{thm:1.7} by showing how to adapt the proofs of \cite[Thm.1.7\;\&\;Prop.8.2]{BaFeRu2014} to our smooth setting. Actually, we also need to extend the proofs of \cite[Prop.6.2\;\&\;8.1]{BaFeRu2014}.

\emph{Fix Nash manifolds $N\subset\R^n$ and $M\subset\R^m$, and Nash sets $X\subset N$ and $Y\subset M$.}

Let $f:X\to Y$ be a map. We say that $f$ is a \emph{$\cinfty$~map} if there exist an open neighborhood $U$ of $X$ in $N$ and a map $F:U\to\R^m$, smooth in the usual sense, such that $F(x)=f(x)$ for all $x\in X$. Let $\cinfty(X,Y)$ be the set of smooth maps from $X$ to $Y$, and let $\Nn(X,Y)$ be the set of Nash maps from $X$ to~$Y$. Evidently, $\Nn(X,Y)\subset\cinfty(X,Y)$.

Let $X_1,\ldots,X_s$ be the Nash irreducible components of $X$. According to \cite[Def.1.5]{BaFeRu2014}, we say that the map $f:X\to Y$ is \emph{$\cc$-Nash} if the restriction $f|_{X_i}:X_i\to Y$ is a Nash map for every $i\in\{1,\ldots,s\}$. Let $^\cc\Nn(X,Y)$ be the set of $\cc$-Nash maps from $X$ to $Y$. Let us extend this definition to $\cinfty$ maps. We say that $f:X\to Y$ is a \emph{$^\cc\cinfty$ map} if the restriction $f|_{X_i}:X_i\to Y$ is a $\cinfty$~map for every $i\in\{1,\ldots,s\}$. Let $^\cc\cinfty(X,Y)$ the set of $^\cc\cinfty$ maps from $X$ to $Y$. Evidently, $\Nn(X,Y)\subset{^\cc}\Nn(X,Y)$ and $\cinfty(X,Y)\subset{^\cc}\cinfty(X,Y)$. If $Y=\R$, we simplify notations by setting $\cinfty(X):=\cinfty(X,\R)$, ${^\cc}\cinfty(X):={^\cc}\cinfty(X,\R)$, $\Nn(X)=\Nn(X,\R)$ and ${^\cc}\Nn(X):={^\cc}\Nn(X,\R)$.

Let us define a topology on $\cinfty(X,Y)$ and one on $^\cc\cinfty(X,Y)$ following \cite[Sects.2.C\;\&\;2D]{BaFeRu2014}. First, we consider the case $Y=\R^m$. Endow $\cinfty(N,\R^m)$ with the usual weak $\cinfty$ topology, see \cite[p.36]{hirsch:difftop}. This topology makes $\cinfty(N,\R^m)$ a topological real vector space with the usual pointwise defined addition and multi\-plication by real scalars. Consider the restriction map $\rho:\cinfty(N,\R^m)\to\cinfty(X,\R^m)$, i.e., $\rho(F):=F|_X$. Using smooth partition of unity, it is immediate to verify that $\rho$ is surjective. Endow $\cinfty(X,\R^m)$ with the quotient topology induced by $\rho$. An important property of this quotient topology is that $\rho$ is an open map. Indeed, if $\mc{U}\subset\cinfty(N,\R^m)$ is open and $\II:=\{F\in\cinfty(N,\R^m):\rho(F)=0\}$, then $\rho^{-1}(\rho(\mc{U}))=\bigcup_{F\in\II}(\mc{U}+F)$ is open as well, {since} the translations of $\cinfty(N,\R^m)$ are homeomorphisms. Identify $\cinfty(X,Y)$ with the subset of $\cinfty(X,\R^m)$ of all $\cinfty$ maps $f:X\to\R^m$ such that $f(X)\subset Y$, and endow $\cinfty(X,Y)$ with the relative topology induced by $\cinfty(X,\R^m)$. Let us define a topology on $^\cc\cinfty(X,Y)$ as well. Consider the topological product $\cinfty(X_1,Y)\times\cdots\times\cinfty(X_s,Y)$ and the multiple restriction map $j:{^\cc}\cinfty(X,Y)\to\cinfty(X_1,Y)\times\cdots\times\cinfty(X_s,Y)$ defined by $j(f):=(f|_{X_1},\ldots,f|_{X_s})$. Endow $^\cc\cinfty(X,Y)$ with the topology making $j$ a homeomorphism onto its image.

Let $\gamma:\cinfty(X,Y)\hookrightarrow{^\cc}\cinfty(X,Y)$ be the inclusion map. By the universal property of the quotient topology, we know that each restriction map $\cinfty(X,\R^m)\to\cinfty(X_i,\R^m)$ is continuous. This implies that $\gamma:\cinfty(X,Y)\hookrightarrow{^\cc}\cinfty(X,Y)$ is continuous as well. Let $J:\cinfty(X,Y)\to\cinfty(X_1,Y)\times\cdots\times\cinfty(X_s,Y)$ be the composition map $j\circ\gamma$, which is continuous.

Given $s\in\N^*$, we denote $\mc{P}(s)$ the power set of $\{1,\ldots,s\}$.

The next result is a smooht variant of \cite[Props.6.2\;\&\;8.1]{BaFeRu2014}.

\begin{prop}\label{prop:6.2}
If $X\subset N$ is a Nash set with monomial singularities, then there exists a continuous linear map $\theta:\cinfty(X,\R^m)\to\cinfty(N,\R^m)$ such that $\theta$ is an extension map, i.e, $\theta(f)|_X=f$ for all $f\in\cinfty(X,\R^m)$. Moreover, the multiple restriction map $J:\cinfty(X,Y)\to\cinfty(X_1,Y)\times\cdots\times\cinfty(X_s,Y)$ is a homeomorphism onto its image.
\end{prop}
\begin{proof}
As in the proof of \cite[Prop.6.2]{BaFeRu2014}, it suffices to consider the case $\R^m=\R=Y$ and to prove the existence of a continuous linear extension map $^\cc\theta:{^\cc}\cinfty(X,\R)\to\cinfty(N,\R)$. This implies at once that $\theta:={^\cc\theta}\circ\gamma$ is the required extension map and $\gamma$ is a homeomorphism, so $J$ is a homeomorphism onto its image. The problem of constructing $^\cc\theta$ is local in nature, {since} $N$ admits smooth partitions of unity subordinate to each of its open covers. This fact and Definition~\ref{def:Nmsing} reduce the problem to the case in which $N=\R^n$ and $X=L_1\cup\ldots\cup L_s$ is a union of coordinate linear subspaces of $\R^n$. In this situation, the proof of \cite[Prop.4.C.1]{BaFeRu2014} gives an explicit formula for $^\cc\theta$. For every $I\in\mc{P}(s)\setminus\{\varnothing\}$, let $L_I:=\bigcap_{i\in I}L_i$, let $i_I:L_I\hookrightarrow X$ be the inclusion map and let $\pi_I:\R^n\to L_I$ be the orthogonal projection of $\R^n$ onto $L_I$. Set $^\cc\theta(f):=-\sum_{I\in\mc{P}(s)\setminus\{\varnothing\}}(-1)^{|I|}(f\circ i_I\circ\pi_I)$, where $|I|$ is the cardinality of $I$. It is immediate to verify that $^\cc\theta(f)|_{L_i}=f|_{L_i}$ for all $i\in\{1,\ldots,s\}$, so $^\cc\theta(f)|_X=f$. Since the composition operation is continuous in the weak $\cinfty$ topology, the map $^\cc\theta$ is continuous. 
\end{proof}

\begin{remark}\label{rem:1.6}
In \cite[Thm.1.6]{BaFeRu2014}, the authors prove the following remarkable fact: if $X\subset N$ is a Nash set with monomial singularities, then $^\cc\Nn(X)=\Nn(X)$. $\sqbullet$
\end{remark}

Let $U$ be an open semialgebraic subset of $N$ and let $S$ be a subset of $U$. Denote {by} $\cinfty(U)$ the ring of smooth functions on $U$, $\Nn(U)$ the ring of Nash functions on $U$, and define the vanishing ideals $\II^\infty_U(S):=\{f\in\cinfty(U):f|_S\equiv0\}$ and $\II^\Nn_U(S):=\{f\in\Nn(U):f|_S\equiv0\}$.

The next result is a smooth variant of \cite[Prop.8.2]{BaFeRu2014}.

\begin{prop}\label{prop:8.2}
Let $L\subset N$ be a Nash set with monomial singularities and let $F:N\to M$ be a smooth map. Then every Nash map $h:L\to M$ which is sufficiently $\cinfty$ close to $F|_L$ has a Nash extension $H:N\to M$ which is arbitrarily $\cinfty$ close to $F$. More precisely, for every neighborhood $\mc{U}$ of $F$ in $\cinfty(N,M)$, there exists a neighborhood $\mc{V}$ of $F|_L$ in $\cinfty(L,M)$ with the following property: for every Nash map $h\in\Nn(L,M)\cap\mc{V}$, there exists a Nash map $H\in\Nn(N,M)\cap\mc{U}$ such that $H|_L=h$. 
\end{prop}
\begin{proof}
Since $M$ has a Nash tubular neighborhood in $\R^m$ (see \cite[Cor.8.9.5]{BCR} or \cite[\S I.3]{Sh}), it suffices to consider the case in which $M=\R^m$. Working component by component, we can further assume that $M=\R$. We now follow the strategy of the proof of \cite[Prop.7.6]{BaFeRu2014}. Let $\mc{U}$ be a neighborhood of $F$ in $\cinfty(N)$. The openness of the restriction map $\rho:\cinfty(N)\to\cinfty(L)$ implies that $\rho(\mc{U})$ is a neighborhood of $F|_L$ in $\cinfty(L)$. Suppose that $h\in\rho(\mc{U})$ and choose $G\in\mc{U}$ such that $\rho(G)=h$, i.e., $G|_L=h$. It remains to show that there exists a Nash function $H:N\to\R$ arbitrarily $\cinfty$ close to $F$ (hence, we may assume that $H\in\mc{U}$) such that $H|_L=G|_L$. By \cite[Thm.1.4]{BaFeRu2014}, there exist a finite family of open semialgebraic subsets $U_1,\ldots,U_\ell$ of $N$ and Nash diffeomorphisms $\{u_i:U_i\to\R^d\}_{i=1}^\ell$, where $d=\dim(N)$, such that $L\subset U_1\cup\ldots\cup U_\ell$ and, for each $i\in\{1,\ldots,\ell\}$, $u_i(L\cap U_i)$ is a union of coordinate linear subspaces of $\R^d$. It follows that $\II^\infty_{\R^d}(u_i(L\cap U_i))\subset\II^\Nn_{\R^d}(u_i(L\cap U_i))\cinfty(\R^d)$ for each $i$. The latter inclusion can be proven exactly as in \cite[Prop.7.3\;\&\;Rmk.7.4]{BaFeRu2014} (however, here the proof is slightly simpler). As an immediate consequence, we have $\II^\infty_{U_i}(L\cap U_i)\subset\II^\Nn_{U_i}(L\cap U_i)\cinfty(U_i)$. Since $L\subset N$ is coherent, it holds $\II^\Nn_{U_i}(L\cap U_i)=\II^\Nn_N(L)\Nn(U_i)$ (see equation (2.2) and Lemma 5.1 of \cite{BaFeRu2014}). It follows that $\II^\infty_{U_i}(L\cap U_i)\subset\II^\Nn_N(L)\cinfty(U_i)$ for all $i\in\{1,\ldots,\ell\}$. Making use of a smooth partition of unity subordinate to the open cover $\{U_1,\ldots,U_\ell,N\setminus L\}$ of $N$, we obtain at once that $\II^\infty_N(L)\subset\II^\Nn_N(L)\cinfty(N)$. Since $G|_L=h\in\Nn(L)$, by \cite[Thm.8.9.12]{BCR}, there exists $\widetilde{G}\in\Nn(N)$ such that $\widetilde{G}|_L=G|_L$, i.e., $G-\widetilde{G}\in\II^\infty_N(L)$. Thus, there exist $e\in\N^*$, functions $\{f_j\}_{j=1}^e$ in $\II^\Nn_N(L)$ and functions $\{\psi_j\}_{j=1}^e$ in $\cinfty(N)$ such that $G-\widetilde{G}=\sum_{j=1}^e\psi_jf_j$ on $N$. Thanks to the Weierstrass approximation theorem, for each $j\in\{1,\ldots,e\}$, there exists a polynomial function $\widetilde{\psi}_j:N\to\R$ arbitrarily $\cinfty$ close to $\psi_j$ on $N$. This proves that the Nash function $H\in\Nn(N)$ defined by $H:=\widetilde{G}+\sum_{j=1}\widetilde{\psi}_jf_j$ is arbitrarily $\cinfty$ close to $G$ and $H|_L=G|_L$.
\end{proof}

Let us restate and prove Theorem \ref{thm:1.7}.

\vspace{.5em}

\noindent \textbf{Theorem \ref{thm:1.7}.} \textit{Let $N\subset\R^n$ and $M\subset\R^m$ be Nash manifolds, let $X\subset N$ and $Y\subset M$ be Nash monomial crossings and let $f:N\to M$ be a smooth map such that $f(X)\subset Y$ and the restriction of $f$ from $X$ to $Y$ preserves Nash irreducible components, i.e., the image of each Nash irreducible component of $X$ under $f$ is contained in some Nash irreducible component of~$Y$. Then there exists a Nash map $g:N\to M$ arbitrarily $\cinfty$ close to $f$ such that $g(X)\subset Y$ and the restriction of $g$ from $X$ to $Y$ preserves Nash irreducible components.}
\begin{proof}
Let us adapt the proof of \cite[Thm.1.7]{BaFeRu2014} to the present situation. Let $X_1,\ldots,X_s$ be the Nash irreducible components of $X$ and let $Y_1,\ldots,Y_t$ be the Nash irreducible components of~$Y$. By hypothesis, all $X_i\subset N$ and $Y_j\subset M$ are Nash manifolds, and there exists a function $\kappa:\{1,\ldots,s\}\to\{1,\ldots,t\}$ such that $f(X_i)\subset Y_{\kappa(i)}$ for all $i\in\{1,\ldots,s\}$. For each $J\in\mc{P}(t)$ and for each $p\in\N^*$, we set $Y_J:=\bigcap_{j\in J}Y_j$, $X_J:=\bigcap_{i\in\kappa^{-1}(J)}X_i$, $\mc{P}(t,p):=\{J\in\mc{P}(t):|J|=p\}$, $X^{\sss(p)}:=\bigcup_{J\in\mc{P}(t,p)}X_J$ and $Y^{\sss(p)}:=\bigcup_{J\in\mc{P}(t,p)}Y_J$. 
By construction, we have that $f(X_J)\subset Y_J$ and $f(X^{\sss(p)})\subset Y^{\sss(p)}$ for all $J\in\mc{P}(t)$ and for all $p\in\N^*$, and $X^{\sss(1)}=X$ and $Y^{\sss(1)}=Y$. By \cite[Prop.8.3]{BaFeRu2014}, we know that all $X_J\subset\R^n$ and $Y_J\subset\R^m$ are Nash manifolds, and all $X^{\sss(p)}\subset N$ and $Y^{\sss(p)}\subset M$ are Nash monomial crossings. Observe that each of the sets $X_J$, $X^{\sss(p)}$, $Y_J$ and $Y^{\sss(p)}$ may be empty. Let $r:=\max\{p\in\N^*:X^{\sss(p)}\neq\varnothing\}$. For every $p\in\{1,\ldots,r\}$, we have that both $X^{\sss(p)}$ and $Y^{\sss(p)}$ are non-empty, {since} $X^{\sss(p)}\supset X^{(r)}\neq\varnothing$ and $Y^{\sss(p)}\supset f(X^{\sss(p)})\neq\varnothing$. Thus, we can define the $\cinfty$ map $f_p:X^{\sss(p)}\to Y^{\sss(p)}$ as the restriction of $f$ from $X^{\sss(p)}$ to $Y^{\sss(p)}$. For each $J\in\mc{P}(t)$ such that $X_J\neq\varnothing$ (and hence $Y_J\supset  f(X_J)\neq\varnothing$), we define the $\cinfty$ map $f_J:X_J\to Y_J$ as the restriction of $f$ from $X_J$ to $Y_J$.
	
Let us prove the following claim: For every $p\in\{1,\ldots,r\}$, there exists a Nash map $g_p:X^{\sss(p)}\to Y^{\sss(p)}$ arbitrarily $\cinfty$ close to $f_p$ such that $g_p(X_J)\subset Y_J$ for all $J\in\mc{P}^*(t,p)$, where $\mc{P}^*(t,p):=\bigcup_{\ell=p}^r\mc{P}(t,\ell)$.
	
Let us proceed by induction on $p=r,r-1,\ldots,1$. Suppose that $p=r$. Observe that $X^{\sss(r)}$ is the disjoint union of the $X_J$'s with $J$ in $\mc{P}(t,r)$. Otherwise, there would exist $J_1,J_2\in\mc{P}(t,r)$ with $J_1\neq J_2$ (and thus $|J_1\cup J_2|>r$) such that $\varnothing\neq X_{J_1}\cap X_{J_2}=X_{J_1\cup J_2}$, contradicting the maximality of $r$. Given any $J\in\mc{P}(t,r)$, we have that $f_r(X_J)=f(X_J)\subset Y_J$. If $X_J\neq\varnothing$ (and thus $Y_J\neq\varnothing$), then the Weierstrass approximation theorem and the existence of Nash tubular neighborhoods of $Y_J$ in $\R^m$ imply that there exists a Nash map $g_J:X_J\to Y_J$ arbitrarily $\cinfty$ close to $f_J$. Since $X^{\sss(r)}$ is the disjoint union of the $X_J$'s with $J\in\mc{P}(t,r)$, the Nash map $g_r:X^{\sss(r)}\to Y^{\sss(r)}$, defined as $g_r(x):=g_J(x)$ if $x\in X_J$ for some (unique) $J\in\mc{P}(t,r)$, is arbitrarily $\cinfty$ close to $f_r$ and $g_r(X_J)\subset Y_J$ for all $J\in\mc{P}(t,r)=\mc{P}^*(t,r)$. 

Suppose now the claim {holds} for some $p\in\{1,\ldots,r\}$. We {may} assume that $p\neq1$, otherwise we are done. Let $g_p:X^{\sss(p)}\to Y^{\sss(p)}$ be a Nash map arbitrarily $\cinfty$ close to $f_p$ such that $g_p(X_J)\subset Y_J$ for all $J\in\mc{P}^*(t,p)$. Let $K\in\mc{P}(t,p-1)$. We have that $g_p(X_K\cap X^{\sss(p)})\subset Y_K$. Indeed, $X_K\cap X^{\sss(p)}=\bigcup_{J\in\mc{P}(t,p)}(X_K\cap X_J)=\bigcup_{J\in\mc{P}(t,p)}X_{K\cup J}$ and thus $g_p(X_K\cap X^{\sss(p)})=\bigcup_{J\in\mc{P}(t,p)}g_p(X_{K\cup J})\subset \bigcup_{J\in\mc{P}(t,p)}Y_{K\cup J}\subset Y_K$. Here we used the fact that, for all $J\in\mc{P}(t,p)$, it holds $K\cup J\in\mc{P}^*(t,p)$ so $g_p(X_{K\cup J})\subset Y_{K\cup J}$. If $X_K\cap X^{\sss(p)}\neq\varnothing$, then we define the $\cinfty$ maps $f_{p,K}:X_K\cap X^{\sss(p)}\to Y_K$ and $g_{p,K}:X_K\cap X^{\sss(p)}\to Y_K$ as the restrictions of $f_p$ and $g_p$ from $X_K\cap X^{\sss(p)}$ to $Y_K$, respectively. Recall that $X_K\subset\R^n$ and $Y_K\subset\R^m$ are Nash manifolds. In addition, $X_K\cap X^{\sss(p)}\subset X_K$ is a Nash monomial crossings by \cite[Prop.8.3]{BaFeRu2014}. If we choose $g_p$ sufficiently $\cinfty$ close to $f_p$, then we {may} assume that, for each $K\in\mc{P}(t,p-1)$,  $g_{p,K}$ is arbitrarily $\cinfty$ close to $f_{p,K}=f_K|_{X_K\cap X^{\sss(p)}}$. By Proposition \ref{prop:8.2}, there exists a Nash map $\widetilde{g}_{p,K}:X_K\to Y_K$ arbitrarily $\cinfty$ close to $f_K$ such that $\widetilde{g}_{p,K}|_{X_K\cap X^{\sss(p)}}=g_{p,K}=g_p|_{X_K\cap X^{\sss(p)}}$. Consider $K,K'\in\mc{P}(t,p-1)$ such that $K\neq K'$ and $X_K\cap X_{K'}\neq\varnothing$, and choose $x\in X_K\cap X_{K'}$. Since $|K\cup K'|>p-1$, we deduce that $X_K\cap X_{K'}=X_{K\cup K'}\subset X^{\sss(p)}$ so $x\in X_K\cap X^{\sss(p)}$ and $x\in X_{K'}\cap X^{\sss(p)}$. It follows that $\widetilde{g}_{p,K}(x)=g_p(x)=\widetilde{g}_{p,K'}(x)$. This proves that the map $g_{p-1}:X^{\sss(p-1)}\to Y^{\sss(p-1)}$, defined by $g_{p-1}(x):=\widetilde{g}_{p,K}(x)$ if $x\in X_K$ for some $K\in\mc{P}(t,p-1)$, is a well-defined $\cc$-Nash map, which is also a Nash map by Remark~\ref{rem:1.6}. Observe that $g_{p-1}(X_K)\subset Y_K$ for all $K\in\mc{P}(t,p-1)$ and $g_{p-1}$ is an extension of $g_p$. Indeed, it coincides with $g_p$ on $\bigcup_{K\in\mc{P}(t,p-1)}(X_K\cap X^{\sss(p)})=\big(\bigcup_{K\in\mc{P}(t,p-1)}X_K\big)\cap X^{\sss(p)}=X^{\sss(p-1)}\cap X^{\sss(p)}=X^{\sss(p)}$. As a consequence, $g_{p-1}(X_K)\subset Y_K$ for all $K\in\mc{P}^*(t,p-1)$. Since $g_{p-1}|_{X_K}$ is arbitrarily $\cinfty$ close to $f_K=f|_{X_K}$ for all $K\in\mc{P}(t,p-1)$, by Proposition \ref{prop:6.2}, it follows that $g_{p-1}$ is arbitrarily $\cinfty$ close to $f_{p-1}$. This proves the preceding claim. In particular, we proved the existence of a Nash map $g_1:X\to Y$ arbitrarily $\cinfty$ close to the restriction of $f$ from $X$ to $Y$. Finally, using again Proposition \ref{prop:8.2}, we obtain a Nash map $g:N\to M$ arbitrarily $\cinfty$ close to $f$ such that $g(x)=g_1(x)$ for all $x\in X$.
\end{proof}

\section{Manifold nature of real algebraic sets `defined over $\Q$'} \label{appendix-C} 
Let $L:=\R$ or $\C$, and let $X$ be an algebraic subset of $L^n$. Suppose we want to give a formal meaning to the assertion: the algebraic set $X\subset\R^n$ is `defined over $\Q$'. The weakest definition seems to be the following: $X$ is `defined over $\Q$' if it can be described in $L^n$ by polynomial equations with rational coefficients, i.e., $X$ is a $\Q$-algebraic subset of $L^n$ in the sense of Definition~\ref{K-algebraic-set}. A stronger definition is the following: $X$ is `defined over $\Q$' if $\II_L(X)=\II_\Q(X)L[\x]$.

As we have seen in Remarks \ref{rem210}$(\mr{ii})$, when $L=\C$, the preceding two definitions coincide. In fact, by Hilbert's Nullstellensatz, every $\Q$-algebraic subset $X$ of $\C^n$ satisfies $\II_\C(X)=\II_\Q(X)\C[\x]$. In the real setting, the situation is completely different. In addition to the definition of real $\Q$-algebraic set given in Definition \ref{K-algebraic-set}, we can provide at least three (indeed, four) {distinct} definitions of real algebraic set `defined over $\Q$'.

The first is what we have already written above.

\begin{defn}[{\cite[Def.3, p.30]{togn:instmat}\;\&\;\cite[Def.3.4.10]{FG}}]
Let $X\subset\R^n$ be an algebraic set. We say that $X$ is \emph{defined over~$\Q$} if $\II_\R(X)=\II_\Q(X)\R[\x]$. $\sqbullet$
\end{defn}

Sections 3.4, 3.5 and 3.6 of \cite{FG} are devoted to the study of the latter concept.

Let $U$ and $S$ be two subsets of $\R^n$ such that $S\subset U$. Endow the set $\reg(U)$ of regular functions on $U$ with the usual ring structure induced by the pointwise addition and multiplication. Denote $\I^r_U(S)$ the ideal of $\reg(U)$ of all regular functions vanishing on $S$ and, given any $F\subset\R[\x]$, denote $F\reg(U)$ the ideal of $\reg(U)$ generated by $\{f|_U\}_{f\in F}$, where $f|_U:U\to\R$ is the regular (indeed, polynomial) function defined by $f|_U(x):=f(x)$ for all $x\in U$.

A second way to introduce the notion of real algebraic set `defined over $\Q$' is the following.

\begin{defn}
Let $X\subset\R^n$ be an algebraic set. We say that $X$ is \emph{weakly defined over $\Q$} if $\I^r_{\R^n}(X)=\I_\Q(X)\reg(\R^n)$. $\sqbullet$
\end{defn}

A third way is the notion of $\Q$-determined $\Q$-algebraic set given in Definition \ref{Q-determined}. A further way to define a real algebraic set `defined over $\Q$' is the following Diophantine notion:

\begin{defn}
Let $X\subset\R^n$ be an algebraic set and let $X(\Q):=X\cap\Q^n$. We say that $X$ is \emph{strongly defined over $\Q$} if $X(\Q)$ is Zariski dense in $X$. $\sqbullet$
\end{defn}

\begin{remark}
Let $X\subset\R^n$ be a $\Q$-irreducible $\Q$-algebraic set of dimension $d$. If there exists a $\Q$-biregular map from a non-empty $\Q$-Zariski open subset of $\R^d$ to a non-empty $\Q$-Zariski open subset of $X$, then obviously $X$ is strongly defined over $\Q$. $\sqbullet$
\end{remark}

A comparison result reads as follows.

\begin{prop}\label{prop:implications}
Let $X\subset\R^n$ be an algebraic set. Consider the following four conditions:
\begin{itemize}
 \item[$(\mr{i})$] $X$ is strongly defined over $\Q$.
 \item[$(\mr{ii})$] $X$ is defined over $\Q$.
 \item[$(\mr{iii})$] $X$ is weakly defined over $\Q$.
 \item[$(\mr{iv})$] $X$ is $\Q$-determined $\Q$-algebraic.
\end{itemize}
The implications $(\mr{i})\Longrightarrow(\mr{ii})\Longrightarrow(\mr{iii})\Longrightarrow(\mr{iv})$ {hold} true.
\end{prop}
\begin{proof}
$(\mr{i})\Longrightarrow(\mr{ii})$ Since the inclusion $\II_\Q(X)\R[\x]\subset\II_\R(X)$ is always true, we only have to show the converse inclusion. Let $f \in\II_\R(X)$ and let $\{u_p\}_{p\in P}$ be a basis of $\R$ as a $\Q$-vector space. As we have seen in the proof of Lemma \ref{K-difference-preparation}, there exist polynomials $f_p\in\Q[\x]$ such that only finitely many of the $f_p$ are non-zero and $f=\sum_{p\in P}f_pu_p$. If $x\in X(\Q)$, then $0=f(x)=\sum_{p\in P}f_p(x)u_p$. Since each value $f_p(x)$ belongs to $\Q$ and $\{u_p\}_{p\in P}$ is a $\Q$-vector basis of $\R$, we deduce that $f_p(x)=0$ for all $x\in X(\Q)$ and $p\in P$. Since $X(\Q)$ is Zariski dense in $X$ by hypothesis, it follows that $f_p\in\II_\Q(X)$ for all $p\in P$, so $f\in\II_\Q(X)\R[\x]$. This proves that $\II_\R(X)=\II_\Q(X)\R[\x]$.

$(\mr{ii})\Longrightarrow(\mr{iii})$ This is evident {since} $\II_\R(X)$ and $\II_\Q(X)\R[\x]$ generate $\II^r_{\R^n}(X)$ and $\II_\Q(X)\reg(\R^n)$ in $\reg(\R^n)$, respectively.

$(\mr{iii})\Longrightarrow(\mr{iv})$ Let $g_1,\ldots,g_s$ be generators of $\II_\Q(X)$ in $\Q[\x]$, let $f\in\II_\R(X)$ and let $a\in\Reg(X)$. As a regular (indeed, polynomial) function, $f$ belongs to $\II^r_{\R^n}(X)=\II_\Q(X)\reg(\R^n)$, so there exists $h_0,\ldots,h_s\in\R[\x]$ such that $\ZZ_\R(h_0)=\varnothing$ and $h_0(x)f(x)=\sum_{i=1}^sh_i(x)g_i(x)$ for all $x\in\R^n$. Thus, $h_0(a)\nabla f(a)=\sum_{i=1}^sh_i(a)\nabla g_i(a)$. It follows that  $X=\bigcap_{f\in\II_\R(X)}\ZZ_\R(f)=\ZZ_\R(g_1,\ldots,g_s)$ is a $\Q$-algebraic subset of $\R^n$ and $T_a(X)=\bigcap_{f\in\II_\R(X)}\nabla f(a)^\perp=\bigcap_{i=1}^s\nabla g_i(a)^\perp=T^\Q_a(X)$. By Proposition \ref{Y}, we know that $a\in\Reg^\Q(X)$. This proves that $X\subset\R^n$ is $\Q$-determined.
\end{proof}

\begin{remark}
Let $X\subset\R^n$ be an algebraic set and let $\overline{X(\Q)}$ be the $\Q$-Zariski closure of $X(\Q)$ in $\R^n$. If $X$ is strongly defined over $\Q$ then $\overline{X(\Q)}=X$, {since} $X$ is (defined over $\Q$ and thus) $\Q$-algebraic and $\overline{X(\Q)}$ contains the Zariski closure of $X(\Q)$ in $\R^n$, which is equal to $X$ by hypothesis, so $\overline{X(\Q)}=X$. The converse implication is also true: if $\overline{X(\Q)}=X$ then $X$ is strongly defined over $\Q$. To prove this, denote $X'$ the Zariski closure of $X(\Q)$ in $\R^n$. Since $X(\Q)\subset X'(\Q)$, we have that $X'\subset\R^n$ is strongly defined over $\Q$ so it is $\Q$-algebraic. It follows that $X=X'$ {since} $X=\overline{X(\Q)}\subset X'\subset X$. $\sqbullet$
\end{remark}

The next examples prove that Proposition \ref{prop:implications} is sharp, i.e., none of the converse implications $(\mr{i})\Longleftarrow(\mr{ii})\Longleftarrow(\mr{iii})\Longleftarrow(\mr{iv})$ {holds} in general.

\begin{examples}
$(\mr{i})$ Let $h$ be a natural number $\geq 3$, let $f_h:=\x_1^{2h}+\x_2^{2h}-2^h\in\Q[\x]:=\Q[\x_1,\x_2]$ and let $F_h$ be the Fermat curve $F_h:=\ZZ_\R(f_h)\subset\R^2$. Since $f_h$ is irreducible in $\R[\x]$ (and thus in $\Q[\x]$) and changes sign in $\R^2$, by \cite[Thm.4.5.1]{BCR} and \cite[Prop.3.2.4]{FG}, we have that $\II_\R(F_h)=(f_h)\R[\x]$ and $\II_\Q(F_h)=(f_h)\Q[\x]$, so $F_h\subset\R^2$ is defined over $\Q$. However, Fermat's Last Theorem implies that $F_h(\Q)=\varnothing$ so $F_h$ is not strongly defined over~$\Q$.

$(\mr{ii})$ Let $f:=\x_1-\sqrt[3]{2}\in\R[\x_1]$, let $g:=\x_1^2+\sqrt[3]{2}\x_1+\sqrt[3]{4}\in\R[\x_1]$ and let $h:=\x_1^3-2\in\Q[\x_1]$. Identify the polynomials $f$, $g$ and $h$ with the corresponding regular functions in $\reg(\R)$, so $fg=h$ on $\R$. Consider the singleton $X:=\{\sqrt[3]{2}\}\subset\R$, which is $\Q$-algebraic since $X=\ZZ_\R(h)$. Evidently, we have $\I^r_\R(X)=(f)\reg(\R)$ and $\I_\Q(X)=(h)\Q[\x_1]$. Since $\ZZ_\R(g)=\varnothing$, $g$ is invertible in $\reg(\R)$ and $\II^r_\R(X)=(f)\reg(\R)=(fg)\reg(\R)=(h)\reg(\R)=\II_\Q(X)\reg(\R)$. Thus, $X\subset\R$ is weakly defined over $\Q$. However, $X\subset\R$ is not defined over~$\Q$ {since} $f\in\I_\R(X)\setminus\I_\Q(X)\R[x_1]$.

$(\mr{iii})$ Let $f:=\x_1^2-\sqrt[3]{2}\x_2^3\in\R[\x]:=\R[\x_1,\x_2]$, let $g:=\x_1^4+\sqrt[3]{2}\x_1^2\x_2^3+\sqrt[3]{4}\x_2^6\in\R[\x]$ and let $h:=\x_1^6-2\x_2^9\in\Q[\x]$. As above, identify the polynomials $f$, $g$ and $h$ with the corresponding regular functions in $\reg(\R^2)$, so $fg=h$ on $\R^2$. Let $X\subset\R^2$ be the $\Q$-algebraic set $X:=\ZZ_\R(f)=\ZZ_\R(h)$. Using again \cite[Thm.4.5.1]{BCR} and \cite[Prop.3.2.4]{FG}, we deduce that $\II_\R(X)=(f)\R[\x]$ and $\II_\Q(X)=(h)\Q[\x]$. Since $\nabla f$ vanishes only at the origin $O$ of $\R^2$, it holds $\Reg(X)=X\setminus\{O\}$. Also $\nabla h$ vanishes only at $O$. Thus, if $a\in X\setminus\{O\}$, then $T_a(X)=T^\Q_a(X)$ {since} $\nabla h(a)=g(a)\nabla f(a)$ and $g(a)\neq0$. By Proposition \ref{Y}, we have $\Reg^\Q(X)=X\setminus\{O\}=\Reg(X)$, so $X\subset\R^2$ is $\Q$-determined. Since $\II_\Q(X)=(h)\Q[\x]$, we have $\II_\Q(X)\reg(\R^2)=(h)\reg(\R^2)$. Evidently, $f$ belongs to $\II^r_{\R^2}(X)$ {as} $f\in\II_\R(X)$ as a polynomial. However, $f$ does not belong to $\II_\Q(X)\reg(\R^2)=(h)\reg(\R^2)$. Otherwise, there would exist an element $\frac{p}{q}\in\reg(\R^n)$, where $p,q\in\R[\x]$ with $\ZZ_\R(q)=\varnothing$, such that $f=\frac{hp}{q}$ or, equivalently, $fq=hp=fgp$ on $\R^2$. Since $\R^2\setminus\ZZ_\R(f)$ is dense in $\R^2$, it would follow that $q=gp$ on $\R^2$, so $0\neq q(O)=g(O)p(O)=0$ {since} $g(O)=0$. This contradiction implies that $f\in\II^r_{\R^2}(X)\setminus\II_\Q(X)\reg(\R^2)$. Therefore, $X\subset\R^n$ is not weakly defined over $\Q$.

$(\mr{iv})$ Recall that the curves $C_1$ and $C_2$ of $\R^2$ defined in Remark \ref{rem1} are examples of $\Q$-algebraic sets that are not $\Q$-determined. $\sqbullet$ 
\end{examples}

In the next result we give a characterization of $\Q$-determined $\Q$-algebraic sets.

\begin{lem}\label{lem:Q-det&Q_def}
Let $X\subset\R^n$ be a $\Q$-algebraic set and let $U$ be an open subset of $\R^n$ such that $X\cap U=\Reg(X)$ (for example $U:=\R^n\setminus\Sing(X)$). The following conditions are~equivalent:  
\begin{itemize}
 \item[$(\mr{i})$] $X$ is $\Q$-determined.
 \item[$(\mr{ii})$] $\II^r_U(\Reg(X))=\II_\Q(X)\reg(U)$.
\end{itemize}
\end{lem}
\begin{proof}
In this proof we identify the polynomials in $\R[\x]$ with the corresponding functions defined on $\R^n$ or on $U$.

$(\mr{i})\Longrightarrow(\mr{ii})\ $ Suppose that $X$ is $\Q$-determined. By Proposition \ref{Y} and \cite[Prop.2.2.11]{akbking:tras}, for every $a\in\Reg(X)=\Reg^\Q(X)$, there exists a Zariski open neighborhood $U_a$ of $a$ in $\R^n$ such that $\II^r_{U_a}(\Reg(X)\cap U_a)=\II_\Q(X)\reg(U_a)$. The set $\Reg(X)$ endowed with the ($\R$-)Zariski topology is (Noetherian and thus) compact, so there exist finitely many points $a_1,\ldots,a_\ell$ of $\Reg(X)$ such that $\Reg(X)\subset\bigcup_{i=1}^\ell U_{a_i}$. Let $f_0\in\Q[\x]$ be such that $X=\ZZ_\R(f_0)$, let $f_1,\ldots,f_s$ be generators of $\II_\Q(X)$ in $\Q[\x]$ and, for every $i\in\{1,\ldots,\ell\}$, let $g_i\in\R[\x]$ be such that $\ZZ_\R(g_i)=\R^n\setminus U_{a_i}$. 

Pick any $f=\frac{p}{q}\in\II^r_U(\Reg(X))$, where $p,q\in\R[\x]$ and $U\cap\ZZ_\R(q)=\varnothing$. Since the restriction $p|_{U_{a_i}}$ belongs to $\II^r_{U_{a_i}}(\Reg(X)\cap U_{a_i})$, there exist $h_i,h_{i1},\ldots,h_{is}\in\R[\x]$ such that $\ZZ_\R(h_i)\subset\R^n\setminus U_{a_i}=\ZZ_\R(g_i)$ and  $p=\sum_{j=1}^s\frac{h_{ij}}{h_i}f_j$ on $U_{a_i}$, so $(g_ih_i)^2p=\sum_{j=1}^s(g_i^2h_ih_{ij})f_j$ on the whole~$\R^n$. We also have $f_0^2p=(f_0p)f_0$. Adding the latter equations member by member, we obtain $q'p=(f_0p)f_0+\sum_{j=1}^s(\sum_{i=1}^\ell g_i^2h_ih_{ij})f_j$ on $\R^n$ and thus on $U$, where $q':=f_0^2+\sum_{i=1}^\ell(g_ih_i)^2$. This proves that $q'p\in\II_\Q(X)\reg(U)$, so $q'f=(q'p)q^{-1}\in\II_\Q(X)\reg(U)$ as well. Observe that $U\cap\ZZ_\R(q')=U\cap\ZZ_\R(f_0)\cap\bigcap_{i=1}^\ell\ZZ_\R(g_ih_i)=\Reg(X)\setminus\bigcup_{i=1}^\ell U_{a_i}=\varnothing$, so $q'$ is invertible in $\reg(U)$ and thus $f=(q')^{-1}(q'f)\in\II_\Q(X)\reg(U)$.

$(\mr{ii})\Longrightarrow(\mr{i})\ $ Let $f_1,\ldots,f_s$ be generators of $\II_\Q(X)$ in $\Q[\x]$ and let $a\in\Reg(X)$. We have to show that $a\in\Reg^\Q(X)$. Pick any $f\in\II_\R(X)$ and denote $f^*\in\II^r_U(\Reg(X))$ the restriction of $f$ to $U$. By $(\mr{ii})$, there exist $h_1,\ldots,h_s\in\reg(U)$ such that $f^*=\sum_{i=1}^s h_if_i|_U$, so $\nabla f(a)=\nabla f^*(a)=\sum_{i=1}^sh_i(a)\nabla f_i(a)$. It follows that $T_a(X)=T^\Q_a(X)$, so $a\in\Reg^\Q(X)$ by Proposition \ref{Y}.
\end{proof}

As a consequence, we have:

\begin{cor}\label{cor:nonsing+Q-det=Q_def}
Let $X\subset\R^n$ be a nonsingular $\Q$-algebraic set. Then $X$ is weakly defined over~$\Q$ if and only if it is $\Q$-determined (or, equivalently, $\Q$-nonsingular).
\end{cor}
\begin{proof} 
By implication $(\mr{iii})\Longrightarrow(\mr{iv})$ of Proposition \ref{prop:implications}, if $X$ is weakly defined over $\Q$, then it is also $\Q$-determined. Suppose $X$ is $\Q$-determined. Since $X$ is nonsingular, we have $\Reg(X)=X$ by definition. Thus, we {may} apply implication $(\mr{i})\Longrightarrow(\mr{ii})$ of Lemma \ref{lem:Q-det&Q_def} with $U=\R^n$, obtaining that $X$ is weakly defined over~$\Q$, as required.
\end{proof}

\begin{ack}
The first author is supported by GNSAGA of INDAM, and partially supported by Spanish STRANO PID2021-122752NB-I00. The second author is supported by GNSAGA of INDAM, by the ANR NewMIRAGE (ANR-23-CE40-0002) and has been supported by IdEX of Universit\'e C\^ote d'Azur during the development of this work.
\end{ack}



\end{document}